\pdfoutput=1

\documentclass[a4paper, 11pt, abstract=true, titlepage, oneside, bibliography=totoc]{scrartcl}

\makeatletter
\DeclareOldFontCommand{\it}{\normalfont\itshape}{\emph}
\DeclareOldFontCommand{\bf}{\normalfont\bfseries}{\textbf}
\DeclareOldFontCommand{\sc}{\normalfont\scshape}{\textsc}
\makeatother

\makeatletter

\usepackage[T1]{fontenc}
\usepackage{lmodern}
\usepackage[utf8]{inputenc}
\usepackage[onehalfspacing]{setspace}
\usepackage[headsepline]{scrlayer-scrpage}
\clearpairofpagestyles
\usepackage{layout}
\usepackage[margin=1.1in]{geometry}
\usepackage{adjustbox}
\usepackage{tcolorbox}
\usepackage{authblk}
\usepackage[titletoc,title]{appendix}
\usepackage{comment}

\usepackage[hyphens]{url}
\usepackage{color}

\usepackage{amsmath,amssymb,amsfonts,amsthm, mathtools, bm}
\usepackage{thmtools, thm-restate}
\usepackage{nicefrac}
\usepackage{array}

\usepackage[normal]{threeparttable}
\usepackage{threeparttablex}
\usepackage{tabularx}
\usepackage{longtable}
\usepackage{booktabs}
\usepackage{colortbl}
\usepackage{multirow}
\usepackage{makecell}
\usepackage{tikz}
\usepackage{pbox}

\usepackage[acronym,hyperfirst=false,toc,nogroupskip,nomain,indexonlyfirst]{glossaries}

\usepackage{subcaption} 
\usepackage[font=footnotesize, labelfont=bf]{caption}
\usepackage{floatrow}
\usepackage{graphicx}

\usepackage{placeins}
\usepackage{float} 
\usepackage{caption}

\usepackage{algorithm}
\usepackage{algpseudocode}

\usepackage{tikz}
\usetikzlibrary{arrows.meta}
\usetikzlibrary{math}
\usetikzlibrary{shapes}
\usepackage{pgfplots}
\pgfplotsset{compat=1.18}
\usepgfplotslibrary{groupplots}
\usepgfplotslibrary{dateplot}
\usepackage{tikzscale}

\usepackage{enumerate}
\usepackage{paralist}
\usepackage{multicol}
\usepackage[author=Patrick]{pdfcomment}
\usepackage{xparse}
\usepackage{twoopt}
\usepackage{etoolbox}

\usepackage{eurosym}
\usepackage{ellipsis}
\usepackage{natbib}
\usepackage{paralist}
\usepackage{ulem}

\graphicspath{{Figures/}}
\usepackage{wrapfig}

\usepackage{optidef}

\usepackage{bm}
\usepackage{bbm}

\usepackage{siunitx} 
\usepackage{flafter}

\usepackage{eurosym}

\usepackage{floatrow}

\usepackage{transparent}

\usepackage{pgfplots}
\usepgfplotslibrary{statistics}
\pgfplotsset{
    compat=newest,
   /pgfplots/warning/plot without coordinates/.code={}, 
   filter discard warning=false
}
\usepackage{tikz}
\usepackage{pgf-pie}

\usepackage{ragged2e}
\newsavebox{\abstractbox}
\renewenvironment{abstract}
{\begin{lrbox}{0}\begin{minipage}{\textwidth}
			\begin{center}\normalfont\sectfont\abstractname\end{center}\quotation}
		{\endquotation\end{minipage}\end{lrbox}%
	\global\setbox\abstractbox=\box0 }

\makeatletter
\expandafter\patchcmd\csname\string\maketitle\endcsname
{\vskip\z@\@plus3fill}
{\vskip\z@\@plus2fill\box\abstractbox\vskip\z@\@plus1fill}
{}{}
\makeatother

\addtokomafont{captionlabel}{\footnotesize\bfseries} % style for caption labels
\addtokomafont{caption}{\footnotesize\bfseries} % and the caption
\bibpunct[, ]{(}{)}{,}{a}{}{,}%
\def\newblock{\ }%
\DeclareMathOperator*{\argmax}{arg\,max}

\counterwithin*{equation}{section}

\newenvironment{myfont}{\fontfamily{ccr}\selectfont}{\par}
\DeclareTextFontCommand{\textmyfont}{\myfont}

\AtBeginDocument{%
	\abovedisplayskip=7.5pt plus 0pt minus 0pt
	\abovedisplayshortskip=7.5pt plus 0pt minus 0pt
	\belowdisplayskip=7.5pt plus 0pt minus 0pt
	\belowdisplayshortskip=7.5pt plus 0pt minus 0pt
}
\newcolumntype{L}[1]{>{\raggedright\let\newline\\\arraybackslash\hspace{0pt}}p{#1}}
\newcolumntype{C}[1]{>{\centering\let\newline\\\arraybackslash\hspace{0pt}}p{#1}}
\newcolumntype{R}[1]{>{\raggedleft\let\newline\\\arraybackslash\hspace{0pt}}p{#1}}
\def\checkmark{\tikz\fill[scale=0.4](0,.35) -- (.25,0) -- (1,.7) -- (.25,.15) -- cycle;}

\renewcommand{\emph}[1]{\textit{#1}}

\newcounter{result}
\newenvironment{result}
  {
    \refstepcounter{result}
    \noindent \textbf{Result \theresult.} \itshape
  }
  {\par}

\setkomafont{section}{\normalfont\Large\bfseries}
\setkomafont{subsection}{\normalfont\large\bfseries}

\RedeclareSectionCommand[
  beforeskip=1.8em,
  afterskip=0.9em,
  indent=0pt,
  afterindent=false
]{section}

\RedeclareSectionCommand[
  beforeskip=1.2em,
  afterskip=0.4em,
  indent=0pt,
  afterindent=false
]{subsection}

\counterwithout{equation}{section}

\newglossarystyle{csuper}{%
  \renewenvironment{theglossary}%
    {\begin{tabular}{@{}l@{\hskip 0.5em}l@{\hskip 0.5em}>{\raggedleft\arraybackslash}p{\dimexpr\textwidth-6em\relax}@{}}}%
    {\end{tabular}}%

  \renewcommand*{\glossentry}[2]{%
    \glsentryitem{##1}\glstarget{##1}{\glsentryshort{##1}} & \dotfill & \Glossentrydesc{##1}##2\tabularnewline
  }%

}
\setglossarystyle{csuper}

\renewcommand*{\glossarysection}[2][]%
  {\par\addsec*{#2}}

\newacronym{acr:pt}{PT}{public transportation}
\newacronym{acr:pts}{PTS}{public transportation system}
\newacronym{acr:ft}{FT}{freight terminal}
\newacronym{acr:htu}{HTU}{hybrid transportation unit}
\newacronym{acr:lsp}{LSP}{logistic service provider}
\newacronym{acr:osm}{OSM}{{O}pen {S}treet {M}ap}
\newacronym{acr:mcnd}{MCND}{multicommodity fixed charge network design problem}
\newacronym{acr:lp}{LP}{linear problem}
\newacronym{acr:bap}{B\&P}{Branch-and-price}
\newacronym{acr:bd}{BD}{{B}enders {D}ecomposition}
\newacronym{acr:lr}{LR}{{Lagrangian {R}elaxation}}
\newacronym{acr:bapac}{B\&P\&C}{Branch-and-price-and-cut}
\newacronym{acr:alns}{ALNS}{{A}daptive {L}arge {N}eighborhood {S}earch}
\newacronym{acr:cg}{CG}{{C}olumn {G}eneration}
\newacronym{acr:bab}{B\&B}{Branch-and-bound}
\newacronym{acr:pab}{P\&B}{Price-and-branch}
\newacronym{acr:rmp}{RMP}{restricted master problem}
\newacronym{acr:vrp}{VRP}{vehicle routing problem}
\newacronym{acr:hlp}{HLP}{hub location problem}
\newacronym{acr:mip}{MIP}{mixed integer program}
\newacronym{acr:spp}{SPP}{shortest path problem}
\newacronym{acr:agv}{AGV}{automated guided vehicle}

\makeatletter
\newcommand{\customesmallfontsize}{%
  \@setfontsize\customsize{7pt}{7pt}%
}
\makeatother

\makeatletter
\newcommand{\anothersmallfontsite}{%
  \@setfontsize\customsize{8pt}{8pt}%
}
\makeatother

\makeatletter
\newcommand*{\DecimalMathComma}{%
  \mathchardef\dec@mcc=\mathcode`,
  \mathcode`,="613B 
  \mathchardef\dec@period=\mathcode`. 
  \mathcode`.="613A 
  \mathcode`.=\dec@period 
  \thinmuskip=0.5mu 
}
\makeatother

\newfloatcommand{capbtabbox}{table}[][\FBwidth]

\newcommand{\revone}[1]{{\leavevmode\color{black}#1}}
\newcommand{\revtwo}[1]{{\leavevmode\color{black}#1}}

\begin{document}
\emergencystretch 4em

% custom commands
\newcommand{\Graph}{G}
\newcommand{\SequenceOfRoute}{L}
\newcommand{\CapaFunction}{Q}
\newcommand{\Terminal}{b}
\newcommand{\Cost}{c}
\newcommand{\Destination}{d}
\newcommand{\RequestStart}{e}
\newcommand{\FreightFlow}{f}
\newcommand{\PassengerFlow}{g}
\newcommand{\Vehicle}{h}
\newcommand{\Layer}{k}
\newcommand{\RequestEnd}{l}
\newcommand{\Origin}{o}
\newcommand{\Path}{p}
\newcommand{\Demand}{q}
\newcommand{\Request}{r}
\newcommand{\Stop}{s}
\newcommand{\TimeStep}{t}
\newcommand{\Vertex}{v}
\newcommand{\DistanceApproximator}{w}
\newcommand{\HTUOps}{x}
\newcommand{\HTUVehicle}{y}
\newcommand{\FreightFlowPath}{z}
\newcommand{\SetOfArcs}{\mathcal{A}}  
\newcommand{\SetOfTerminals}{\mathcal{B}}
\newcommand{\SetOfDestinations}{\mathcal{D}}
\newcommand{\SpannedArcsFilter}{\mathcal{G}}
\newcommand{\SetOfVehicles}{\mathcal{H}}
\newcommand{\SpannedArcs}{\mathcal{I}}
\newcommand{\NotSpannedArcs}{\mathcal{J}}
\newcommand{\SetOfLayers}{\mathcal{K}}
\newcommand{\SetOfSequences}{\mathcal{L}}
\newcommand{\SetOfFlatStops}{\mathcal{M}}
\newcommand{\Neighbors}{\mathcal{N}}
\newcommand{\SetOfOrigins}{\mathcal{O}}
\newcommand{\SetOfPaths}{\mathcal{P}}
\newcommand{\SetOfFlatTerminals}{\mathcal{Q}}
\newcommand{\SetOfRequests}{\mathcal{R}}
\newcommand{\SetOfStops}{\mathcal{S}}
\newcommand{\SetOfTimeSteps}{\mathcal{T}}
\newcommand{\SetOfVertices}{\mathcal{V}}
\newcommand{\SetOfRoutes}{\mathcal{W}}
\newcommand{\DualAlpha}{\alpha}
\newcommand{\MapToFlat}{\beta}
\newcommand{\DualGamma}{\gamma}
\newcommand{\DualDelta}{\delta}
\newcommand{\DualEta}{\eta}
\newcommand{\PathArcRelation}{\theta}
\newcommand{\NumConnections}{\iota}
\newcommand{\NumUnits}{\kappa}
\newcommand{\UnitCapacity}{\lambda}
\newcommand{\SpanningMapping}{\mu}
\newcommand{\DualNu}{\nu}
\newcommand{\VertexDemand}{\xi}
\newcommand{\DualPi}{\pi}
\newcommand{\DualRho}{\rho}
\newcommand{\CostConversionFactor}{\sigma}
\newcommand{\DualTau}{\tau}
\newcommand{\DualUpsilon}{\upsilon}
\newcommand{\PricingStrength}{\phi}
\newcommand{\ServiceLevel}{\chi}
\newcommand{\UnitConversionFactor}{\psi}
\newcommand{\VehicleArcRelation}{\omega}
\newcommand{\Passenger}{\textsc{P}}
\newcommand{\Freight}{\textsc{F}}
\newcommand{\Arc}{(i,j)}
\newcommand{\UnbracedArc}{i,j}
\newcommand{\Penalty}{\textsc{Pen}}
\newcommand{\Routing}{\textsc{Routing}}
\newcommand{\Investment}{\textsc{Investment}}
\newcommand{\Access}{\textsc{A}}
\newcommand{\Egress}{\textsc{E}}
\newcommand{\Dummy}{\textsc{D}}
\newcommand{\Transit}{\textsc{T}}
\newcommand{\Hold}{0}
\newcommand{\VehicleLayer}{\textsc{V}}
\newcommand{\PathSegment}{\textsc{F}}
\newcommand{\Exp}{}
\newcommand{\Flat}{\textsc{F}}
\newcommand{\Simplified}{\textsc{S}}
\newcommand{\UB}{\textsc{UB}}
\newcommand{\LB}{\textsc{LB}}
\newcommand{\PotentialFreightArcs}{\textsc{C}}
\newcommand{\PassDemandScale}{\varrho}
\newcommand{\PenaltyCostScale}{\vartheta}

% then the title
\newgeometry{head=16.49675pt, top=0.5in}

\title{\large Dynamic capacity allocation of hybrid transportation units for cargo-hitching in urban public transportation systems}

\author[a]{\normalsize Paul Bischoff\thanks{Corresponding author}}
\author[a]{\normalsize Benedikt Lienkamp}
\author[b]{\normalsize Tarun Rambha}
\author[a,c]{\normalsize Maximilian Schiffer}

\affil[a]{
\small School of Management, Technical University of Munich, Germany

\scriptsize paul.bischoff@tum.de, benedikt.lienkamp@tum.de
}
\affil[b]{
\small Centre for Infrastructure, Sustainable Transportation and Urban Planning (CiSTUP) \& Robert Bosch Centre for Cyber-Physical Systems (RBCCPS), Indian Institute of Science, Bengaluru, India

\scriptsize tarunrambha@iisc.ac.in
}
\affil[c]{
\small Munich Data Science Institute, Technical University of Munich, Germany

\scriptsize schiffer@tum.de
}

\date{}

\lehead{\pagemark}
\rohead{\pagemark}

\begin{abstract}
\begin{singlespace}
{\small\noindent 
To improve the utilization of public transportation systems (PTSs) during off-peak hours, we present an algorithmic framework that designs PTSs with hybrid transportation units (HTUs), which can transport passengers or freight by leveraging a flexible interior. 
Against this background, we study a capacitated network design problem to enable cargo-hitching in existing PTSs. Specifically, we study a setting with fixed vehicle routes and timetables in which vehicles can be equipped with HTUs to enable cargo-hitching.  
We optimize the network design from a total cost perspective to account for normalized network design costs tied to the investment in HTUs and freight routing costs. We present an algorithmic framework that encodes some of the problem's constraints in a spatially and temporally expanded, layered graph, and solves the resulting network design problem with a price-and-branch algorithm. 
We apply this framework to a case study based on the subway network in the city of Munich. Our algorithm outscales commercial solvers by a factor of six and yields integer feasible solutions with a median integrality gap of less than $1.56$\% for all instances. We show that cargo-hitching with HTUs increases the utilization of PTSs, especially during off-peak hours, without cannibalizing passenger service level and quality. \revone{We quantify the value of \glspl{acr:htu} at up to $3.2$\% of the total cost}. Moreover, we present a sensitivity analysis that indicates that cargo-hitching is worthwhile if truck-based transport occurs at an externality cost of more than \euro{1.5} per vehicle and kilometer and loading and unloading costs of less than \euro{2.0} per passenger equivalent.\\
\smallskip}
\end{singlespace}
\end{abstract}

\maketitle
\restoregeometry

\section{Introduction} \label{sec:introduction}
\noindent As urban populations grow and cities become more interconnected, the demand for efficient public transit rises \citep{UN2018}. Moreover, this population increase leads to significant growth of e-commerce transactions whose transportation contributes up to $15\%$ of urban road transport \citep{Dablanc2011}.
As a consequence, cities suffer from overloaded transportation systems, whose negative externalities cause environmental harm via CO\textsubscript{2} emissions, health dangers via particulate matters and NO\textsubscript{x}, and economic harm through working hours lost in congestion \citep{LevyBuonocoreEtAl2010, FattahMorshedEtAl2022}. Focusing on freight transportation, electric vehicles and city freighters allow to reduce emission-related externalities. Focusing on passenger transportation, \gls{acr:pt} yields low costs per trip. Additionally, \glspl{acr:pts} offer sustainable mobility solutions as the emissions per passenger in a highly utilized \gls{acr:pts} are significantly lower compared to individual mobility solutions \citep{NoussanEtAl2022}. So far, concepts discussed to realize sustainable transportation often focus either exclusively on freight or passenger transportation but share a central characteristic: the sustainability of each concept increases with greater utilization of its transportation modes. Still, for both freight and passenger transportation, concepts that allow to permanently maintain a high utilization are missing. \Glspl{acr:pts} in European cities show off-peak utilizations below 40\% between 10 a.m. and 4 p.m. \citep{ChengGuo2018,McKinsey2020} and freight transport by design contains dead-headed driving, particularly when trucks or city freighters return to a depot. 

To mitigate low \gls{acr:pts} utilization during off-peak hours and relief heavily congested road networks partially occupied by freight trucks, this paper studies the concept of \textit{cargo-hitching}, where a municipality equips its \gls{acr:pts} such that it accommodates intermodal freight transportation without cannibalizing its primary purpose of offering convenient passenger transportation services. The concept promises a utilization increase of the \gls{acr:pts} at zero additional installed capacity by using spare capacity available predominantly during off-peak hours. Furthermore, the utilization increase comes hand-in-hand with relieving congestion on roads because conventional truck-based deliveries can be reduced, which underlines the concept's win-win nature. 

\begin{figure}[!bt]
\begin{floatrow}
\ffigbox{
    \scriptsize
    \def\svgwidth{\columnwidth} 
    \includegraphics[width=0.9\linewidth]{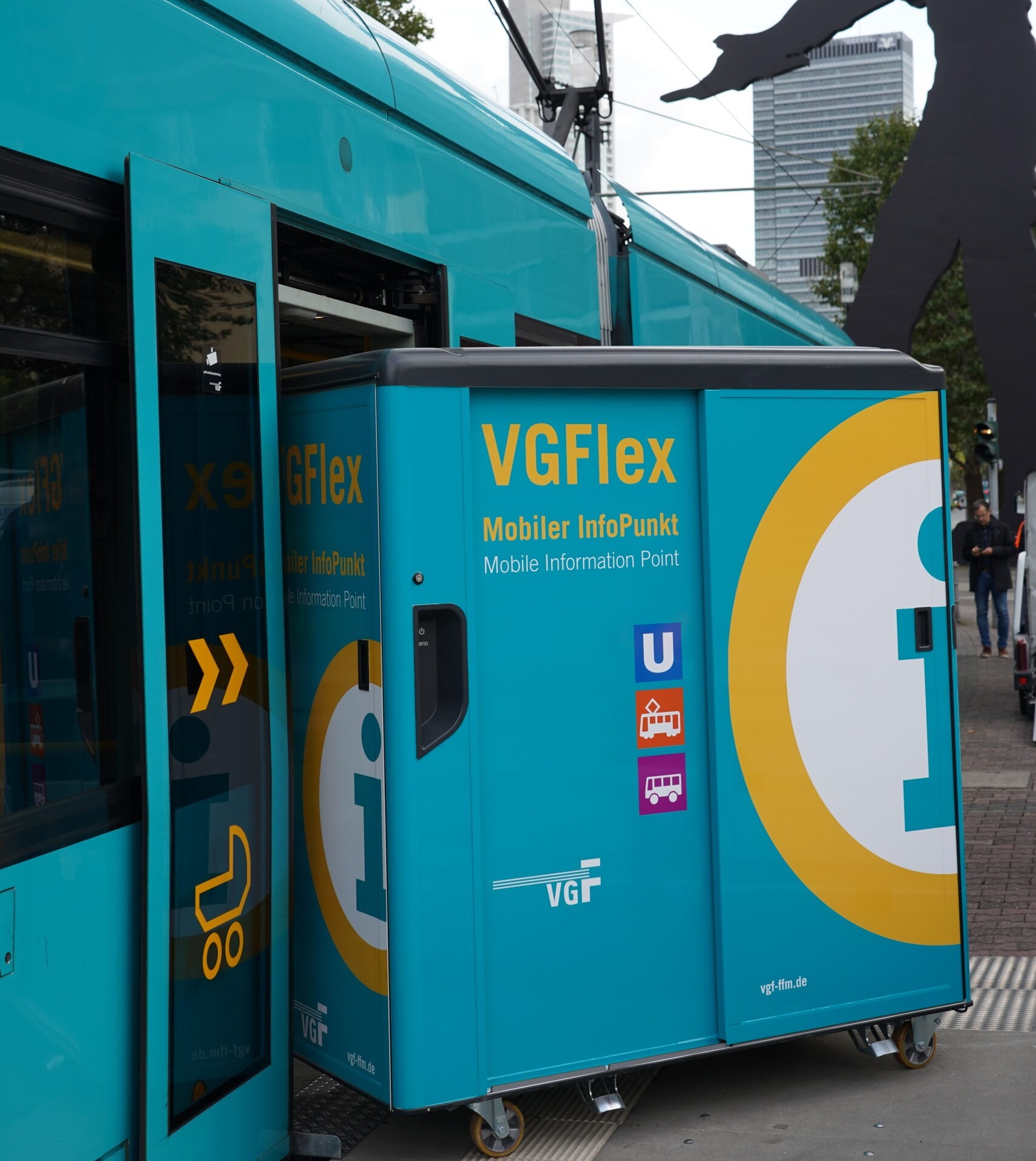}
}{
    \caption{Practical implementation of cargo-hitching \citep{Onomotion2021}}
    \label{fig:dhl-paketbahn}
}
\ffigbox{
    \scriptsize
    \begin{subfigure}[b]{\linewidth}
        \centering
        \def\svgwidth{\columnwidth}
        \includegraphics[width=0.9\linewidth]{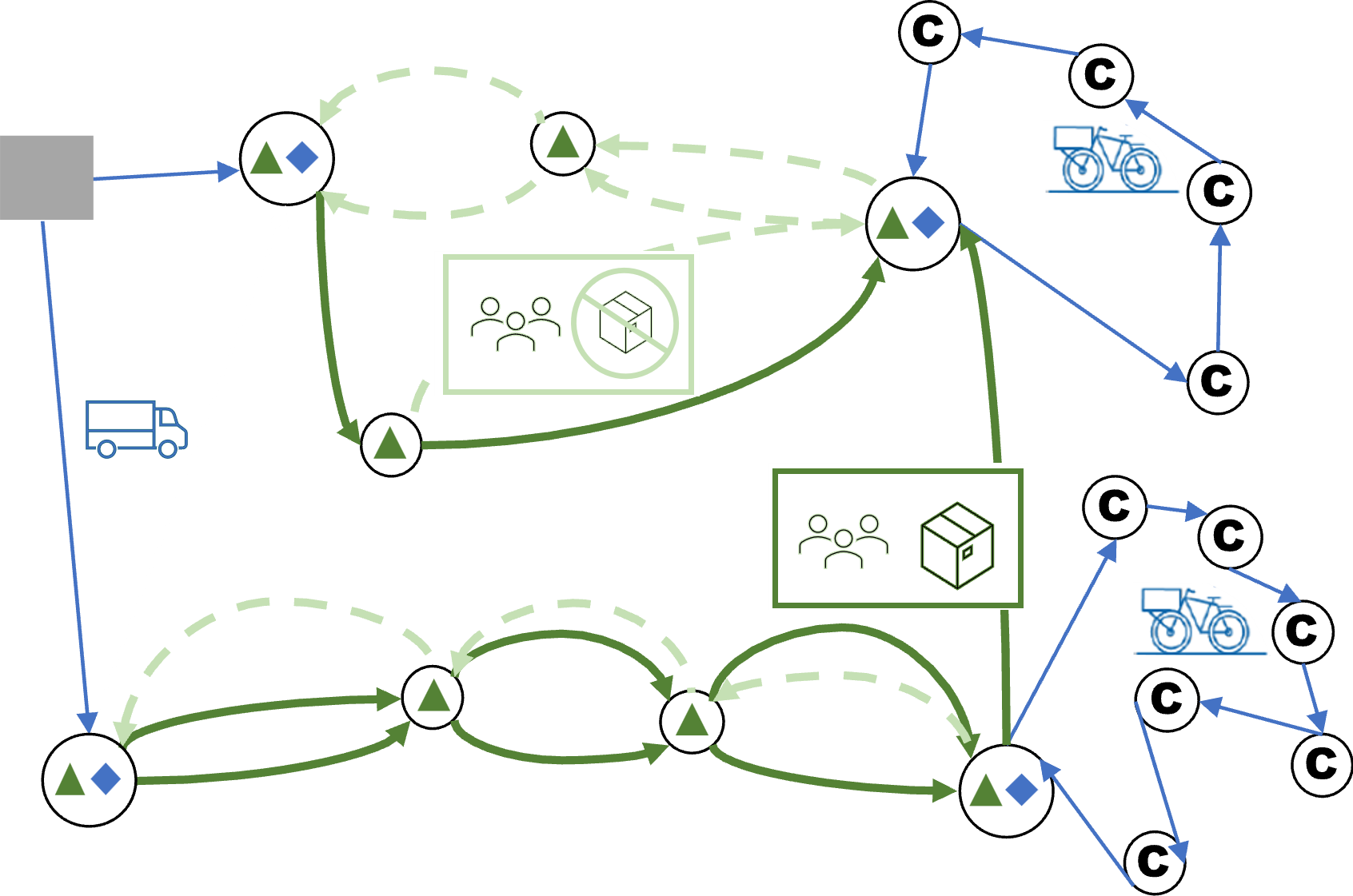}
        \caption*{} 
        \label{fig:htu-concept}
    \end{subfigure}
    
    \vspace{-0.4cm} 
    
    \begin{subfigure}[b]{\linewidth}
        \centering
        \def\svgwidth{0.9\columnwidth}
        \input{Figure_2_Legend.tex}
        \caption*{} 
    \end{subfigure}
}{
    \caption{Cargo-hitching system as three-echelon freight delivery via \glspl{acr:pts}}
    \label{fig:system}
}
\end{floatrow}
\end{figure}

Cargo-hitching has attracted practitioners' attention over the last two decades. Although the first notable implementation \textit{City Cargo} in Amsterdam was stopped primarily due to financing issues during the 2008 economic crisis \citep{ArvidssonBrowne2013}, large urban \gls{acr:pts} operators believe in the potential of the concept and fund its development and implementation \citep{VGF2021}. For example, Figure~\ref{fig:dhl-paketbahn} shows a recent and ongoing project, the \textit{Gütertram} in Frankfurt am Main, Germany. 

Figure~2 schematically shows a system in which a \gls{acr:pts} operator can equip selected \gls{acr:pt} vehicles to accommodate both passenger and freight transportation, and can select \gls{acr:pt} stops to be used to exchange freight between \glspl{acr:lsp}' vehicles and the equipped \gls{acr:pt} vehicles. As a result, freight deliveries pass a three-echelon system consisting of truck delivery to \gls{acr:pt} stops performed by \glspl{acr:lsp}, transportation in \revtwo{a selected subset of} the \gls{acr:pt} vehicles performed by the transit system operator, and last-mile delivery via city freighters again performed by \glspl{acr:lsp}. More specifically, \glspl{acr:lsp} transport truckloads of freight to selected \gls{acr:pt} stops that are equipped to handle freight and store freight shipments temporarily. The \gls{acr:pts} operator receives the freight, and is responsible for their timely transportation via the \gls{acr:pt} network to the freight shipment's broader destination area. The operator transports the freight shipments in the scheduled \gls{acr:pt} vehicles \revtwo{that are equipped to transport freight} by utilizing available excess capacity. Finally, after passing the \gls{acr:pt} network, the freight is available for transportation in the broader destination area: \glspl{acr:lsp} pick the freight shipments up at the respective \gls{acr:pt} stops, and subsequently perform the last-mile delivery via city freighters.

\noindent\textbf{Contribution:} We propose a novel urban cargo-hitching design problem that is, to the best of our knowledge, the first work to consider the dynamic allocation of \gls{acr:pts} capacity in a multi-line \glspl{acr:pts}. To solve real-world instances, we present a \gls{acr:pab} framework that relies on the following problem specific modifications to ensure scalability: first, it leverages a problem specific partially temporal and spatial graph expansion that allows to encode the complexity of some problem constraints in its graph structure. It further relies on an efficient pricing scheme that allows to decompose pricing problems by requests into \glspl{acr:spp} that can be solved via the $A^{*}$ algorithm. 
Moreover, we introduce multiple preprocessing techniques involving graph pruning and arc contraction that effectively decrease the cardinality of the arc set. To benchmark the performance of this algorithmic framework, we further provide a \gls{acr:mip} formulation as well as a \gls{acr:bap} algorithm. Our framework utilizing the \gls{acr:pab} approach solves instances with up to $3{,}000$ requests on the subway network of Munich to a median integrality gap of $1.56\%$ within computation times of $90$ minutes. The proposed \gls{acr:bap} algorithm improves the obtained median integrality gaps by up to $0.38$ percentage points at the price of significantly increased \revone{pricing effort}. Both of our algorithmic solutions outscale a commercial solver by more than a factor of 6. 

Beyond providing evidence on the computational efficiency of our framework, we derive several managerial insights based on our case study for the city of Munich. First, we show that cargo-hitching can offer a utilization increase at zero additional installed capacity. In this context, our algorithmic framework provides solutions that predominantly utilize the \gls{acr:pts}’s off-peak hours to transport freight requests. Second, we show that the amount of freight transported via cargo-hitching is sensitive to the externality cost for truck-based delivery and the cost for (un-)loading freight into the \gls{acr:pts}. In this context, our sensitivity analysis indicates that the full potential of cargo-hitching in the Munich subway network is realized if truck-based transport has externality costs of more than \euro{1.5} per vehicle and kilometer and loading and unloading costs are less than \euro{2.0} per passenger equivalent. Finally, we show that the savings realized by cargo-hitching depend on the spare capacity left within the \gls{acr:pts} as well as the amount of freight requests that can be shipped through it. In this context, relying on hybrid transportation units (HTUs) to realize cargo-hitching is particularly beneficial if the amount of freight requests is high but the spare capacity left in the \gls{acr:pts} fluctuates over the day due to passenger flow peaks. \revone{Moreover, allowing the dynamic allocation of capacity supports the acceptance of cargo-hitching and saves up to $3.2$\% of the total cost compared to static allocations.}

\noindent \textbf{Organization:} The remainder of this paper \revone{--— Appendix~\ref{app:abbrevations} lists the abbreviations that we use throughout ---} is as follows. In Section~\ref{sec:literature}, we summarize the state of the art. We then elaborate on our problem setting in Section~\ref{sec:problem-description} and develop our algorithmic framework in Section~\ref{sec:algorithmic-framework}. Section~\ref{sec:case-study} describes our case study based on the subway network in Munich, Germany. We present numerical results that show the efficiency of our algorithmic framework and derive managerial insights in Section~\ref{subsection:comp-results}. We conclude in Section~\ref{sec:conclusion} by summarizing our main findings.

\section{State of the art} \label{sec:literature}
\noindent Exploring the potential of cargo-hitching relates back to the seminal work of \mbox{\cite{TrentiniMalhene2012}}.
Early publications studied problem settings from an \gls{acr:lsp} perspective. Some works focussed on routing freight through a \gls{acr:pts} without considering inter-dependencies to the preceding and succeeding vehicle routing tasks. In particular, \cite{FatnassiChaouachiEtAl2015, BehiriBelmokhatar-BerrafEtAl2018}, and \mbox{\cite{OzturkPatrick2018}} studied the scheduling of vehicles for freight transportation on a given network with fixed routes. \mbox{\cite{ChengGuo2018}} studied the matching of freight to services. \mbox{\cite{MachadoDeSousa2023}} extended the assignment problem to a matching between requests on the one hand and stations and services on the other hand. Furthermore, \cite{MaZhangEtAl2023} considered a single-line co-modal urban \gls{acr:pts} and provided operative time-invariant equilibrium conditions, such as fare prices and capacities.  
In contrast, \cite{GhilasDemirEtAl2016} introduced the pickup and delivery problem with scheduled lines adapting the \gls{acr:lsp} perspective. They extended their work in further publications by providing an exact solution approach based on \gls{acr:bap} \mbox{\citep{GhilasCordeauEtAl2018}} as well as a heuristic approach based on \gls{acr:alns} \mbox{\citep{GhilasDemirEtAl2016b}}, and considered stochastic demands \citep{GhilasDemirEtAl2016c}. Other works studied related operational problem settings. Specifically, \cite{MassonTrentiniEtAl2017} solved a two-tier \gls{acr:vrp} via ALNS in which they considered the \gls{acr:pts} and subsequently the last-mile delivery via city freighters, \cite{MandalArchetti2023} studied a three-tier \gls{acr:vrp} in which they additionally considered the transportation to the \gls{acr:pts} and applied a decomposition method to solve it. We refer to \mbox{\cite{MouradPuchingerEtAl2019}} and \mbox{\cite{ElbertRentschler2022}} for more details and recent advances on the operational aspects of cargo-hitching and related city logistic concepts. Instead, we focus our discussion on the strategic aspects of node-based network design decisions, and the tactical and operational aspects of arc-based mode choices in the following. 

\noindent \textbf{Network design studies:} \cite{ZhaoLiEtAl2018} and \cite{JiZhengEtAl2020} formulated \glspl{acr:hlp} in order to determine suitable \gls{acr:pt} stops to handle freight and demonstrated their approaches on the Shanghai network. However, both works neglected capacity restrictions at hubs.  
In contrast, \cite{ElOuadiErroussoEtAl2022} assign customers to suburban or urban bundling hubs with restricted hub capacities. However, the \gls{acr:pt} lines' capacities are unlimited, and flows are only considered at the hubs but not in the \gls{acr:pts}. They applied machine learning to cluster zones and predict demands. \cite{AzcuyAgatzEtAl2021} studied a two-tier delivery system with a given \gls{acr:pt} capacity allocation, and minimized the expected travel distance performed by the last-mile vehicles. Thus, they derived strategic insights on \gls{acr:pt} stop locations from the solution of the operative \gls{acr:vrp} problem. Similarly, \mbox{\cite{DonneAlfandariEtAl2023}} determined suitable \gls{acr:pt} stops and \gls{acr:pt} lines, but still considered a simplified problem setting as the selected stops determine which \gls{acr:pt} lines can transport freight, i.e., no explicit capacity allocation decision happens at the vehicle level. Moreover, the problem setting remains time-invariant and ignores important transshipment and synchronization constraints. \cite{Nieto-IsazaFontaineEtAl2022} and \cite{KildizYildiz2023} studied crowd-shipping delivery systems and strategically determined the locations of mini depots and satellites based on two-stage stochastic network design problems with stochastic demands. Although focusing on crowd-shipping, their problem setting shows parallels to our problem as they determine freight routes based on given network layouts and time-tables. However, they ignore the system operator's option to allocate capacities, and assume that \gls{acr:pts} capacity is sufficient. Furthermore, \cite{ElbertRentschlerEtAl2023} studied a combined hub location and service network design problem for long-haul rail transportation but discarded passenger flow-related constraints and assumed constant capacities. 

\noindent \textbf{Mode choice studies: } Some works investigated the system operator's mode choices in co-modal \glspl{acr:pts} but did not allow for the dynamic re-allocation of capacity. In this context, some works determined the sharing mode in cargo-hitching systems \mbox{\citep{DiYangEtAl2022}}, studied the scheduling of freight vehicles on fixed networks \mbox{\citep{HoerstingCleophas2023}}, or assigned freight to fixed services \mbox{\citep{MachadoPimentel2023}}. 
\mbox{\cite{DiYangEtAl2022}} studied the joint optimization of train carriage arrangement and flow control. They determine the capacity allocation in terms of number of freight units to attach to every \gls{acr:pt} vehicle specifically during off-peak hours and, consequently, ignore the dynamic vehicle capacity re-allocation during operations. Moreover, \cite{MachadoPimentel2023} assumed that the demand for passenger transportation is known apriori and studied a stochastic problem in which uncertain freight demands are dropped into the left-over capacity of the \gls{acr:pt} bus system. 
Mostly, existing works consider two sharing modes: a \textit{sharing-train mode} that allows passengers and freight to be transported in separate units of the same vehicle, and a \textit{sharing-carriage mode} that allows passengers and shipments to share the same unit. \mbox{\cite{HoerstingCleophas2023}} considered the sharing mode as an exogenous feature, and restricted the constant capacity accordingly. \cite{LiZhuEtAl2023} determined the mode of each vehicle and the freight routes through the \gls{acr:pts}. Although similar to our problem setting with respect to the \textit{sharing-train mode}, they allowed split freight routes which requires expensive unloading, sorting, and loading operations, and limits the works applicability in real-world settings. Moreover, their problem setting relies on a penalty to prevent the reduction of passenger level of service. Hence, they allow unlimited passenger service cannibalization if the benefits outweigh the penalties. \cite{LiZhuEtAl2024} and \cite{LinZhang2024} neglect all strategic design aspects and investigate settings closest to our problem setting. \cite{LiZhuEtAl2024} study the scheduling of passenger and freight underground units with semi-dynamic allocation of capacity. In their setting, a system operator can change the composition of underground trains in between trips at terminal stations. However, they do not enforce the prioritization of passengers over freight. Moreover, \cite{LinZhang2024} determine the number of transportation units in an underground system, their assignment to vehicles and the allocation of their capacity. Their work remains limited to a single line \gls{acr:pt}.  

\noindent \textbf{Conclusion: } Table~\ref{tbl:literature-table} shows the characteristics of the closest related works on network design and mode choices for cargo-hitching. As can be seen, no work exists that considers dynamic capacity allocation in a shared-vehicle setting when determining the network design of an urban cargo-hitching system. Notably, \mbox{\cite{DiYangEtAl2022}} and \cite{LiZhuEtAl2023} considered a constant capacity allocation task by determining the sharing mode on the \gls{acr:pt} vehicle level but discarded other important characteristics. \cite{LiZhuEtAl2024} and \cite{LinZhang2024} increased the flexibility of the capacity allocation task but discard all design aspects. Other works about cargo-hitching discarded the capacity allocation task even on a tactical level where capacity allocation is constant.  
\def\checkmark{\tikz\fill[scale=0.4](0,.35) -- (.25,0) -- (1,.7) -- (.25,.15) -- cycle;} 
\newcommand{\rt}[1]{\rotatebox{90}{#1}}
\newcommand{\lb}[1]{\pbox{12mm}{#1}}
\newcolumntype{Y}[1]{>{\centering\arraybackslash}p{#1}}
\begin{table}[!b]
\scriptsize
\begin{tabular}{l*{12}{>{\centering\arraybackslash}m{0.38cm}}>{\centering\arraybackslash}m{0.5cm}}
\hline
\textbf{} & \rotatebox{60}{\makebox[3.8cm][l]{\raggedright\cite{BehiriBelmokhatar-BerrafEtAl2018}}} & \rotatebox{60}{\makebox[3.8cm][l]{\raggedright\cite{JiZhengEtAl2020}}} & \rotatebox{60}{\makebox[3.8cm][l]{\raggedright\cite{AzcuyAgatzEtAl2021}}} & \rotatebox{60}{\makebox[3.8cm][l]{\raggedright\cite{DiYangEtAl2022}}} & 
 \rotatebox{60}{\makebox[3.8cm][l]{\raggedright\cite{DonneAlfandariEtAl2023}}} & \rotatebox{60}{\makebox[3.8cm][l]{\raggedright\cite{Nieto-IsazaFontaineEtAl2022}}} & \rotatebox{60}{\makebox[3.8cm][l]{\raggedright\cite{MachadoPimentel2023}}} & \rotatebox{60}{\makebox[3.8cm][l]{\raggedright\cite{ElbertRentschlerEtAl2023}}} & \rotatebox{60}{\makebox[3.8cm][l]{\raggedright\cite{LiZhuEtAl2023}}} & \rotatebox{60}{\makebox[3.8cm][l]{\raggedright\cite{HoerstingCleophas2023}}} & \rotatebox{60}{\makebox[3.8cm][l]{\raggedright\cite{LiZhuEtAl2024}}}  & \rotatebox{60}{\makebox[3.8cm][l]{\raggedright\cite{LinZhang2024}}} & \rotatebox{60}{\makebox[3.8cm][l]{\raggedright\textbf{Our work}}} \\ \hline
Heterogeneous \gls{acr:pt} vehicles & $\checkmark$ & - & - & - & $\checkmark$ & - & $\checkmark$ & - & $\checkmark$ & - & - & - & $\checkmark$\\
Freight transshipments & - & - & - & - & - & $\checkmark$ & - & - & $\checkmark$ & - & $\checkmark$  & - & $\checkmark$\\
Time synchronization & $\checkmark$ & - & $\checkmark$ & $\checkmark$ & - & $\checkmark$ & $\checkmark$ & $\checkmark$ & $\checkmark$ & $\checkmark$ & $\checkmark$ & $\checkmark$ & $\checkmark$\\
Limited capacities & $\checkmark$ & - & $\checkmark$ & $\checkmark$ & $\checkmark$ & $\checkmark$ & $\checkmark$ & $\checkmark$ & $\checkmark$ & $\checkmark$ & $\checkmark$ & $\checkmark$ & $\checkmark$\\
Optional cargo-hitching & - & - & $\checkmark$ & $\checkmark$ & - & $\checkmark$ & - & $\checkmark$ & - & $\checkmark$ & $\checkmark$ & - &  $\checkmark$ \\
Passenger service level & $\checkmark$ & - & - & - & $\checkmark$ & $\checkmark$ & $\checkmark$ & - & - & $\checkmark$ & - & $\checkmark$ & $\checkmark$\\
Capacity allocation & - & - & - & $\checkmark$ & - & - & - & - & $\checkmark$ & - & $\checkmark$ & $\checkmark$ & $\checkmark$\\ \hline 
\end{tabular}
\caption{Related works on quantitative strategic or tactical research on cargo-hitching.} \label{tbl:literature-table}
\end{table}

\section{Problem setting} \label{sec:problem-description}
\noindent This work develops an algorithmic framework for the strategic planning tasks of a municipality to enable cargo-hitching in their \gls{acr:pts}. To do so, the municipality needs to transform their \gls{acr:pts} by adding two novel elements: 

\noindent \textbf{\Glspl{acr:ft}} are designated \gls{acr:pt} stops that allow the exchange of freight between a \gls{acr:pts} and other means of transportation. At an \gls{acr:ft}, a truck can unload freight, which is then transported on a leg of the \gls{acr:pts}, and unloaded at a different \gls{acr:ft} for last-mile delivery. Additionally, \glspl{acr:ft} allow storing and transshipping freight deliveries. The operations in \glspl{acr:ft} are automated with \glspl{acr:agv} and are spatially separated from the passengers.  

\noindent \textbf{\Glspl{acr:htu}} allow transporting both freight and passengers. We focus on \glspl{acr:htu} with a flexible interior that can be changed between trips to accommodate freight or passengers --- but not both at the same time. For example, an \gls{acr:htu} can be a specifically designed subway train wagon \citep[cf.][]{KellyMarinov2017}.

By replacing \revone{a subset of the} conventional wagons with \glspl{acr:htu}, a municipality decides on the share of a \gls{acr:pt} vehicle that can be flexibly used for either freight or passenger transport. We note that this setting describes a special case of the \textit{shared-vehicle} approach, which is particularly amenable for rail-based systems in which passengers and freight share the vehicle but not the wagon \citep{ElbertRentschler2022}. \revone{The replacement of some conventional wagons with HTUs is a strategic planning decision that determines which modular parts of the fleet are equipped with flexible interiors. This decision determines the fleet’s potential flexibility but does not dictate its operational use. Accordingly, the system operator determines the operational mode — transporting either freight or passengers — of these \glspl{acr:htu} in their short-term planning because this decision depends on constraints such as the current demand.} Switching an \gls{acr:htu} from passenger to freight transport or vice versa, e.g., by unfolding seats, requires little set-up time and can be automated. We restrict such mode changes to happen only at \glspl{acr:ft} due to practical constraints but allow for multiple functional switches of the \glspl{acr:htu} during a day, such that the municipality can vary the share of freight capacities in the \gls{acr:pts} between peak and off-peak hours to account for passenger demand.

We solve the municipalities' strategic planning problem of deciding which \gls{acr:pt} lines and vehicles to equip with how many \glspl{acr:htu} to minimize total freight transportation cost by allowing for cargo-hitching. 
We do not address the selection of suitable \glspl{acr:ft} but instead assume that the appropriate subset of \gls{acr:pt} stops has already been determined. The selected stops should offer sufficient space for freight operations and must be built at stops where \gls{acr:pt} schedules have a slightly longer stopover time.
We explicitly consider passenger and freight transportation but prioritize passenger transportation to reflect concerns about limited acceptance of the concept that might arise if the current passenger service levels cannot be maintained. The resulting planning problem resembles a capacitated multi-commodity fixed-charge network design problem with an additional layer of complexity. The additional complexity arises from the capacity allocation decision between passengers and freight that leads to an additional integral dependency in which an \gls{acr:htu} can accommodate freight or passengers but not both simultaneously. 
The operator determines the assignment of \glspl{acr:htu} to \gls{acr:pt} vehicles and the dynamic allocation of flexible capacity from assigned \glspl{acr:htu} to either passenger or freight. Additionally, the operator decides on the subset of freight requests that are accepted for transportation via the \gls{acr:pts}, the subset of freight requests that are rejected, and the paths on which the passenger requests and the accepted freight requests are routed through the \gls{acr:pts}. 

We formally define the resulting planning problem in Section~\ref{subsection:problem-definition}. Then, Section~\ref{subsection:graph-construction} formalizes the construction of our expanded graph representation that allows to encode certain problem characteristics, and Section~\ref{subsection:milp} provides a \gls{acr:mip} formulation of our problem.

\subsection{Problem definition} \label{subsection:problem-definition}
\noindent Formally, we consider a set of requests $\SetOfRequests = \SetOfRequests^{\Passenger} \cup \SetOfRequests^{\Freight}$, which is the union of two distinct subsets: passenger requests $\SetOfRequests^{\Passenger}$ and freight requests $\SetOfRequests^{\Freight}$. Every request is defined as a quintuple $\Request= (\Origin^{r}, \Destination^{\Request}, \Demand^{\Request}, \RequestStart^{\Request}, \RequestEnd^{\Request})$. Here, $\Origin^{r}$ denotes the request's origin, $\Destination^{\Request}$ its destination, $\Demand^{\Request}$ its demand, and $\RequestStart^{\Request}$ as well as $\RequestEnd^{\Request}$ define the time interval $[e^{\Request}, l^{\Request}]$ in which the request must be processed, with $\RequestStart^{\Request}$ being the earliest start time and $\RequestEnd^{\Request}$ marking the latest service completion time. The constraints on the service time apply to both passenger and freight requests. By adjusting a time window in one direction --- specifically, by relaxing either $\RequestStart^{\Request}$ or $\RequestEnd^{\Request}$ --- we can adapt the less restrictive assumptions commonly used in works on flow assignment in public transit where either only a departure time or only an arrival time is given. 

\revone{The \gls{acr:pts} consists of a set of stops $\Stop \in \SetOfFlatStops$. A subset of stops $\SetOfFlatTerminals \subseteq \SetOfFlatStops$ serves as \glspl{acr:ft}. Moreover, a fleet of transit vehicles, denoted by $\SetOfVehicles$, operates on this network with every transit vehicle $\Vehicle \in \SetOfVehicles$ following a specific route, defined as an ordered sequence of stop-time tuples $\SequenceOfRoute_{\Vehicle} = \langle(\Stop_1,\TimeStep_1),\ldots,(\Stop_{n},\TimeStep_{n})\rangle$, where $n$ represents the number of stops on that route. Let $\SetOfSequences := \{\SequenceOfRoute_{\Vehicle}: \Vehicle \in \SetOfVehicles \}$ denote the set of routes. Moreover, we denote the set of all scheduled stop-time tuples by $\SetOfRoutes := \bigcup_{\SequenceOfRoute \in \SetOfSequences} \{ (\Stop, \TimeStep): (\Stop, \TimeStep) \in \SequenceOfRoute \}$. Every vehicle $\Vehicle$ consist of $\NumUnits_{\Vehicle}$ identical units with a unit capacity of $\UnitCapacity_{\Vehicle}$. For subsequent discussions, we define $\SetOfTimeSteps(\Stop):= \langle \TimeStep^{\Stop}_{1}, \TimeStep^{\Stop}_{2}, ... \rangle$ as the sorted list of arrival times of all vehicles at an arbitrary stop $\Stop \in \SetOfFlatStops$. In this context, $\TimeStep^{\Stop}_{u} \in \{ \TimeStep : (\Stop^{\prime}, \TimeStep) \in \SetOfRoutes,\ \Stop^{\prime} = \Stop \}$ and $\TimeStep^{\Stop}_{u} < \TimeStep^{\Stop}_{u+1}$ for all $u$.}
 
\noindent\textbf{Decisions:} 
The municipality services a given request $\Request \in \SetOfRequests$ by sending it from its origin $\Origin^{\Request}$ through the \gls{acr:pts} to its destination $\Destination^{\Request}$. Here, any possible connection from $\Origin^{\Request}$ to $\Destination^{\Request}$ is called a path $\Path \in \SetOfPaths(\Request)$ where $\SetOfPaths(\Request)$ denotes the set of all feasible paths servicing~$\Request$. \revone{Every path can be represented as a sequence $\Path = \langle (\Origin^{\Request}, \RequestStart^{\Request}), (\Stop^{\Path}_{1}, \TimeStep^{\Path}_{1}),\ldots,(\Stop^{\Path}_{|\Path|-2}, \TimeStep^{\Path}_{|\Path|-2}), (\Destination^{\Request}, \RequestEnd^{\Request})\rangle$ where $|\Path|$ denotes the length of path $\Path$.} The municipality takes the following decisions:
\begin{compactitem}
    \item[i)] assigning an integer number $\HTUVehicle_{\Vehicle} \le \NumUnits_{\Vehicle}, \HTUVehicle_{\Vehicle} \in \mathbb{N}_0, \Vehicle \in \SetOfVehicles$ of \glspl{acr:htu} to every \gls{acr:pt} vehicle. The remaining units can only transport passengers.
    \item[ii)] allocating an integer number $\HTUOps_{(\Stop_i,\TimeStep_i),(\Stop_j, \TimeStep_j)} \le \HTUVehicle_{\Vehicle}, \HTUOps_{(\Stop_i,\TimeStep_i),(\Stop_j, \TimeStep_j)} \in \mathbb{N} \cup \{0\}$ of \glspl{acr:htu} to transporting freight between two consecutive freight terminals \mbox{$(\Stop_i,\TimeStep_i),(\Stop_j, \TimeStep_j) \in \SequenceOfRoute_{\Vehicle}, \Stop_i, \Stop_j \in \SetOfFlatTerminals$} on every vehicles' route $\SequenceOfRoute_{\Vehicle}$. We implicitly assume either zero or a minimum of $2$ \glspl{acr:ft} on every vehicle's route.
    \item[iii)] selecting a subset of freight requests accepted for transportation via the \gls{acr:pts}.
    \item[iv)] assigning flows $0 \le \PassengerFlow^{\Request}_{\Path} \le \Demand^{\Request}, \Request \in \SetOfRequests^{\Passenger}, \Path \in \SetOfPaths(\Request)$ that determine one or multiple paths to partially service the respective passenger request. Note that we allow the split of flows to reflect the various journey patterns of individual passengers.
    \item[v)] assigning flows  $\FreightFlow^{\Request}_{\Path} \in \{0,\Demand^{\Request}\}, \Request \in \SetOfRequests^{\Freight}, \Path \in \SetOfPaths(\Request)$ that determine a singular feasible path for every accepted freight request.
\end{compactitem}

\noindent\textbf{Solution: } A feasible, well-defined solution adheres to the following constraints:
\begin{compactitem}
    \item[i)] it preserves a passenger service level $\ServiceLevel \in [0,1]$ ensuring that an exogenously given share of passenger demand is serviced, i.e.,  $\sum_{\Request \in \SetOfRequests^{\Passenger}} \sum_{\Path \in \SetOfPaths(\Request)} \PassengerFlow^{\Request}_{\Path} \ge \ServiceLevel  \sum_{\Request \in \SetOfRequests^{\Passenger}} \Demand^{\Request}$ holds.
    \item[ii)] it respects the passenger capacity of every \gls{acr:pt} vehicle at all times. The passenger capacity is determined by the status quo adjusted by the capacity allocated to freight transportation, i.e., $\UnitCapacity_{\Vehicle} (\NumUnits_{\Vehicle} - \HTUOps_{(\Stop_i,\TimeStep_i),(\Stop_j, \TimeStep_j)}), \Vehicle \in \SetOfVehicles$. Here, $\Stop_i, \Stop_j$ denote any pair of consecutive \glspl{acr:ft} on the respective vehicle's route $\SequenceOfRoute_{\Vehicle}$.
    \item[iii)] it respects the freight capacity of every \gls{acr:pt} vehicle at all times. The freight capacity is determined by the allocated \gls{acr:htu} decision, i.e., the capacity between any pair of consecutive \glspl{acr:ft} $\Stop_i, \Stop_j$ on a vehicle's route $\SequenceOfRoute_{\Vehicle}$ is $\UnitCapacity_{\Vehicle} \HTUOps_{(\Stop_i,\TimeStep_i),(\Stop_j, \TimeStep_j)}, \Vehicle \in \SetOfVehicles$.
\end{compactitem}

\noindent\textbf{Objective:} The municipality aims to minimize the total system cost with respect to cargo-hitching adoption. This cost entails multiple components:
\begin{compactitem}
    \item[i)] a design cost $\Cost_{\Vehicle} > 0$ per \gls{acr:htu} that is assigned to vehicle $\Vehicle \in \SetOfVehicles$ representing the normalized investment cost that is scaled to the investigated time period. 
    \item[ii)] a penalty cost $\Cost^{\Request}_{\Penalty} > 0$ for every freight request $\Request \in \SetOfRequests^{\Freight}$ rejected by the municipality and representing the cost of negative externalities due to conventional truck delivery.
    \item[iii)] a routing cost $\Cost_{\Stop_i, \Stop_{i+1}} > 0$ per unit freight that the \gls{acr:pts} transports between consecutive stops $\Stop_i, \Stop_{i+1} \in \SetOfFlatStops$ on any vehicle $\Vehicle$'s route $\SequenceOfRoute_{\Vehicle}$.
    \item[iv)] \revone{a cost per unit freight for handling operations that are required to unload freight from trucks into the \gls{acr:pts}, load and unload freight to and from \gls{acr:pt} vehicles, and load freight from the \gls{acr:pts} into city-freighters to perform the last mile.}
    \item[v)] \revone{a last-mile cost per unit freight for every request $\Request \in \SetOfRequests^{\Freight}$ that represents the externality costs that occur due to the last-mile delivery of the request from the last \gls{acr:pt} stop on its selected path $\Stop^{\Path}_{|\Path|-2}$ to its destination $\Destination^{\Request}$.}
\end{compactitem}

Note that we do not account for passenger transportation costs in the objective because cost differences between passenger paths of acceptable quality are marginal.

\noindent\revone{\textbf{Discussion: } Three remarks on our problem setting and its justification are in order.

First, we acknowledge that in our solutions a very small fraction of passenger requests may remain formally unserved. This reflects a deliberate design choice and is consistent with \gls{acr:pts} operations in practice where passenger demand is either comfortably accommodated by available capacity or subject to temporary bottlenecks during peak hours. In the latter case, unserved requests rarely correspond to complete service denial but rather manifest in short delays and more transfers. While our model does not explicitly capture this temporal adjustment, it provides a realistic justification for allowing a negligible proportion of requests to remain served. Importantly, we enforce a strict service level $\ServiceLevel = 99.9$\% in our computational study, meaning that at most $0.1$\% of requests may remain unserved. This level of tolerance is within the typical range in reliability-oriented service planning. Moreover, it implies that our modeling framework is conservative with respect to freight integration and any potential effects of prioritizing freight over passengers remain negligible. Future research could extend the model by explicitly accounting for temporal displacement effects, for example through schedule-based or queuing formulations.

Second, while our model includes an approximate representation of last-mile costs, it does not explicitly account for first-mile externalities. This choice reflects the geographical and operational context: last-mile deliveries typically occur within the vaster city center area, where externalities such as congestion and air pollution are more pronounced, whereas first-mile operations typically originate from \glspl{acr:lsp}’ depots located in the urban periphery, where such impacts are less severe. Moreover, the approximation error associated with modeling last-mile costs is likely to be smaller than that for first-mile costs, which involve longer distances and more variation due to the use of heavier vehicles and diverse depot locations. Nonetheless, our findings represent an upper bound on the cost reductions achievable with shared passenger–freight systems.

Third, although our model imposes a system-wide passenger service-level constraint, it does not explicitly ensure an equitable distribution of the unserved proportion across individual lines, user-groups, or regions. This limitation could lead to service imbalances, where less profitable routes receive disproportionately lower levels of passenger service. To address this concern, several model extensions could be interesting for future research. On the one hand, the model could be extended to account for more granular service constraints that extend the service level requirement across given subsets of passenger demands $\SetOfRequests^{\Passenger}$. On the other hand, increased fairness could be imposed by considering an equity-based component in the objective, e.g., a fairness penalty for large disparities in service levels or a second objective within a multi-objective framework that explicitly determines the Pareto frontier between economic sustainability and fairness without the need to monetize disparities.
}
\subsection{Expanded graph construction} \label{subsection:graph-construction}
\begin{figure}[bt!]
\centering{
    \begin{subfigure}[b]{0.48\textwidth}
            \centering
            \def\svgwidth{\textwidth}
            \resizebox{.85\textwidth}{!}{\input{Figure_3a.tex}}
            \caption{Partial temporal expansion (different vehicle routes separated by different colors)}
            \label{fig:time-expansion}
    \end{subfigure}
    \quad
    \begin{subfigure}[b]{0.48\textwidth} 
            \def\svgwidth{\textwidth}
            %% Creator: Inkscape 1.2.2 (732a01da63, 2022-12-09), www.inkscape.org
%% Accompanies image file '20230227 multi layered graph.pdf' (pdf, eps, ps)
\begingroup%
  \makeatletter%
  \providecommand\color[2][]{%
    \errmessage{(Inkscape) Color is used for the text in Inkscape, but the package 'color.sty' is not loaded}%
    \renewcommand\color[2][]{}%
  }%
  \providecommand\transparent[1]{%
    \errmessage{(Inkscape) Transparency is used (non-zero) for the text in Inkscape, but the package 'transparent.sty' is not loaded}%
    \renewcommand\transparent[1]{}%
  }%
  \providecommand\rotatebox[2]{#2}%
  \newcommand*\fsize{\dimexpr\f@size pt\relax}%
  \newcommand*\lineheight[1]{\fontsize{\fsize}{#1\fsize}\selectfont}%
  \ifx\svgwidth\undefined%
    \setlength{\unitlength}{530.98614671bp}%
    \ifx\svgscale\undefined%
      \relax%
    \else%
      \setlength{\unitlength}{\unitlength * \real{\svgscale}}%
    \fi%
  \else%
    \setlength{\unitlength}{\svgwidth}%
  \fi%
  \global\let\svgwidth\undefined%
  \global\let\svgscale\undefined%
  \makeatother%
  \begin{picture}(1,0.4896836)%
    \lineheight{1}%
    \setlength\tabcolsep{0pt}%
    \put(0.82421857,0.05263769){\color[rgb]{0,0,0}\makebox(0,0)[lt]{\lineheight{1.25}\smash{\begin{tabular}[t]{l}$\mathcal{S}_{\Hold}$\end{tabular}}}}%
    \put(0.82421857,0.23709778){\color[rgb]{0,0,0}\makebox(0,0)[lt]{\lineheight{1.25}\smash{\begin{tabular}[t]{l}$\mathcal{S}_{1}$\end{tabular}}}}%
    \put(0.82421857,0.42155789){\color[rgb]{0,0,0}\makebox(0,0)[lt]{\lineheight{1.25}\smash{\begin{tabular}[t]{l}$\mathcal{S}_{2}$\end{tabular}}}}%
    \put(0,0){\includegraphics[width=\unitlength,page=1]{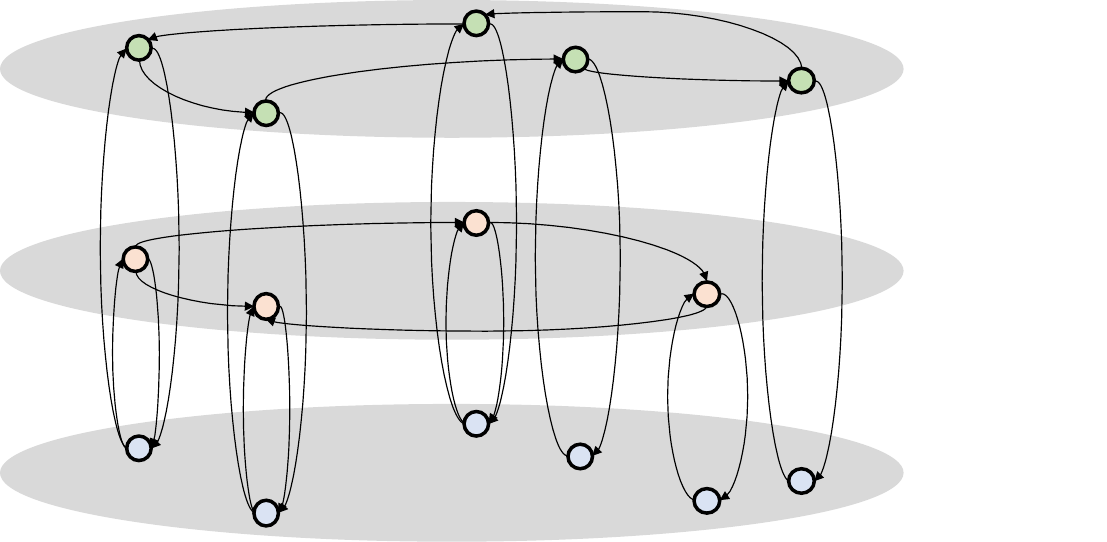}}%
  \end{picture}%
\endgroup%

            \caption{Spatial expansion}
            \label{fig:spatial-expansion}
    \end{subfigure}%
}
\caption{Schematic illustration of the spatial and temporal graph expansion elements} 
\label{fig:network-expansion}
\end{figure}

\noindent To devise an effective algorithm, we encode some of the problem's temporal and spatial complexity by using a problem-specific graph representation. Specifically, we use a temporal graph expansion in which vertices represent a combination of location and time (cf. Figure~\ref{fig:time-expansion}) and combine it with a spatial expansion in which we separate different vehicles' routes through the \gls{acr:pts} into $|\SetOfVehicles| + 1$ different graph layers (cf. Figure~\ref{fig:spatial-expansion}). The expanded graph contains one separate layer of temporally expanded vertices for every vehicle $\Vehicle \in \SetOfVehicles$, and one additional layer that we refer to as holding layer. We denote the resulting multi-layered graph with its vertex and arc sets as $\Graph=(\SetOfVertices, \SetOfArcs)$. We finish the formal construction of $\Graph$ by providing an example at the end of this section.

\noindent \textbf{Vertex set construction: }Let \revone{$\SetOfTimeSteps:= \bigcup_{\Request \in \SetOfRequests} (\RequestStart^{\Request} \cup \RequestEnd^{\Request}) \cup \bigcup_{\Stop \in \SetOfFlatStops} \{ \TimeStep: \TimeStep \in \SetOfTimeSteps(\Stop) \}$} be the set of all relevant timestamps in a given problem instance. Furthermore, let the vertex set $\SetOfVertices:= \SetOfStops \cup \SetOfOrigins \cup \SetOfDestinations$ be the union of the subsets $\SetOfStops$, $\SetOfOrigins$, and $\SetOfDestinations$, which we define one-by-one in the remainder of this section. 

First, we obtain $\SetOfOrigins$ and $\SetOfDestinations$ by expanding each request's origin and destination into the time dimension. Thus, $\SetOfOrigins:= \bigcup_{\Request \in \SetOfRequests} (\Origin^{\Request}, \RequestStart^{\Request})$ and $\SetOfDestinations:= \bigcup_{\Request \in \SetOfRequests} (\Destination^{\Request}, \RequestEnd^{\Request})$. 

Second, we describe a vertex $\Vertex \in \SetOfStops$ as a temporal node represented by a triplet that links a physical location $\Stop \in \SetOfFlatStops$ to a specific timestamp $\TimeStep \in \SetOfTimeSteps$ and a specific graph layer $\{\Hold\} \cup \SetOfVehicles$. Accordingly, we define sets $\SetOfStops_{\Vehicle} := \{(\Stop, \TimeStep, \Vehicle): (\Stop, \TimeStep) \in \SequenceOfRoute_{\Vehicle}\}, \medspace \Vehicle \in \SetOfVehicles$ as the vertices at which vehicle $\Vehicle$ arrives at stop $\Stop$ at timestamp $\TimeStep$. This is a partial time expansion in which only the relevant points in time are expanded \citep[cf.][]{BolandHewittEtAl2017}. Here, every set $\SetOfStops_{\Vehicle}$ denotes the vertices of a different vehicle layer. Additionally, we introduce a specific holding layer vertex set $\SetOfStops_{\Hold} := \{ (\Stop, \TimeStep, \Hold): (\Stop, \TimeStep) \in \SetOfRoutes \}$ that contains one additional copy per stop $\Stop$ and timestamp $\TimeStep$ in which a vehicle arrives at $\Stop$. The holding layer connects the vehicle layers and orchestrates the time synchronization of transfers and transshipments. Finally, we define $\SetOfStops:= \SetOfStops_{\Hold} \medspace \cup \medspace \bigcup_{\Vehicle \in \SetOfVehicles} \SetOfStops_{\Vehicle}$ as the set of all temporally expanded vertices that represent stops across all graph layers. In the expanded graph, we refer to the set of \gls{acr:ft} representations as $\SetOfTerminals := \{ (\Stop, \cdot, \cdot) \in \SetOfStops: \Stop \in \SetOfFlatTerminals \}$.

\noindent \textbf{Arc set construction: }We create the resulting graph's arc set $\SetOfArcs$ in a top-down manner. More specifically, we derive the set of arcs $\SetOfArcs$ in the expanded graph $\Graph$ as the union of multiple disjoint arc subsets, i.e., $\SetOfArcs:= \SetOfArcs^{\VehicleLayer} \cup \SetOfArcs^{\Hold} \cup \SetOfArcs^{\Transit} \cup \SetOfArcs^{\Access} \cup \SetOfArcs^{\Egress}$. 

Vehicle arcs $\Arc \in \SetOfArcs^{\VehicleLayer}$ complete the vehicle layer vertex sets in the multi-layered graph~$\Graph$ and represent the \gls{acr:pt} vehicles' routes. We define the set $\SetOfArcs^{\VehicleLayer} := \bigcup_{\Vehicle \in \SetOfVehicles} \SetOfArcs_{\Vehicle}$ as the union of temporal arc sets $\SetOfArcs_{\Vehicle}$ that contain the arcs representing the route of vehicle $\Vehicle \in \SetOfVehicles$ in its corresponding graph layer. Here, we construct the arcs $\Arc \in \SetOfArcs_{\Vehicle}$ such that they connect consecutive stops $(\Stop_{l}, \TimeStep_{l}), (\Stop_{l+1}, \TimeStep_{l+1})$ on a vehicle's route $\SequenceOfRoute_{\Vehicle}$. More formally, we construct arcs $\Arc \in \SetOfArcs_{\Vehicle}$ such that $i = (\Stop_{l}, \TimeStep_{l}, \Vehicle) \in \SetOfStops_{\Vehicle}, \medspace j = (\Stop_{l+1}, \TimeStep_{l+1}, \Vehicle) \in \SetOfStops_{\Vehicle}, \medspace l \in \{1, ..., n-1 \}$. 

The arc set $\SetOfArcs^{\Hold}$ is the set of holding arcs that enable holding requests at stops of the \gls{acr:pts}, i.e., passengers or freight waiting at a stop. Let $u$ be the index of the ordered \revone{list} of times $\SetOfTimeSteps(\Stop)$ in which any vehicle arrives at stop $\Stop$. Then, we create holding arcs $\Arc \in \SetOfArcs^{\Hold}$ where \revone{$i=(\Stop, \TimeStep^{\Stop}_{u}, \Hold) \in \SetOfStops_{\Hold}$ and $j = (\Stop, \TimeStep^{\Stop}_{u+1}, \Hold) \in \SetOfStops_{\Hold}$, for all $\Stop \in \SetOfFlatStops$, $\TimeStep^{\Stop}_{u}, \TimeStep^{\Stop}_{u+1} \in \SetOfTimeSteps(\Stop)$,} and $u=1,...,|\SetOfTimeSteps(\Stop)|-1$. Thus, we connect vertices in the holding layer vertex set that are copies of the same physical stop location $\Stop \in \SetOfFlatStops$ such that two connected vertices represent two consecutive timestamps \revone{$\TimeStep^{\Stop}_{u}, \TimeStep^{\Stop}_{u+1} \in \SetOfTimeSteps(\Stop)$} in which a \gls{acr:pt} vehicle arrives at the respective stop.

To connect the disjunct vertex sets of the different layers, we add transit arcs $\Arc~\in~\SetOfArcs^{\Transit}$. Here, we connect temporal vertices $i=(\Stop, \TimeStep, \Vehicle) \in \SetOfStops_{\Vehicle}, \medspace \Vehicle \in \SetOfVehicles$ with their corresponding representation in the holding layer $j=(\Stop, \TimeStep, \Hold) \in \SetOfStops_{\Hold}$. This is a many-to-one mapping as multiple vertices in the vehicle layer vertex set may share the same representation in the holding layer, e.g., if they represent the same physical stop at the same timestamp. To establish this mapping in a bidirectional fashion, we further add the inverse arc $(j, i)$. 

Moreover, we connect the temporal vertices $i=(\Origin^{\Request}, \RequestStart^{\Request}) \in \SetOfOrigins, \Request \in \SetOfRequests^{\Passenger}$ with the \gls{acr:pts} at vertices $j=(\Stop, \TimeStep, \Hold) \in \SetOfStops_{\Hold}$ by arcs $\Arc \in \SetOfArcs^{\Access}$. Accordingly, we connect the temporal vertices $i=(\Origin^{\Request}, \RequestStart^{\Request}) \in \SetOfOrigins, \Request \in \SetOfRequests^{\Freight}$ with the \gls{acr:pts} at vertices $j=(\Stop, \TimeStep, \Hold) \in \SetOfStops_{\Hold} \cap \SetOfTerminals$ by additional arcs $\Arc \in \SetOfArcs^{\Access}$. Finally, we construct arcs $\Arc \in \SetOfArcs^{\Egress}$ that connect the \gls{acr:pts} vertices $i=(\Stop, \TimeStep, \Hold) \in \SetOfStops_{\Hold}$ with the temporal destinations $j=(\Destination^{\Request}, \RequestEnd^{\Request}) \in \SetOfDestinations, \Request \in \SetOfRequests^{\Passenger}$ and additional arcs $\Arc \in \SetOfArcs^{\Egress}$ that connect the \gls{acr:pts} vertices $i=(\Stop, \TimeStep, \Hold) \in \SetOfStops_{\Hold} \cap \SetOfTerminals$ with the temporal destinations $j=(\Destination^{\Request}, \RequestEnd^{\Request}) \in \SetOfDestinations, \Request \in \SetOfRequests^{\Freight}$. Here, we prune the graph based on distance and time thresholds as outlined in Appendix~\ref{app:implementation}. 

\noindent \textbf{Preprocessing: }We apply multiple preprocessing steps to encode problem characteristics in $\Graph$, which reduces the size and the computational complexity of the \gls{acr:mip} formulation in Section~\ref{subsection:milp}. 

First, we encode the system operator's acceptance decisions on transporting freight requests $\Request \in \SetOfRequests^{\Freight}$ into the expanded and multi-layered graph by constructing dummy arcs \mbox{$\Arc \in \SetOfArcs^{\Dummy} \subset \SetOfArcs$} such that the decision to reject a freight request $\Request$ corresponds to routing it through the network on a dummy arc. Formally, we add dummy arcs $\Arc \in \SetOfArcs^{\Dummy}$, with $i=(\Origin^{\Request}, \RequestStart^{\Request})$ and $j=(\Destination^{\Request}, \RequestEnd^{\Request})$ for all $\Request \in \SetOfRequests^{\Freight}$ and assign the arc cost $\Cost_{\UnbracedArc}=\frac{1}{\Demand^{\Request}} \Cost^{\Request}_{\Penalty}$ such that the encoded routing cost equals the penalty cost from rejecting a request. Because we enforce binary freight flows, we can decode the decision to accept or reject a request from the network flow without further intricacies. 

Second, to reduce the problem's complexity, we aim at decreasing the cardinality of the arc set. Furthermore, the required capacity for freight transportation remains stable at \gls{acr:pt} stops that are not \glspl{acr:ft}, i.e., $\forall \Stop \in \SetOfFlatStops \setminus \SetOfFlatTerminals$ as no freight can enter or leave the \gls{acr:pts} at such stops. We leverage this observation and further abstract the \gls{acr:pts} by constructing freight path segment arcs that connect consecutive \glspl{acr:ft} on every vehicle's route --- thereby contracting multiple arcs into a single arc. Let $m \ge 1$, where $m-1$ is the number of \gls{acr:pt} stops between the two consecutive freight terminals. Furthermore, let $i=(\Stop_{l}, \TimeStep_{l}, \Vehicle) \in \SetOfStops_{\Vehicle} \cap \SetOfTerminals$ and $j=(\Stop_{l+m}, \TimeStep_{l+m}, \Vehicle) \in \SetOfStops_{\Vehicle} \cap \SetOfTerminals$ be the two vertices in the vehicle layer of a vehicle $\Vehicle$ representing the two consecutive \glspl{acr:ft} on the vehicle's route. Formally, $(\Stop_{l}, \TimeStep_{l}), (\Stop_{l+m}, \TimeStep_{l+m}) \in \SequenceOfRoute_{\Vehicle}$ such that $\nexists \hspace{0.1cm} (\Stop_{l+p}, \TimeStep_{l+p}) \in \SequenceOfRoute_{\Vehicle}: 0 < p < m, \hspace{0.1cm} \Stop_{l+p} \in \SetOfFlatTerminals$. For all such $\UnbracedArc$, we add the freight path segment arcs $\Arc \in \SetOfArcs^{\PathSegment} \subset \SetOfArcs$. Here, we make sure that costs are consistent by setting $\Cost_{\UnbracedArc} = \sum^{m}_{p=1} \Cost_{\Stop_{l+p-1},\Stop_{l+p}}$. 

Third, we pre-compute sets of passenger paths $\SetOfPaths(\Request), \Request \in \SetOfRequests^{\Passenger}$ \citep[cf.][]{LiZhuEtAl2024}. Hence, we reduce the problem's computational complexity by converting the minimum cost network flow problem for passengers into a relaxed fractional set covering problem. Thus, we can ensure suitable passenger service criteria based on the pre-computed sets by filtering for the chosen criteria, e.g., filtering for paths that induce a maximum number of transfers or using specific modes in a multi-modal setting. \revone{Filtering based on a set of precomputed paths therefore allows us to mitigate concerns about increased transfers in integrated passenger-freight systems by explicitly controlling for transfer frequency and passenger travel experience.} 

\noindent \textbf{Graph expansion example: }The following provides an illustrative example of a \gls{acr:pts} and its graph representation. In this example, the \gls{acr:pts} consists of two vehicles $\SetOfVehicles = \{1,2\}$ operating the following routes: $\SetOfRoutes_{1} = \langle(s_1,2), (s_2,3), (s_3,4), (s_4,6)\rangle$ and $\SetOfRoutes_{2} = \langle(s_5,1), (s_2,2), (s_3,3), (s_6,4)\rangle$. 
We assume the set of \glspl{acr:ft} to be $\SetOfFlatTerminals = \{s_{1}, s_{2}, s_{4}, s_{5}, s_{6} \}$, and consider two requests $\SetOfRequests = \{1,2\}$, both of which come with a demand of $\Demand^{1}=\Demand^{2}=1$. 
Freight request $1 \in \SetOfRequests^{\Freight}$ can leave its origin no earlier than time $0$ and must be fulfilled with $t=6$. 
Request $2 \in \SetOfRequests^{\Passenger}$ is a passenger request, and the corresponding group of passengers wants to start their itinerary no earlier than time $\TimeStep=0$ and finish it before or at $\TimeStep=5$. 
Figure \ref{fig:example-expanded-graph} shows the corresponding expanded graph. In the preprocessing, we have added dummy arc $((\Origin^{1}, 0),(\Destination^{1}, 5)) \in \SetOfArcs^{\Dummy}$ as well the freight path segments $((s_1, 2, 1),(s_2, 3, 1)), ((s_2, 3, 1),(s_4, 6, 1)), ((s_5, 1, 2),(s_2, 2, 2)), ((s_2, 2, 2),(s_6, 4, 2)) \in \SetOfArcs^{\PathSegment}$.

\begin{figure}[!t]
  \begin{subfigure}{0.77\textwidth}
    \centering
    \customesmallfontsize
    \DecimalMathComma
    \def\svgwidth{\textwidth}
    \input{Figure_4.tex}
  \end{subfigure}%
  \begin{subfigure}{0.23\textwidth}
    \centering
    \scriptsize
    \def\svgwidth{\textwidth}
    \input{Figure_4_Legend.tex}
  \end{subfigure}
  \caption{Illustrative example of a partial temporal expanded, multi-layered, and preprocessed graph}
  \label{fig:example-expanded-graph}
\end{figure}

\subsection{MIP formulation}\label{subsection:milp}
\noindent Based on the expanded, multi-layered, and pre-processed graph $\Graph=(\SetOfVertices, \SetOfArcs)$ introduced in Section~\ref{subsection:graph-construction}, we formulate the problem as a \gls{acr:mip} in this section. In Appendix~\ref{app:notation}, we provide a tabular summary of notation.
For the ease of notation, we first explicitly define the arc set $\SetOfArcs^{\PotentialFreightArcs} \subseteq \SetOfArcs$ that we consider when determining freight flows. This arc set contains the freight path segments $\SetOfArcs^{\PathSegment}$, the dummy arcs $\SetOfArcs^{\Dummy}$, the connections between \gls{acr:ft} representations from $\SetOfArcs^{\Transit}$ and $\SetOfArcs^{\Hold}$, and the relevant connections of origins and destinations with the \gls{acr:pts}. More formally, 
\begin{align*}
  \SetOfArcs^{\PotentialFreightArcs} := \SetOfArcs^{\PathSegment} \hspace{0.1cm} \cup \hspace{0.1cm} \SetOfArcs^{\Dummy} \hspace{0.1cm} \cup \hspace{0.1cm} & \{ ((\Origin^{\Request}, \RequestStart^{\Request}),(\Stop, \TimeStep, \Hold)) \in \SetOfArcs^{\Access}: \Request \in \SetOfRequests^{\Freight} \} \hspace{0.1cm} \cup \notag \\ 
  &\{ ((\Stop, \TimeStep, \Hold),(\Destination^{\Request}, \RequestEnd^{\Request})) \in \SetOfArcs^{\Egress}: \Request \in \SetOfRequests^{\Freight} \} \hspace{0.1cm} \cup \notag \\ 
  &\{ \Arc \in \SetOfArcs^{\Hold} \cup \SetOfArcs^{\Transit}: i,j \in \SetOfTerminals \} 
\end{align*} 

Note that a graph decomposition by commodity type, i.e., passenger and freight, is not straightforward because both commodities share a total capacity that needs explicit allocation. We define the vertex demand $\VertexDemand^{\Request}_{i}, \medspace i \in \SetOfVertices, \Request \in \SetOfRequests^{\Freight}$ describing the difference between total inflow and total outflow of a specific request in a vertex as
\begin{equation*}
    \VertexDemand^{\Request}_{i} = \begin{cases}
        1, & \text{if} \quad i = (\Origin^{\Request}, \RequestStart^{\Request}),\\%\in \SetOfOrigins^{\Exp},\\
        -1,& \text{if} \quad i = (\Destination^{\Request}, \RequestEnd^{\Request}),\\%\in \SetOfDestinations^{\Exp},\\
        0 & \text{otherwise.}
    \end{cases} 
\end{equation*} 

In this context, we restrict the vertex sets $\Neighbors(i)$ to neighboring vertices of vertex $i$ that are directly connected via arcs from $\SetOfArcs^{\PotentialFreightArcs}$, i.e., arcs that allow for freight transportation. Formally, $\Neighbors^{+}(i):= \{ j \in \SetOfVertices: \Arc \in \SetOfArcs^{\PotentialFreightArcs}\}$ and $\Neighbors^{-}(i):= \{ j \in \SetOfVertices: (j,i) \in \SetOfArcs^{\PotentialFreightArcs}\}$, respectively. 
We refer by $\SpanningMapping: \SetOfArcs^{\VehicleLayer} \rightarrow \SetOfArcs^{\PathSegment} \cup \hspace{2pt} \{ \emptyset \}$ to the many-to-one mapping function that assigns a contracted vehicle arc to its corresponding freight path segment. 
Some arcs in the vehicle layers might not be contracted.  
Therefore, we differentiate between arcs from the vehicle layer arc set that are contracted and arcs that are not being contracted. Here, $\SpannedArcs := \{ \Arc \in \SetOfArcs^{\VehicleLayer}: \SpanningMapping(\UnbracedArc) \neq \emptyset \}$ denotes the contracted arcs, and vice versa $\NotSpannedArcs = \{ \Arc \in \SetOfArcs^{\VehicleLayer}: \SpanningMapping(\UnbracedArc) = \emptyset \}$ denotes the arcs that are not contracted. In Figure~\ref{fig:arc-contraction}, we present a simplified illustration of this formalization with six temporal vertices, of which only three represent \glspl{acr:ft}. 
Table~\ref{tbl:arc-contraction} demonstrates the resulting arc sets. 

Finally, using the introduced notation, we can formulate the problem as a \gls{acr:mip} as follows
% \vspace{10pt}
\begin{mini!}|s|[3]<b> % mini! = minimize 
    {\HTUVehicle, \HTUOps, \PassengerFlow, \FreightFlow}  % optimization variable
    {\sum_{\Vehicle \in \SetOfVehicles} \Cost_{\Vehicle} \HTUVehicle_{\Vehicle} + 
    \sum_{\Request \in \SetOfRequests^{\Freight}} \Demand^{\Request} \sum_{\Arc \in \SetOfArcs^{\PotentialFreightArcs}} \Cost_{\UnbracedArc} \FreightFlow^{\Request}_{\UnbracedArc} \label{obj:func-primal}} 
    {\label{MIP}} % label for optimization problem
    {} % optimization result
    \addConstraint{\sum_{\Request \in \SetOfRequests^{\Passenger}} \sum_{\Path \in \SetOfPaths(\Request)} \Demand^{\Request} \PassengerFlow^{\Request}_{\Path}}{\geq \ServiceLevel \sum_{\Request \in \SetOfRequests^{\Passenger}} \Demand^{\Request} \label{cons:service-level}}
    \addConstraint{\sum_{j \in \Neighbors^{+}(i)} f^{\Request}_{\UnbracedArc} - \sum_{j \in \Neighbors^{-}(i)} \FreightFlow^{\Request}_{j,i}}{= \VertexDemand^{\Request}_{i},\quad} {\forall \Request \in \SetOfRequests^{\Freight}, i \in \SetOfOrigins \cup \SetOfDestinations \cup \SetOfTerminals \label{cons:flow-conservation}}  
    \addConstraint{\sum_{\Request \in \SetOfRequests^{\Passenger}} \sum_{\Path \in \SetOfPaths(\Request)} \Demand^{\Request} \PassengerFlow^{\Request}_{\Path} \PathArcRelation^{\Path}_{\UnbracedArc}}{\leq \sum_{\Vehicle \in \SetOfVehicles} \VehicleArcRelation^{\Vehicle}_{\UnbracedArc} \UnitCapacity_{\Vehicle} (\NumUnits_{\Vehicle} - \HTUOps_{\SpanningMapping(\UnbracedArc)}),\quad} {\forall \Arc \in \SpannedArcs \label{cons:passenger-capacity-spanned}} 
    \addConstraint{\sum_{\Request \in \SetOfRequests^{\Passenger}} \sum_{\Path \in \SetOfPaths(\Request)} \Demand^{\Request} \PassengerFlow^{\Request}_{\Path} \PathArcRelation^{\Path}_{\UnbracedArc}}{\leq \sum_{\Vehicle \in \SetOfVehicles} \VehicleArcRelation^{\Vehicle}_{\UnbracedArc} \UnitCapacity_{\Vehicle} \NumUnits_{\Vehicle},\quad} {\forall \Arc \in \NotSpannedArcs \label{cons:passenger-capacity-not-spanned}}
    \addConstraint{\sum_{\Request \in \SetOfRequests^{\Freight}} \Demand^{\Request} \FreightFlow^{\Request}_{\UnbracedArc}}{\leq \sum_{\Vehicle \in \SetOfVehicles} \VehicleArcRelation^{\Vehicle}_{\UnbracedArc} \UnitCapacity_{\Vehicle} \HTUOps_{\UnbracedArc},\quad}{\forall \Arc \in \SetOfArcs^{\PathSegment} \label{cons:freight-capacity}} 
    \addConstraint{\sum_{\Path \in \SetOfPaths(\Request)} \PassengerFlow^{\Request}_{\Path}}{\leq 1,\quad}{\forall \Request \in \SetOfRequests^{\Passenger} \label{cons:pass-convexity}} 
    \addConstraint{\HTUOps_{\UnbracedArc}}{\leq \sum_{\Vehicle \in \SetOfVehicles} \VehicleArcRelation^{\Vehicle}_{\UnbracedArc} \HTUVehicle_{\Vehicle},\quad}{\forall \Arc \in \SetOfArcs^{\PathSegment} \label{cons:cap-assignment}}
    \addConstraint{\HTUVehicle_{\Vehicle}}{\leq \NumUnits_{\Vehicle},\quad}{\forall \Vehicle \in \SetOfVehicles \label{cons:cap-total}}
    \addConstraint{\PassengerFlow^{\Request}_{\Path}}{\geq 0,\quad}{\forall \Request \in \SetOfRequests^{\Passenger}, \Path \in \SetOfPaths(\Request) \label{cons:domain-pass-flow}}
    \addConstraint{\FreightFlow^{\Request}_{\UnbracedArc}}{\in \{0,1\},\quad}{\forall \Request \in \SetOfRequests^{\Freight}, \Arc \in \SetOfArcs^{\PotentialFreightArcs} \label{cons:domain-freight-flow}}
    \addConstraint{\HTUVehicle_{\Vehicle}}{\in \mathbb{N}_0,\quad}{\forall h \in \SetOfVehicles \label{cons:domain-investment}}
    \addConstraint{x_{\UnbracedArc}}{\in \mathbb{N}_0,\quad}{\forall \Arc \in  \SetOfArcs^{\PathSegment}} \label{cons:domain-mode}
\end{mini!}

\noindent \revone{where $\PathArcRelation^{\Path}_{\UnbracedArc} = 1$ if $\Arc \in \SetOfArcs$ belongs to path $\Path \in \SetOfPaths$, and $\PathArcRelation^{\Path}_{\UnbracedArc} = 0$ otherwise. Similarly, $\VehicleArcRelation^{\Vehicle}_{\UnbracedArc} = 1$ if vehicle $\Vehicle \in \SetOfVehicles$ operates on arc $\Arc \in \SetOfArcs^{\Vehicle}$, and $\VehicleArcRelation^{\Vehicle}_{\UnbracedArc} = 0$ otherwise.}

\begin{figure}[!b]
\begin{floatrow}
\ffigbox{
    \scriptsize
    \def\svgwidth{0.98\columnwidth}
    \input{Figure_5.tex}
}{
    \caption{Arc contraction illustration}
    \label{fig:arc-contraction}
}
\capbtabbox{%
    \begin{tabular}{lc}
Set              & Elements              \\ \hline
\toprule
$\SpannedArcs$ & $(2,3), (3,4), (4,5)$ \\ \hline
$\NotSpannedArcs$ & $(1,2), (5,6)$        \\ \hline
\toprule
\caption{Arc contraction sets}
\label{tbl:arc-contraction}
\end{tabular}
}{}
\end{floatrow}
\end{figure}

The Objective~\eqref{obj:func-primal} minimizes freight transportation costs, which entail investment cost for \glspl{acr:htu} to design a network with flexible capacities, constant penalty terms for rejected freight requests, encoded as a routing cost on the arc subset $\SetOfArcs^{\Dummy} \subset \SetOfArcs^{\PotentialFreightArcs}$, and variable costs for every unit of freight flow transported via the \gls{acr:pts}. Constraint~\eqref{cons:service-level} ensures the demand-weighted passenger service level. Constraints~\eqref{cons:flow-conservation} impose classical flow conservation for every freight request $\Request \in \SetOfRequests^{\Freight}$. Constraints~\eqref{cons:passenger-capacity-spanned} and \eqref{cons:passenger-capacity-not-spanned} ensure that the system capacity in terms of passengers is respected on the contracted arcs and the non-contracted arcs, respectively. Constraints~\eqref{cons:freight-capacity} restrict the freight flow per arc $\Arc \in \SetOfArcs^{\PathSegment}$ depending on the number of \glspl{acr:htu} whose capacities are allocated to freight transportation on the respective freight path segment. Constraints~\eqref{cons:pass-convexity} restrict the passenger flow per request to one. Constraints~\eqref{cons:cap-assignment} and \eqref{cons:cap-total} limit the number of \glspl{acr:htu} allocated for freight transportation. Specifically, Constraints~\eqref{cons:cap-assignment} restrict this number to the \glspl{acr:htu} assigned to vehicle $\Vehicle \in \SetOfVehicles$, while Constraints~\eqref{cons:cap-total} limit the number of \glspl{acr:htu} to the maximum number of units $\NumUnits_{\Vehicle}$ available. Finally, Constraints~\eqref{cons:domain-pass-flow} - \eqref{cons:domain-mode} define the domains of the decision variables.

\section{Algorithm}\label{sec:algorithmic-framework}
\noindent We propose a \gls{acr:pab} solution method to find integer feasible solutions to Problem~\ref{MIP}. First, we reformulate Problem~\ref{MIP} as a path-based formulation in Section~\ref{subsection:path-based-model}. Second, we detail our \gls{acr:pab} approach in Section~\ref{subsection:p&b}. Finally, we provide a \gls{acr:bap} method as an alternative to our \gls{acr:pab} algorithm in Section~\ref{sec:branch-and-price}.

\subsection{Path-based reformulation}\label{subsection:path-based-model}
\noindent Let $\FreightFlowPath^{\Request}_{\Path}$ denote the fraction of demand $\Demand^{\Request}, \medspace \Request \in \SetOfRequests^{\Freight}$ that is transported via path $\Path \in \SetOfPaths(\Request)$. The path-based formulation required to apply \gls{acr:cg} reads 
\begin{mini!}|s|[3]<b>  
    {\HTUVehicle, \HTUOps, \PassengerFlow, \FreightFlowPath} 
    {\sum_{\Vehicle \in \SetOfVehicles} \Cost_{\Vehicle} \HTUVehicle_{\Vehicle} + 
    \sum_{\Request \in \SetOfRequests^{\Freight}} \sum_{\Path \in \SetOfPaths(\Request)} \sum_{\Arc \in \SetOfArcs^{\PotentialFreightArcs}} \Demand^{\Request} \Cost_{\UnbracedArc} \PathArcRelation^{\Path}_{\UnbracedArc} \FreightFlowPath^{\Request}_{\Path} \label{obj:path-based-formulation}} 
    {\label{problem:path-based-formulation}} 
    {} 
    \addConstraint{(\DualAlpha) \quad \sum_{\Request \in \SetOfRequests^{\Freight}} \sum_{\Path \in \SetOfPaths(\Request)} \Demand^{\Request} \PathArcRelation^{\Path}_{\UnbracedArc} \FreightFlowPath^{\Request}_{\Path}}{\leq \sum_{\Vehicle \in \SetOfVehicles} \VehicleArcRelation^{\Vehicle}_{\UnbracedArc} \UnitCapacity_{\Vehicle} \HTUOps_{\UnbracedArc},\quad}{\forall \Arc \in \SetOfArcs^{\PathSegment} \label{cons:freight-capacity-path-based}} 
    \addConstraint{(\DualEta) \quad \sum_{\Path \in \SetOfPaths(\Request)} \FreightFlowPath^{\Request}_{\Path}}{= 1,\quad}{\forall \Request \in \SetOfRequests^{\Freight} \label{cons:freight-convexity}} 
    \addConstraint{\quad \quad \FreightFlowPath^{\Request}_{\Path}}{\in \{0,1\},\quad}{\forall \Request \in \SetOfRequests^{\Freight}, \Path \in \SetOfPaths(\Request) \label{cons:domain-freight-flow-path-based}}
\end{mini!}
adhering further to Constraints~\ref{cons:service-level}, \ref{cons:passenger-capacity-spanned} - \ref{cons:passenger-capacity-not-spanned}, \ref{cons:pass-convexity} - \ref{cons:domain-pass-flow}, and \ref{cons:domain-investment}-\ref{cons:domain-mode}. In the continuous relaxation of Problem~\ref{problem:path-based-formulation} each set of constraints is associated with a set of dual variables. 
Here, the dual variables $\DualAlpha_{\UnbracedArc} \in \mathbb{R}^{-}_{0}, \medspace \Arc \in \SetOfArcs^{\PathSegment}$ are associated with Constraints~\ref{cons:freight-capacity-path-based} limiting the freight flow on the respective arcs. Constraints~\ref{cons:freight-convexity} enforce the sum of all freight flows per request to be equal to one. Thus, the associated dual variables $\DualEta^{\Request} \in \mathbb{R}, \medspace \Request \in \SetOfRequests^{\Freight}$ are free. Moreover, the dual variable $\DualGamma \in \mathbb{R}^{+}_{0}$ is associated with the service level Constraint~\ref{cons:service-level}. The dual variables $\DualUpsilon_{\UnbracedArc} \in \mathbb{R}^{-}_{0}, \medspace \Arc \in \SpannedArcs$ and $\DualNu_{\UnbracedArc} \in \mathbb{R}^{-}_{0}, \medspace \Arc \in \NotSpannedArcs$ are associated with the passenger capacity limiting Constraints~\ref{cons:passenger-capacity-spanned}, and Constraints~\ref{cons:passenger-capacity-not-spanned} respectively. Additionally, dual variables $\DualDelta^{\Request} \in \mathbb{R}^{-}_{0}, \medspace \Request \in \SetOfRequests^{\Passenger}$ are associated with Constraints~\ref{cons:pass-convexity}, dual variables $\DualPi_{\UnbracedArc} \in \mathbb{R}^{-}_{0}, \medspace \Arc \in \SetOfArcs^{\PathSegment}$ are associated with Constraints~\ref{cons:cap-assignment}, and dual variables $\DualTau_{\Vehicle} \in \mathbb{R}^{-}_{0}, \medspace \Vehicle \in \SetOfVehicles$ are linked to Constraints~\ref{cons:cap-total}. Based on the path-based Problem~\ref{problem:path-based-formulation}, we outline the components of our \gls{acr:pab} algorithm in the next section.

\subsection{Price-and-branch}\label{subsection:p&b}
\noindent Algorithm~\ref{algo:price-and-branch} shows a pseudocode of our \gls{acr:pab} approach. Contrary to \gls{acr:bap} where we iterate between pricing and branching, in \gls{acr:pab} we price once and then enforce integer feasible solutions without generating new columns.
\begin{figure}[bt] 
\vspace{-15pt}
    \begin{minipage}{0.98\linewidth} 
        \begin{algorithm}[H]
            \caption{\textbf{Price-and-branch}} \label{algo:price-and-branch}
            \footnotesize
            \begin{algorithmic}[1]
                \Require Path-based formulation \eqref{problem:path-based-formulation} 
                \State $\texttt{relaxation} \gets \texttt{ContinuousRelaxation}$\label{l:solve_formulation}
                \State $\texttt{rmp} \gets \texttt{InitializeRMP(relaxation)}$ \label{l:initrmp}\Comment{Ensures feasibility throughout \gls{acr:cg}}
                \State $\texttt{LB}, \texttt{UB} \gets 0, \infty$
                \While{$\texttt{OptimalityGap} > \epsilon$ and $\texttt{SolveTime} < \texttt{TimeLimit}$} \label{l:stopppingcriteria}
                    \State $\texttt{duals} \gets \texttt{SolveRMP(rmp)}$ \label{l:solvermp}\Comment{Warmstarting at previous solution}
                    \State $\texttt{cols} \gets \texttt{Price(duals)}$ \label{l:price}
                    \State $\texttt{rmp} \gets \texttt{AddColumns(rmp, cols)}$ \label{l:newcolumns}
                    \State $\texttt{UB} \gets \texttt{UpdateBounds(rmp)}$ \label{l:ub} \Comment{Solution value of RMP}
                    \If {$\texttt{FullPricingIteration}$} \label{l:lb1}\Comment{No update in partial pricing iterations}
                        \State $\texttt{LB} \gets \texttt{UpdateBounds(rmp, duals)}$\label{l:lb2} 
                    \EndIf
                \EndWhile \label{l:pricing-over-pab}
                \State $\texttt{solution} \gets \texttt{BranchAndCut(rmp)}$ \label{l:bac}\Comment{No further updates of lower bound} \\
                \Return $\texttt{solution}$
            \end{algorithmic}
        \end{algorithm}
    \end{minipage}
    \vspace{-5pt}
\end{figure}
Algorithm~\ref{algo:price-and-branch} solves the continuous relaxation of the given path-based formulation (l.~\ref{l:solve_formulation}) via \gls{acr:cg} with partial pricing. More specifically, the algorithms initializes a \gls{acr:rmp} (l.~\ref{l:initrmp}), and then iteratively solves the \gls{acr:rmp} (l.~\ref{l:solvermp}) and a pricing problem (l.~\ref{l:price}) in order to add new columns (l.~\ref{l:newcolumns}) to the \gls{acr:rmp}. Every $5$ iterations, we enhance the partial pricing by conducting a full pricing iteration that allows to update not only the upper bound (l.~\ref{l:ub}), but also the lower bound (l.~\ref{l:lb1}-\ref{l:lb2}). After the \gls{acr:cg} has converged or the time limit has been reached (l.~\ref{l:stopppingcriteria}), the algorithm branches on the obtained continuous solution in order to enforce integer feasible solutions (l.~\ref{l:bac}). In the following paragraphs, we detail each algorithmic component.

\noindent \textbf{Restricted master problem: }The \gls{acr:cg} procedure in our algorithm solves the continuous relaxation of Problem~\ref{problem:path-based-formulation}. However, in Problem~\ref{problem:path-based-formulation} the number of feasible paths per request $\Request \in \SetOfRequests^{\Freight}$ is intractable even for medium-sized networks. Accordingly, we solve the \gls{acr:rmp} considering only a subset of these paths, i.e., a subset of $\Tilde{\SetOfPaths}(\Request) \subseteq \SetOfPaths(\Request), \medspace \Request \in \SetOfRequests^{\Freight}$. Therefore, we initialize the \gls{acr:rmp} such that the problem is feasible in the very first iteration of the \gls{acr:cg}. Specifically, we initialize $\Tilde{\SetOfPaths}(\Request):= \{ \langle (\Origin^{\Request}, \RequestStart^{\Request}), (\Destination^{\Request}, \RequestEnd^{\Request}) \rangle \}, \Request \in \SetOfRequests^{\Freight}$. Thus, in the initial solution to the \gls{acr:rmp} all freight requests are rejected, i.e., sent via their dummy arcs. Then, we add additional variables $\FreightFlowPath^{\Request}_{\Path}, \Request \in \SetOfRequests^{\Freight}, \Path \in \Tilde{\SetOfPaths}(\Request)$ dynamically until the algorithm terminates. 
We solve the \gls{acr:rmp} with a standard commercial solver by warmstarting from the solution of the previous \gls{acr:cg} iteration. 

\noindent \textbf{Pricing problems: }The pricing problems identify the variables $\FreightFlowPath^{\Request}_{\Path}, \Request \in \SetOfRequests^{\Freight}, \Path \in \Tilde{\SetOfPaths}(\Request)$ that we add to the \gls{acr:rmp}. In every pricing problem, we identify a variable $\FreightFlowPath$ that represents the column yielding the maximum primal solution value improvement, and therefore comes with the smallest reduced cost for a given request. The exact solution to the pricing problem is required to obtain valid lower bounds and the reduced cost of a variable depends on the dual problem. To this end, the dual problem of the continuous relaxation of Problem~\ref{problem:path-based-formulation} is
\begin{maxi!}|s|[3]<b> 
{\DualGamma, \DualAlpha, \DualUpsilon, \DualNu, \DualDelta, \DualPi, \DualTau, \DualEta}
{\sum_{\Request \in \SetOfRequests^{\Passenger}} (\ServiceLevel \Demand^{\Request} \DualGamma  + \DualDelta^{\Request}) + \sum_{\Request \in \SetOfRequests^{\Freight}} \DualEta^{\Request} + \sum_{\Vehicle \in \SetOfVehicles} \NumUnits_{\Vehicle} \DualTau_{\Vehicle} \notag}
{\label{problem:dual}}
{} 
\breakObjective{+ \sum_{\Arc \in \SpannedArcs} \sum_{\Vehicle \in \SetOfVehicles} \VehicleArcRelation^{\Vehicle}_{\UnbracedArc} \UnitCapacity_{\Vehicle} \NumUnits_{\Vehicle} \DualUpsilon_{\UnbracedArc} + \sum_{\Arc \in \NotSpannedArcs} \sum_{\Vehicle \in \SetOfVehicles} \VehicleArcRelation^{\Vehicle}_{\UnbracedArc} \UnitCapacity_{\Vehicle} \NumUnits_{\Vehicle} \DualNu_{\UnbracedArc}  \label{obj:dual}}
\addConstraint{(\HTUVehicle) \quad \DualTau_{\Vehicle} - \sum_{\Arc \in \SetOfArcs^{\PathSegment}} \VehicleArcRelation^{\Vehicle}_{\UnbracedArc} \DualPi_{\UnbracedArc}}{\le \Cost_{\Vehicle},}{\hspace{-1.5cm} \forall \Vehicle \in \SetOfVehicles}
\addConstraint{(\FreightFlowPath) \quad \DualEta^{\Request} + \sum_{\Arc \in \SetOfArcs^{\PathSegment}} \Demand^{\Request} \PathArcRelation^{\Path}_{\UnbracedArc} \DualAlpha_{\UnbracedArc}}{\le \sum_{\Arc \in \SetOfArcs^{\PotentialFreightArcs}} \Demand^{\Request} \Cost_{\UnbracedArc} \PathArcRelation^{\Path}_{\UnbracedArc}, }{\hspace{-1.5cm} \forall \Request \in \SetOfRequests^{\Freight}, \Path \in \SetOfPaths(\Request) \label{prob:dual:cons:freight-flow}}
\addConstraint{(\PassengerFlow) \quad  \Demand^{\Request} \DualGamma + \DualDelta^{\Request} + \sum_{\Arc \in \SpannedArcs} \Demand^{\Request} \PathArcRelation^{\Path}_{\UnbracedArc} \DualUpsilon_{\UnbracedArc} + \sum_{\Arc \in \NotSpannedArcs} \Demand^{\Request} \PathArcRelation^{\Path}_{\UnbracedArc} \DualNu_{\UnbracedArc}}{\le 0,}{\hspace{-1.5cm} \forall \Request \in \SetOfRequests^{\Passenger}, \Path \in \SetOfPaths(\Request)}
\addConstraint{(\HTUOps) \quad \DualPi_{\UnbracedArc} - \DualAlpha_{\UnbracedArc} \sum_{\Vehicle \in \SetOfVehicles} \VehicleArcRelation^{\Vehicle}_{\UnbracedArc} \UnitCapacity_{\Vehicle} + \sum_{(i',j') \in \SpannedArcsFilter(\UnbracedArc)} \DualUpsilon_{i',j'} \sum_{\Vehicle \in \SetOfVehicles} \VehicleArcRelation^{\Vehicle}_{\UnbracedArc} \UnitCapacity_{\Vehicle}}{\le 0,}{\hspace{-1.5cm} \forall \Arc \in \SetOfArcs^{\PathSegment}}
\addConstraint{\quad \quad \DualGamma \ge 0; \quad \DualAlpha, \DualUpsilon, \DualNu, \DualDelta, \DualPi, \DualTau \le 0; \quad \DualEta \text{ free} }{}{}
\end{maxi!}
where $\SpannedArcsFilter(\UnbracedArc):= \{ (i^{\prime}, j^{\prime}) \in \SetOfArcs^{\VehicleLayer}: \SpanningMapping(i^{\prime}, j^{\prime}) = \Arc \}$ denotes the set of all vehicle layer arcs that are contracted into the given freight path segment $\Arc \in \SetOfArcs^{\PathSegment}$. Then, we obtain the respective reduced cost by re-arranging Constraints~\ref{prob:dual:cons:freight-flow}:

\begin{equation} \label{eq:reduced-cost}
    \bar{\Cost}^{\Request}_{\Path} = \Demand^{\Request} \Bigl[ \sum_{\Arc \in \SetOfArcs^\PotentialFreightArcs} \PathArcRelation^{\Path}_{\UnbracedArc} \Cost_{\UnbracedArc} - \sum_{\Arc \in \SetOfArcs^{\PathSegment}} \PathArcRelation^{\Path}_{\UnbracedArc} \DualAlpha_{\UnbracedArc} \Bigl] - \DualEta^{\Request}, \quad \forall \Request \in \SetOfRequests^{\Freight}, \Path \in \SetOfPaths(\Request)
\end{equation}

and the corresponding pricing problems for every $\Request \in \SetOfRequests^{\Freight}$ are independent and read 

\begin{mini!}|s|[0]<b> 
    {\FreightFlow}  
    { \Demand^{\Request} \Bigl[ \sum_{\Arc \in \SetOfArcs^{\PotentialFreightArcs}} \Cost_{\UnbracedArc} \FreightFlow^{\Request}_{\UnbracedArc} - \sum_{\Arc \in \SetOfArcs^{\PathSegment}} \DualAlpha_{\UnbracedArc} \FreightFlow^{\Request}_{\UnbracedArc} \Bigl] - \DualEta^{\Request} \label{obj:func-pricing}} 
    {\label{problem:pricing-problem}} 
    {} 
    \addConstraint{\sum_{j \in \Neighbors^{+}(i)} f^{\Request}_{\UnbracedArc} - \sum_{j \in \Neighbors^{-}(i)} \FreightFlow^{\Request}_{j,i}}{= \VertexDemand^{\Request}_{i}}{\forall  i \in \SetOfOrigins \cup \SetOfDestinations \cup \SetOfTerminals \label{cons:flow-conservation-pricing}} 
    \addConstraint{\FreightFlow^{\Request}_{\UnbracedArc}}{\in \{0,1\},\quad}{\forall \Arc \in \SetOfArcs^{\PotentialFreightArcs} \label{cons:domain-freight-flow-pricing}}
\end{mini!}

Solving the pricing problems~\ref{problem:pricing-problem} is equivalent to solving a \gls{acr:spp} with adapted arc costs $\Cost_{\UnbracedArc} - \DualAlpha_{\UnbracedArc}, \Arc \in \SetOfArcs^{\PathSegment}$ on the subgraph of $\Graph$ that is defined by the arc set $\SetOfArcs^{\PotentialFreightArcs}$ for every freight request. Note that the decomposed structure of the pricing problems allows their parallel computation. We can solve these \glspl{acr:spp} by standard approaches such as Dijkstra's algorithm, or, more efficiently, with the A* algorithm which is a label setting algorithm that prioritizes paths that are more likely to be optimal and thereby reduces unnecessary exploration. However, A* requires an admissible distance approximation which we can find from computing \glspl{acr:spp} on the unexpanded network, i.e., by discarding the problem's time dimension as follows: 
First, we compute lower bounds to the total cost on the minimum cost path between all pairs of \glspl{acr:ft} in the unexpanded network. We denote the resulting cost mapping as a function \mbox{$\DistanceApproximator^{\prime}: \SetOfFlatTerminals \times \SetOfFlatTerminals \rightarrow [0, \infty)$}. Due to the static nature of this cost mapping, we can compute it once in an offline fashion and use it in all pricing iterations. 
Finally, we complete the computation of the distance approximation $\DistanceApproximator$ by finding the connection of the request destination representations that yields the lowest approximation of total path cost. Specifically, we determine 
\begin{equation*}
    \DistanceApproximator(i, j) = \min_{i^{\prime} \in \Neighbors^{-}(j)} \DistanceApproximator^{\prime}(\MapToFlat(i), \MapToFlat(i^{\prime})) + \Cost_{i^{\prime}, j} , \medspace \forall \medspace i \in \SetOfTerminals, \medspace j \in \SetOfDestinations
\end{equation*}
where $\MapToFlat$ maps temporally expanded vertices to their unexpanded representation, i.e., $\MapToFlat(\Vertex) = \Stop, \medspace \forall \medspace \Vertex = (\Stop, \cdot, \cdot) \in \SetOfStops$. This gives an admissible distance approximation, which allows us to use the A* algorithm to speed up \gls{acr:spp} computations. We refer to \cite{LienkampSchiffer2022} for more details. 

\noindent \textbf{Partial pricing: }We apply partial pricing and thus, do not solve all pricing problems in every \gls{acr:cg} iteration. Instead, we only solve a subset of pricing problems in order to reduce computation time and promote heterogeneity in the generated column, which are more likely to be jointly selected in an optimal solution to Problem~\ref{problem:path-based-formulation}.  
Formally, we set the pricing strength that determines the maximum number of columns to add per iteration to $\PricingStrength \le 1$. After solving all pricing problems in the first \gls{acr:cg} iteration, we pop freight requests from the priority queue and solve the corresponding pricing problems until $\PricingStrength \medspace |\SetOfRequests^{\Freight}|$ variables with negative reduced cost have been found or all pricing problems have been solved. Here, we maintain the order of requests in the priority queue across pricing iterations. 
However, we regularly perform a full pricing iteration as suggested in \cite{KleinSchiffer2023} because partial pricing impedes the computation of tight lower bounds \citep[for a general introduction, see, e.g.,][]{UchoaPessoaEtAl2024}. Furthermore, we fall back to full pricing iterations if \gls{acr:cg} convergence slows down and the optimality gap has not improved beyond some threshold for multiple consecutive iterations. Specifically, we compute the average optimality gap reduction over the last $5$ iterations and conduct a full pricing iteration if this average reduction is below $0.0001$.

\noindent \textbf{Branching: }After solving the continuous relaxation of Problem~\ref{problem:path-based-formulation} via \gls{acr:cg} we obtain a solution that is potentially fractional in $\HTUVehicle, \HTUOps$, and $\FreightFlowPath$. 
We fix the set of variables $\FreightFlowPath$ to those that are in the current \gls{acr:rmp} and utilize a commercial solver's state-of-the-art branch-and-cut implementation to derive an integer feasible solution to Problem~\ref{problem:path-based-formulation}. By relying on this simple approach, we not only avoid initializing a second model, but effectively provide the solution to the root node of the \gls{acr:bab} by starting the commercial solver at the potentially fractional solution obtained from the \gls{acr:cg}.
Since we generated only a subset of all feasible paths --- specifically those required to solve the continuous relaxation at the root node --- integer solutions that we obtain while branching remain upper bounds to Problem~\ref{MIP}. 
\subsection{Branch-and-price}\label{sec:branch-and-price}
\noindent As an alternative to our \gls{acr:pab} algorithm, we provide a basic \gls{acr:bap} approach outlined as shown in Algorithm~\ref{algo:branch-and-price} to find optimal solutions to Problem~\ref{problem:path-based-formulation}. Contrary to the \gls{acr:pab} (cf. Algorithm~\ref{algo:price-and-branch}), the \gls{acr:bap} algorithm allows to continuously improve the lower bound by iteratively switching between branching and pricing until optimality is proven or a time limit is reached.

\begin{figure}[!bt] 
    \vspace{-15pt}
    \begin{minipage}{0.98\linewidth} 
        \begin{algorithm}[H]
            \caption{\textbf{Branch-and-price}} \label{algo:branch-and-price}
            \footnotesize
            \begin{algorithmic}[1]
                \Require Path-based formulation \eqref{problem:path-based-formulation}
                \State $\texttt{GlobalLB}, \texttt{GlobalUB}, \texttt{Incumbent} \gets 0, \infty, \texttt{None}$ 
                \State $\texttt{ActiveNodeQueue} \gets \texttt{InitializeNodeQueue}(\texttt{ContinuousRelaxation})$ \label{l:nodequeue}
                \While{$\texttt{OptimalityGap} > \epsilon$ and $\texttt{SolveTime} < \texttt{TimeLimit}$}
                    \State $\texttt{node} \gets \texttt{ActiveNodeQueue}.pop()$ \label{l:nodepop} \Comment{Sorted by parent's lower bound} 
                    \State $\texttt{ContinuousSolution} \gets \texttt{ColumnGeneration(node)}$ \label{l:cg}
                    \If {$\texttt{ContinuousSolution}.value \ge \texttt{GlobalUB}$} \label{l:prune1}
                        \State $\text{\textbf{continue}}$ \label{l:prune2} \Comment{Prune search tree}
                    \EndIf
                    \If {\textbf{not} $\texttt{IsFractional(ContinuousSolution)}$}
                        \State $\texttt{IntegerSolution} \gets \texttt{ContinuousSolution}$
                    \Else 
                        \State $\texttt{IntegerSolution} \gets \texttt{NodeUpperBound(ContinuousSolution)}$ \label{l:nodeub} \Comment{By standard solver} 
                        \State $\texttt{Left, Right} \gets  \texttt{BranchingRule(ContinuousSolution)}$ \label{l:branching} \Comment{cf. Equation~\eqref{eq:branching-rule}}
                        \State $\texttt{ActiveNodeQueue}.push(\texttt{Left, Right})$ \label{l:pushnode}
                    \EndIf
                    \If {$\texttt{IntegerSolution}.value < \texttt{GlobalUB}$} \label{l:updateglobalub1}
                        \State $\texttt{Incumbent} \gets \texttt{IntegerSolution}$ \label{l:updateglobalub2}
                        \State $\texttt{GlobalUB} \gets \texttt{Incumbent}.value$ \label{l:updateglobalub3}
                    \EndIf
                    \State $\texttt{GlobalLB} \gets \texttt{UpdateLB(ContinuousSolution, ActiveNodeQueue)}$ \label{l:updategloballb}\Comment{Min. across active nodes}
                \EndWhile \\
                \Return $\texttt{Incumbent}$
            \end{algorithmic}
        \end{algorithm}
    \end{minipage}
\end{figure}

Algorithm~\ref{algo:branch-and-price} initializes a queue of active nodes with the continuous relaxation of Problem~\ref{problem:path-based-formulation} (l.~\ref{l:nodequeue}). Then, in every iteration, we pop the first node from this queue and solve its continuous relaxation by \gls{acr:cg} (l.~\ref{l:nodepop}-\ref{l:cg}). To this end, we apply the same \gls{acr:cg} with partial pricing \revone{to $\epsilon$-precise solutions}. After solving the node, the algorithm first checks if it can prune the search tree based on the obtained solution (l.~\ref{l:prune1}-\ref{l:prune2}). Otherwise, if the obtained solution is fractional in any design variable $\HTUVehicle$, the algorithm derives a node-based upper bound (l.~\ref{l:nodeub}). Furthermore, it identifies the design variable to branch based on Equation~\eqref{eq:branching-rule}, creates two new child nodes  (l.~\ref{l:branching}), and adds them to the queue of active nodes (l.~\ref{l:pushnode}). Two comments are in order. First, we initialize new nodes with a node-based lower bound that equals their parent's lower bound and sort them in the queue accordingly in ascending order. Second, new nodes can be infeasible. We detect infeasibility after popping the node from the queue when solving the \gls{acr:rmp} in the \gls{acr:cg} procedure (l.~\ref{l:cg}) and proceed with the next iteration in this case. If the node-based upper bound derived in Line~\ref{l:nodeub}, yields an improvement compared to the solution value of the current incumbent, we update the incumbent and the global upper bound accordingly (l.~\ref{l:updateglobalub1}-\ref{l:updateglobalub3}). Finally, we update the global lower bound (l.~\ref{l:updategloballb}) which is given by the minimum node-based lower bounds of all active nodes. In the following, we briefly explain our \gls{acr:bap} approach concerning the applied branching rule and the derivation of upper bounds. 

\noindent \textbf{Branching strategy}: We branch based on the design variables $\HTUVehicle_{\Vehicle}, \medspace \Vehicle \in \SetOfVehicles$ and select the design variable with the highest fractional value in the final solution of the \gls{acr:cg} applied to solve the continuous relaxation of the respective search tree node. More specifically, after solving a node in the search tree, we determine
\begin{equation}\label{eq:branching-rule}
    \argmax_{\HTUVehicle_{\Vehicle}, \Vehicle \in \SetOfVehicles} \hspace{0.25cm} \min \{ \HTUVehicle_{\Vehicle} - \lfloor \HTUVehicle_{\Vehicle} \rfloor, 1 - (\HTUVehicle_{\Vehicle} - \lfloor \HTUVehicle_{\Vehicle} \rfloor ) \}
\end{equation}
and add an inequality to the \gls{acr:rmp} that reflects the branching. This branching on the design variables does not change the pricing problem. \revone{The applied branching rule does not guarantee integer variables $\HTUOps$, and $\FreightFlowPath$, and we rely on the node-based upper bounds to ensure integrality across these variables. In this context, we note that branching on variables $\FreightFlowPath$ could lead to large branching trees and additionally changes the pricing problem. Thus, the implementation of such a branching rule requires further adjustment of the algorithmic framework. Moreover, w}e populate columns found in any search tree node to all currently active nodes and do not prune the set of columns. 

\noindent\textbf{Node-based upper bounds}: Following the observation that commercial solvers are capable of branching on a continuous solution of Problem~\ref{problem:path-based-formulation} effectively, we assign a restrictive time limit to the commercial solver (e.g., $60$ seconds) and utilize its branch-and-cut algorithm to derive upper bounds in every node of the search tree.

\section{Experimental design} \label{sec:case-study}
\begin{figure}[!b]
    \centering
    \includegraphics[width=\textwidth]{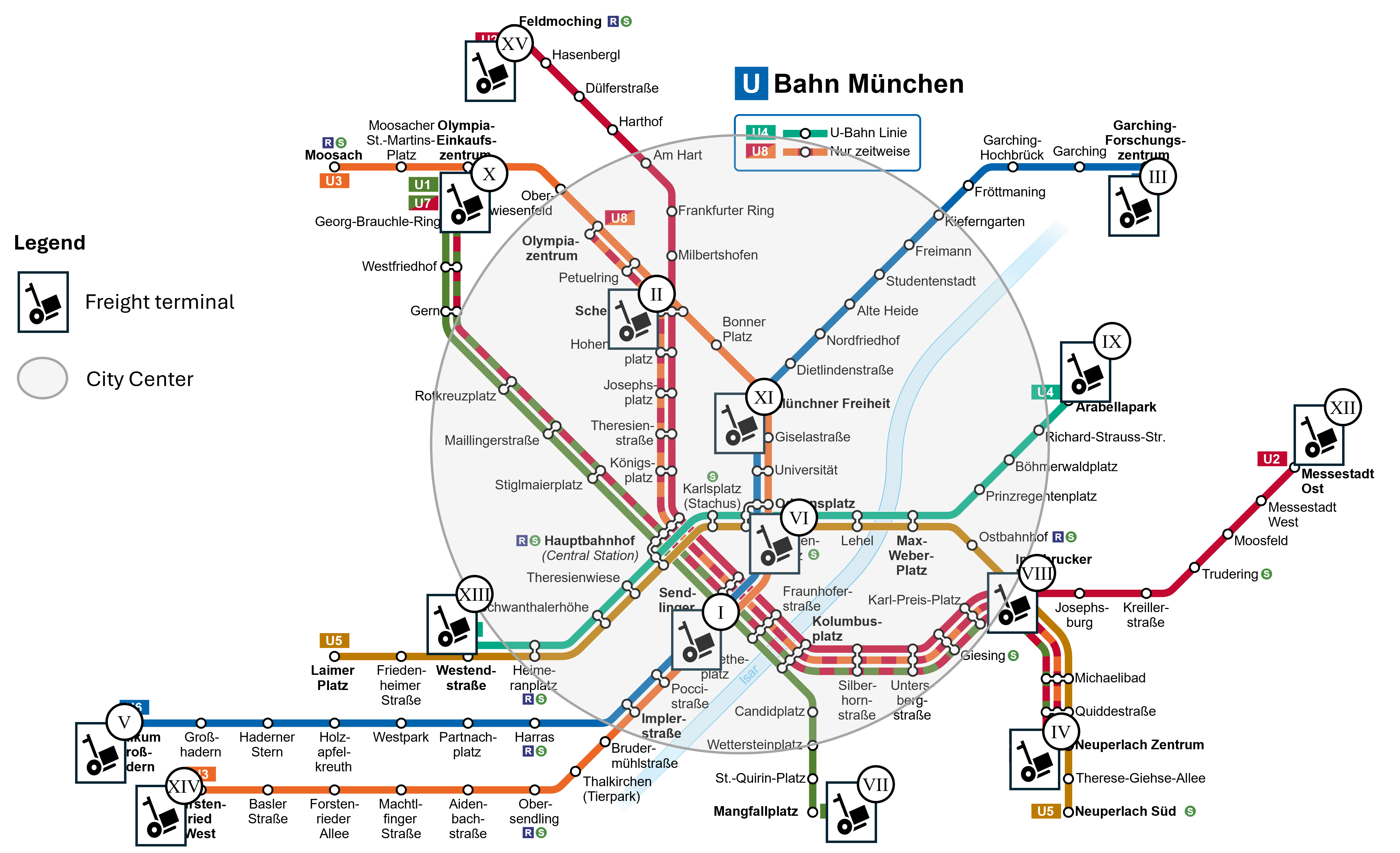}
    \caption{Munich subway network \citep{Netzplan2020} with assumed \gls{acr:ft} layout}
    \label{fig:munich-network}
\end{figure}

\noindent We test our algorithm on the subway network of Munich (cf. Figure~\ref{fig:munich-network}) that we enrich by an assumption about the \gls{acr:ft} layout, and generate instances based on publicly available data where possible. 

In the following, we describe our case study and summarize the sources used and the resulting parameters in Table~\ref{tab:exp_params}. Following the General Transit Feed Specification (GTFS) for the Munich subway network, trip data includes physical locations and timetable information. We concatenate trips from a representative day to reach a reasonable assumption on \revone{85} \gls{acr:pt} vehicles' routes through the network during the considered time period. Specifically, we assume that two trips are performed by the same vehicle if they end and start at the same \gls{acr:pt} stop, and no other trip starts or ends at the same stop in the time between the two considered vehicles. Based on information by the operator \citep{MVG2023}, we set heterogeneous vehicle capacities from $\{870, 910, 936 \}$ by randomly assigning vehicle types according to the probabilities $\{0.52, 0.13, 0.35 \}$ that reflect the current vehicle type distribution in Munich's \gls{acr:pts}.

We sample passenger requests between 6 a.m. and 11 a.m. based on the urban travel demand simulation tool MITO \citep[cf.][]{MoeckelKuehnelEtAl2020}. Here, we scale the demand per passenger request to $\Demand^{\Request}=24.71, \medspace \Request \in \SetOfRequests^{\Passenger}$ such that the complete set of $10{,}000$ sampled requests represents the demand in the evaluated time period, which we derive based on the operator's data for an entire year \citep[cf.][]{MVG2023}. We \revone{require a service level of $\ServiceLevel = 0.999$ and} pre-compute passengers' potential paths as described in Section~\ref{subsection:graph-construction}. 

Real-world data on urban freight shipments is notoriously hard to get. 
As an alternative, we sample freight request destination locations based on the population distribution and income per capita distribution per city district and assume that every freight request originates in one of $20$ randomly located \gls{acr:lsp} distribution centers in the city's outskirts. In this context, we assume that the city's outskirts lie within a radius of $8$ -- $10$ km of the city center. We set the homogeneous demand of freight requests such that the accumulated demand reflects the daily parcel delivery volume in Munich derived based on volume per capita (cf. Table~\ref{tab:exp_params}). To convert the resulting freight demand (in number of parcels) to passenger equivalents, we convert $12$ parcels to $1$ passenger equivalent. This conversion is based on the required space per passenger, the dimensions of a trolley, the parcel volume, and a trolley utilization of $75$\%. Furthermore, we connect every origin and destination representation to $\NumConnections = 1$ \gls{acr:ft} representation as outlined in Section~\ref{app:implementation}. 
Finally, we parameterize costs as follows. \revone{First, we derive the daily value of the total design cost per \gls{acr:htu} for transit vehicles comprising six cars based on $25$ years of usage. Specifically, $\Cost_{\Vehicle}=\frac{\num{1.51e6}}{25 \times 365}=165.48 \text{ \euro{}}, \medspace \Vehicle \in \SetOfVehicles$. We increase the design cost proportionally if the public transit vehicle consists of less cars}. 

\begin{table}[!b]
    \centering
    \scriptsize
    \caption{Case Study Parameters}
    \label{tab:exp_params}
    \begin{tabular}{llll}
        \hline
        \textbf{Parameter} & \textbf{Unit} & \textbf{Value} & \textbf{Source} \\
        \hline
        Subway passenger demand & passengers/year & $353$ Mio & \cite{MVG2023}  \\
        German parcel volume & parcels/year & $4220$ Mio & \cite{Bundesnetzagentur2022}  \\
        Design cost per \gls{acr:htu} & \euro{}/\gls{acr:htu} & $1.515$ Mio & \cite{MVG2020}  \\
        Externality cost (truck) & \euro{}/vehicle \& km   & $[0.05, 11.71]$ & \cite{DeLanghe2017}   \\
        Externality cost (cargo bike \footnote{We assume similar externality cost as for electric mopeds}) & \euro{}/vehicle \& km   & $0.115$ & \cite{SchroederKirnEtAl2023}  \\
        Subway capacity & Passengers/vehicle & $\{ 870, 912, 936\}$ &  \cite{MVG2023}  \\
        Passenger space requirement & $\text{m}^{2}$ & $0.25$ & \cite{VDV} \\
        Duration of fleet usage & years & 25 & \cite{AfA} \\
        Base rate & $\%$ & $3.62$ & \cite{Bundesbank} \\
        Truck typical tour length & km  & $80$ & \cite{OW2021}  \\
        Truck capacity & parcels/vehicle   & $100$ & \cite{OW2021}  \\
        Cargo-Bike typical tour length & km & $12.2$ & \cite{KoingConwayEtAl2016} \\
        Cargo-Bike capacity & parcels/vehicle & $20$ & \cite{LlorcaMoeckel2021}  \\
        Parcel volume & $\text{m}^{3}$ & $0.027$ & \cite{DHL2024} \\
        Trolley dimensions (H x W x D) & m & $1.8\times1.2\times0.8$ & \cite{Wanzl2024} \\
        Working days per year & days & 255 & - \\
        \hline
    \end{tabular}
\end{table}

Let the externality cost of conventional truck-based delivery be \euro{0.2} per vehicle and kilometer. Furthermore, we assume a truck tour length of $80$ km and a delivery capacity of $100$ parcels per tour. According to the conversion factor that we assume, a unit of demand equals $12$ parcels. Thus, we set $\Cost^{\Request}_{\Penalty} = \frac{0.2 \times 80 \times \Demand_{\Request} \times 12}{100}, \medspace \Request \in \SetOfRequests^{\Freight}$  Similarly, let the externality cost of cargo-bike delivery be \euro{0.115} per vehicle and kilometer, the average tour length be $12.2$ km, and the delivery capacity be $20$ parcels per tour. Then, we set $\Cost_{\UnbracedArc} = \frac{0.115 \times 12.2 \times 12}{20}, \medspace \Arc \in \SetOfArcs^{\Egress}$. 

Third, we chose the routing cost $\Cost_{\UnbracedArc} = 0.0406 \times d_{\UnbracedArc}, \medspace \Arc \in \SetOfArcs^{\VehicleLayer}$ proportionally to the kilometers of distance $d_{\UnbracedArc}$ between the stops that $i$ and $j$ represent, scaled by the externality cost of the transported freight. 
% instance variations
All instances share the \gls{acr:pts} network and passenger demand. The sizes of the instances differ by the number of freight requests we consider, and we generate $n=15$ experiments with different seeds for every instance size. Specifically, we generate instances of the sizes $|\SetOfRequests^{\Freight}| \in [250, 500, 1000, 2000, 3000]$. 

All experiments have been conducted single-threaded on a standard desktop computer with an Intel Core i9-9900, 3.1 GHz CPU, and 16 GB of RAM, running Ubuntu 20.04. We have used Python 3.10.2 with CPLEX 22.1 to solve the \gls{acr:rmp} in the \gls{acr:cg} and perform the subsequent branching. We used the DOcplex library as a modelling interface, allow CPLEX to use its presolve capabilities, and configure CPLEX to scale the coefficient matrix aggressively in the \gls{acr:rmp} of the \gls{acr:pab} algorithm. We have run all experiments with a maximum runtime of $90$ minutes. In our \gls{acr:pab} algorithm, we reserve $15$ minutes for the branching and stop the \gls{acr:cg} otherwise at an optimality tolerance of $\epsilon = 0.001$.

\section{Results} \label{subsection:comp-results}
\noindent In the following Section~\ref{sec:comp-result}, we show the efficiency of our algorithmic framework. In this context, we determine the value of partial pricing by providing results with partial pricing of varying degree, i.e., varying number of pricing problems solved per iteration. We show that partial pricing decreases the required number of pricing iterations. Furthermore, we compare our \gls{acr:pab} algorithm to a \gls{acr:mip} and show that we increase the solvable instance size significantly. Moreover, we provide results of the presented \gls{acr:bap} algorithm. In Section~\ref{subsection:managerial-insights}, we present a sensitivity analysis on unknown cost factors, and show that our framework successfully increases the utilization of Munich's \gls{acr:pts} during off-peak hours as well as its capability to allocate capacity in a dynamic fashion respecting passenger demand peaks and scheduling freight transportation around those peaks. 

\subsection{Computational results}\label{sec:comp-result}
\noindent We run our \gls{acr:pab} approach with partial pricing and show the value of partial pricing in Table~\ref{tab:partial-pricing} by comparing different pricing strength parameters $\PricingStrength \in \{0.1, 0.2, 0.3, 0.4, 0.5, 1.0\}$. 

\begin{table}[!b]
    \centering
    \caption{The value of partial pricing ($n=15$)}
    \footnotesize
    \begin{tabular}{lcccccc}
        \toprule
        \textbf{Instance Sizes} & $\PricingStrength=0.1$ & $\PricingStrength=0.2$ & $\PricingStrength=0.3$ & $\PricingStrength=0.4$ & $\PricingStrength=0.5$ & $\PricingStrength=1.0$ \\
        \midrule
        \midrule
        \multicolumn{7}{l}{\textbf{Median number of variables added per request}} \\
        \midrule
        250   & 30.84 & 31.20 & 36.00 & 42.17 & 45.44 & 54.66 \\
        500   & 22.33 & 25.36 & 29.04 & 35.84 & 42.02 & 50.27 \\
        1,000 & 19.20 & 22.01 & 26.92 & 31.64 & 35.99 & 47.87 \\
        2,000 & 18.49 & 23.66 & 28.61 & 35.18 & 38.63 & 50.41 \\
        3,000 & 18.83 & 23.85 & 28.13 & 34.45 & 39.10 & 50.12 \\
        \midrule
        \multicolumn{7}{l}{\textbf{Median integrality gap}} \\
        \midrule
        250   & 1.56\% & 1.45\% & 1.35\% & 1.43\% & 1.20\% & 1.15\% \\
        500   & 0.93\% & 0.87\% & 0.87\% & 0.83\% & 0.79\% & 0.75\% \\
        1,000 & 1.06\% & 0.96\% & 0.97\% & 0.87\% & 0.79\% & 0.93\% \\
        2,000 & 1.19\% & 1.17\% & 1.04\% & 0.99\% & 0.93\% & 0.90\% \\
        3,000 & 0.97\% & 0.98\% & 0.96\% & 0.97\% & 1.10\% & 1.72\% \\
        \bottomrule
    \end{tabular}
    \label{tab:partial-pricing}
\end{table}

As can be seen, solving only $10\%$ of all pricing problems (i.e., $\PricingStrength=0.1$) yields solutions with a median integrality gap below $1.56\%$. At the same time, partial pricing with $\PricingStrength=0.1$ saves the creation of every second path variable compared to an approach with full pricing, i.e., $\PricingStrength=1.0$. We observe results that are similar in quality but require our algorithm to generate significantly less path variables. This result indicates that partial pricing indeed leads to more heterogeneous columns that \revone{utilize different graph components and} are more likely to be jointly selected in an integer solution. We display the distribution of integrality gaps across different instance sizes and pricing strengths in Appendix~\ref{app:partial-pricing}.

\vspace{2pt}
\begin{result}
    Partial pricing increases the heterogeneity of paths such that the median number of generated columns decreases from $54.66$ to $30.84$ in small instances and from $50.12$ to $18.83$ in large instances while the median integrality gaps are less than $1.56\%$ in all instances. 
\end{result}
\vspace{2pt}

Table~\ref{tab:performance} compares the \gls{acr:pab} with the \gls{acr:mip}~\ref{MIP} and Figure~\ref{fig:benchmark-boxplots} provides a complementary visualization of the reported integrality gaps and solve times. 

\begin{table}[t]
  \centering
  \fontsize{8.5}{9.5}\selectfont
  \caption{Benchmark results ($n=15$). \revone{We compute gaps relative to the maximum of the two lower bounds.}}
  \begin{tabular}{ccccccc}
    \toprule
    \multicolumn{1}{c}{\textbf{Instance}} & \multicolumn{2}{c}{\textbf{Median integrality}} & \multicolumn{2}{c}{\textbf{Median solve time until}} & \multicolumn{2}{c}{\textbf{Solved}} \\
    \multicolumn{1}{c}{\textbf{size}} & \multicolumn{2}{c}{\textbf{gap [\%]}} & \multicolumn{2}{c}{\textbf{first feasible solution [s]}} & \multicolumn{2}{c}{\textbf{instances}} \\
    \cmidrule(lr){2-3} \cmidrule(lr){4-5} \cmidrule(lr){6-7}
          & \textbf{MIP} & \textbf{P\&B} & \textbf{MIP} & \textbf{P\&B} & \textbf{MIP} & \textbf{P\&B} \\
    \midrule
    \midrule
    250   & 1.41\% & 1.33\%        & 2905.84 & 486.06       & 15 & 15 \\
    500   & 1.83\% & 0.93\%        & 3445.36 & 568.74       & 13 & 15 \\
    1,000  & -       & 1.06\%       & -       & 696.30       & 0    & 15 \\
    2,000  & -       & 1.19\%       & -       & 1152.22       & 0    & 15 \\
    3,000  & -       & 0.97\%       & -       & 1800.69       & 0    & 15 \\
    \bottomrule
  \end{tabular}%
  \label{tab:performance}%
\end{table}

In this setting, the commercial solver benefits from the presented preprocessing techniques and the implemented graph pruning to the same extent as our \gls{acr:pab} algorithm. Solving \gls{acr:mip}~\ref{MIP} in a large-scale setting is time-consuming and runs into memory bounds quickly. In particular, while the commercial solver provides a solution with a median integrality gap of less than $2\%$ for all instances of size $250$ freight requests and most instances of size $500$ freight requests, it consistently runs into memory bounds when extending the setting to a larger scale. 

\begin{figure}[!b]
  \begin{subfigure}[t]{0.48\textwidth}
    \centering
    % This file was created with tikzplotlib v0.10.1.
\begin{tikzpicture}

\definecolor{darkgray176}{RGB}{176,176,176}
\definecolor{darkgreen0640}{RGB}{80,200,120}
\definecolor{lightgray204}{RGB}{204,204,204}

\definecolor{gray}{RGB}{128,128,128}
\definecolor{skyblue147196222}{RGB}{147,196,222}
\definecolor{blue}{RGB}{8,73,145}

\begin{axis}[
legend cell align={left},
legend style={fill opacity=0.8, draw opacity=1, text opacity=1, draw=lightgray204},
tick align=outside,
tick pos=left,
unbounded coords=jump,
x grid style={darkgray176},
xlabel={Instance size [\# Requests]},
xmin=-0.5, xmax=6,
xtick style={color=black},
xtick={0.15,1.45,2.75,4.05,5.35},
xticklabels={250,500,1000,2000,3000},
y grid style={darkgray176},
ylabel={Time to first feas. solution [s]},
ymajorgrids,
ymin=130.337016761303, ymax=5466.25222307444,
ytick style={color=black},
width=0.75\linewidth,
height=0.62\linewidth,
scale only axis,
trim axis left, trim axis right,
]
\path [draw=black, fill=blue]
(axis cs:-0.15,450.249617695808)
--(axis cs:0.15,450.249617695808)
--(axis cs:0.15,494.189997792244)
--(axis cs:-0.15,494.189997792244)
--(axis cs:-0.15,450.249617695808)
--cycle;
\addplot [black, forget plot]
table {%
0 450.249617695808
0 427.268850803375
};
\addplot [black, forget plot]
table {%
0 494.189997792244
0 522.536777257919
};
\addplot [black, forget plot]
table {%
-0.075 427.268850803375
0.075 427.268850803375
};
\addplot [black, forget plot]
table {%
-0.075 522.536777257919
0.075 522.536777257919
};
\addplot [black, mark=o, mark size=3, mark options={solid,fill opacity=0}, only marks, forget plot]
table {%
0 372.878617048264
};
\path [draw=black, fill=darkgreen0640]
(axis cs:0.15,2673.96479475498)
--(axis cs:0.45,2673.96479475498)
--(axis cs:0.45,3067.72981226444)
--(axis cs:0.15,3067.72981226444)
--(axis cs:0.15,2673.96479475498)
--cycle;
\addplot [black, forget plot]
table {%
0.3 2673.96479475498
0.3 2245.32673215866
};
\addplot [black, forget plot]
table {%
0.3 3067.72981226444
0.3 3215.56966304779
};
\addplot [black, forget plot]
table {%
0.225 2245.32673215866
0.375 2245.32673215866
};
\addplot [black, forget plot]
table {%
0.225 3215.56966304779
0.375 3215.56966304779
};
\addplot [black, mark=o, mark size=3, mark options={solid,fill opacity=0}, only marks, forget plot]
table {%
0.3 3995.7768175602
};
\path [draw=black, fill=blue]
(axis cs:1.15,505.218910098076)
--(axis cs:1.45,505.218910098076)
--(axis cs:1.45,597.15109705925)
--(axis cs:1.15,597.15109705925)
--(axis cs:1.15,505.218910098076)
--cycle;
\addplot [black, forget plot]
table {%
1.3 505.218910098076
1.3 453.299555778503
};
\addplot [black, forget plot]
table {%
1.3 597.15109705925
1.3 699.328678607941
};
\addplot [black, forget plot]
table {%
1.225 453.299555778503
1.375 453.299555778503
};
\addplot [black, forget plot]
table {%
1.225 699.328678607941
1.375 699.328678607941
};
\addplot [black, mark=o, mark size=3, mark options={solid,fill opacity=0}, only marks, forget plot]
table {%
1.3 736.800552129745
};
\path [draw=black, fill=darkgreen0640]
(axis cs:1.45,3045.12731647492)
--(axis cs:1.75,3045.12731647492)
--(axis cs:1.75,4266.31641817093)
--(axis cs:1.45,4266.31641817093)
--(axis cs:1.45,3045.12731647492)
--cycle;
\addplot [black, forget plot]
table {%
1.6 3045.12731647492
1.6 1862.83612799644
};
\addplot [black, forget plot]
table {%
1.6 4266.31641817093
1.6 5223.71062278748
};
\addplot [black, forget plot]
table {%
1.525 1862.83612799644
1.675 1862.83612799644
};
\addplot [black, forget plot]
table {%
1.525 5223.71062278748
1.675 5223.71062278748
};
\path [draw=black, fill=blue]
(axis cs:2.45,654.209554791451)
--(axis cs:2.75,654.209554791451)
--(axis cs:2.75,765.8266248703)
--(axis cs:2.45,765.8266248703)
--(axis cs:2.45,654.209554791451)
--cycle;
\addplot [black, forget plot]
table {%
2.6 654.209554791451
2.6 599.73455119133
};
\addplot [black, forget plot]
table {%
2.6 765.8266248703
2.6 782.603931427002
};
\addplot [black, forget plot]
table {%
2.525 599.73455119133
2.675 599.73455119133
};
\addplot [black, forget plot]
table {%
2.525 782.603931427002
2.675 782.603931427002
};
\addplot [black, mark=o, mark size=3, mark options={solid,fill opacity=0}, only marks, forget plot]
table {%
2.6 377.062116384506
2.6 1025.6811144352
2.6 976.587995052338
2.6 1095.29742050171
};
\path [draw=black, fill=darkgreen0640]
;
\addplot [black, forget plot]
table {%
2.9 nan
2.9 nan
};
\addplot [black, forget plot]
table {%
2.9 nan
2.9 nan
};
\addplot [black, forget plot]
table {%
2.825 nan
2.975 nan
};
\addplot [black, forget plot]
table {%
2.825 nan
2.975 nan
};
\path [draw=black, fill=blue]
(axis cs:3.75,1045.65884745121)
--(axis cs:4.05,1045.65884745121)
--(axis cs:4.05,1297.94282722473)
--(axis cs:3.75,1297.94282722473)
--(axis cs:3.75,1045.65884745121)
--cycle;
\addplot [black, forget plot]
table {%
3.9 1045.65884745121
3.9 884.196155071259
};
\addplot [black, forget plot]
table {%
3.9 1297.94282722473
3.9 1532.94029855728
};
\addplot [black, forget plot]
table {%
3.825 884.196155071259
3.975 884.196155071259
};
\addplot [black, forget plot]
table {%
3.825 1532.94029855728
3.975 1532.94029855728
};
\addplot [black, mark=o, mark size=3, mark options={solid,fill opacity=0}, only marks, forget plot]
table {%
3.9 616.802316904068
};
\path [draw=black, fill=darkgreen0640]
;
\addplot [black, forget plot]
table {%
4.2 nan
4.2 nan
};
\addplot [black, forget plot]
table {%
4.2 nan
4.2 nan
};
\addplot [black, forget plot]
table {%
4.125 nan
4.275 nan
};
\addplot [black, forget plot]
table {%
4.125 nan
4.275 nan
};
\path [draw=black, fill=blue]
(axis cs:5.05,1410.36974561214)
--(axis cs:5.35,1410.36974561214)
--(axis cs:5.35,1919.32543003559)
--(axis cs:5.05,1919.32543003559)
--(axis cs:5.05,1410.36974561214)
--cycle;
\addplot [black, forget plot]
table {%
5.2 1410.36974561214
5.2 760.431878805161
};
\addplot [black, forget plot]
table {%
5.2 1919.32543003559
5.2 2438.23425316811
};
\addplot [black, forget plot]
table {%
5.125 760.431878805161
5.275 760.431878805161
};
\addplot [black, forget plot]
table {%
5.125 2438.23425316811
5.275 2438.23425316811
};
\addplot [black, mark=o, mark size=3, mark options={solid,fill opacity=0}, only marks, forget plot]
table {%
5.2 2687.38654303551
};
\path [draw=black, fill=darkgreen0640]
;
\addplot [black, forget plot]
table {%
5.5 nan
5.5 nan
};
\addplot [black, forget plot]
table {%
5.5 nan
5.5 nan
};
\addplot [black, forget plot]
table {%
5.425 nan
5.575 nan
};
\addplot [black, forget plot]
table {%
5.425 nan
5.575 nan
};
\addplot [black, forget plot]
table {%
-0.15 486.057582378388
0.15 486.057582378388
};
\addplot [black, forget plot]
table {%
0.15 2905.8402326107
0.45 2905.8402326107
};
\addplot [black, forget plot]
table {%
1.15 568.743363142014
1.45 568.743363142014
};
\addplot [black, forget plot]
table {%
1.45 3445.36229634285
1.75 3445.36229634285
};
\addplot [black, forget plot]
table {%
2.45 696.303322792053
2.75 696.303322792053
};
\addplot [black, forget plot]
table {%
2.75 nan
3.05 nan
};
\addplot [black, forget plot]
table {%
3.75 1152.21557402611
4.05 1152.21557402611
};
\addplot [black, forget plot]
table {%
4.05 nan
4.35 nan
};
\addplot [black, forget plot]
table {%
5.05 1800.68681788445
5.35 1800.68681788445
};
\addplot [black, forget plot]
table {%
5.35 nan
5.65 nan
};

% manual legend
\addlegendimage{area legend, draw=black, fill=blue}
\addlegendentry{P\&B}

\addlegendimage{area legend, draw=black, fill=darkgreen0640}
\addlegendentry{MIP}

\end{axis}

\end{tikzpicture}
    \vspace{-0.5cm}
    \caption{Solve time}
    \label{fig:pp-solve-time}
  \end{subfigure}
  \begin{subfigure}[t]{0.48\textwidth}
    \centering
    % This file was created with tikzplotlib v0.10.1.
\begin{tikzpicture}

\definecolor{darkgray176}{RGB}{176,176,176}
\definecolor{darkgreen0640}{RGB}{80,200,120}
\definecolor{lightgray204}{RGB}{204,204,204}

\definecolor{gray}{RGB}{128,128,128}
\definecolor{skyblue147196222}{RGB}{147,196,222}
\definecolor{blue}{RGB}{8,73,145}

\begin{axis}[
legend cell align={left},
legend style={fill opacity=0.8, draw opacity=1, text opacity=1, draw=lightgray204},
tick align=outside,
tick pos=left,
unbounded coords=jump,
x grid style={darkgray176},
xlabel={Instance size [\# Requests]},
xmin=-0.5, xmax=6,
xtick style={color=black},
xtick={0.15,1.45,2.75,4.05,5.35},
xticklabels={250,500,1000,2000,3000},
y grid style={darkgray176},
ylabel={Integrality gap [\%]},
ymajorgrids,
ymin=0, ymax=10,
ytick style={color=black},
width=0.75\linewidth,
height=0.62\linewidth,
scale only axis,
trim axis left, trim axis right,
]
\path [draw=black, fill=blue]
(axis cs:-0.15,1.08)
--(axis cs:0.15,1.08)
--(axis cs:0.15,1.49)
--(axis cs:-0.15,1.49)
--(axis cs:-0.15,1.08)
--cycle;
\addplot [black, forget plot]
table {%
0 1.08
0 0.84
};
\addplot [black, forget plot]
table {%
0 1.49
0 1.75
};
\addplot [black, forget plot]
table {%
-0.075 0.84
0.075 0.84
};
\addplot [black, forget plot]
table {%
-0.075 1.75
0.075 1.75
};
\path [draw=black, fill=darkgreen0640]
(axis cs:0.15,1.27)
--(axis cs:0.45,1.27)
--(axis cs:0.45,1.495)
--(axis cs:0.15,1.495)
--(axis cs:0.15,1.27)
--cycle;
\addplot [black, forget plot]
table {%
0.3 1.27
0.3 1.01
};
\addplot [black, forget plot]
table {%
0.3 1.495
0.3 1.54
};
\addplot [black, forget plot]
table {%
0.225 1.01
0.375 1.01
};
\addplot [black, forget plot]
table {%
0.225 1.54
0.375 1.54
};
\addplot [black, mark=o, mark size=3, mark options={solid,fill opacity=0}, only marks, forget plot]
table {%
0.3 2.02
};
\path [draw=black, fill=blue]
(axis cs:1.15,0.69)
--(axis cs:1.45,0.69)
--(axis cs:1.45,1.065)
--(axis cs:1.15,1.065)
--(axis cs:1.15,0.69)
--cycle;
\addplot [black, forget plot]
table {%
1.3 0.69
1.3 0.61
};
\addplot [black, forget plot]
table {%
1.3 1.065
1.3 1.36
};
\addplot [black, forget plot]
table {%
1.225 0.61
1.375 0.61
};
\addplot [black, forget plot]
table {%
1.225 1.36
1.375 1.36
};
\path [draw=black, fill=darkgreen0640]
(axis cs:1.45,1.37)
--(axis cs:1.75,1.37)
--(axis cs:1.75,2.26)
--(axis cs:1.45,2.26)
--(axis cs:1.45,1.37)
--cycle;
\addplot [black, forget plot]
table {%
1.6 1.37
1.6 0.94
};
\addplot [black, forget plot]
table {%
1.6 2.26
1.6 2.62
};
\addplot [black, forget plot]
table {%
1.525 0.94
1.675 0.94
};
\addplot [black, forget plot]
table {%
1.525 2.62
1.675 2.62
};
\addplot [black, mark=o, mark size=3, mark options={solid,fill opacity=0}, only marks, forget plot]
table {%
1.6 3.73
1.6 3.7
};
\path [draw=black, fill=blue]
(axis cs:2.45,0.795)
--(axis cs:2.75,0.795)
--(axis cs:2.75,1.22)
--(axis cs:2.45,1.22)
--(axis cs:2.45,0.795)
--cycle;
\addplot [black, forget plot]
table {%
2.6 0.795
2.6 0.71
};
\addplot [black, forget plot]
table {%
2.6 1.22
2.6 1.58
};
\addplot [black, forget plot]
table {%
2.525 0.71
2.675 0.71
};
\addplot [black, forget plot]
table {%
2.525 1.58
2.675 1.58
};
\addplot [black, mark=o, mark size=3, mark options={solid,fill opacity=0}, only marks, forget plot]
table {%
2.6 2.17
};
\path [draw=black, fill=darkgreen0640]
;
\addplot [black, forget plot]
table {%
2.9 nan
2.9 nan
};
\addplot [black, forget plot]
table {%
2.9 nan
2.9 nan
};
\addplot [black, forget plot]
table {%
2.825 nan
2.975 nan
};
\addplot [black, forget plot]
table {%
2.825 nan
2.975 nan
};
\path [draw=black, fill=blue]
(axis cs:3.75,0.855)
--(axis cs:4.05,0.855)
--(axis cs:4.05,1.34)
--(axis cs:3.75,1.34)
--(axis cs:3.75,0.855)
--cycle;
\addplot [black, forget plot]
table {%
3.9 0.855
3.9 0.58
};
\addplot [black, forget plot]
table {%
3.9 1.34
3.9 1.51
};
\addplot [black, forget plot]
table {%
3.825 0.58
3.975 0.58
};
\addplot [black, forget plot]
table {%
3.825 1.51
3.975 1.51
};
\addplot [black, mark=o, mark size=3, mark options={solid,fill opacity=0}, only marks, forget plot]
table {%
3.9 2.32
};
\path [draw=black, fill=darkgreen0640]
;
\addplot [black, forget plot]
table {%
4.2 nan
4.2 nan
};
\addplot [black, forget plot]
table {%
4.2 nan
4.2 nan
};
\addplot [black, forget plot]
table {%
4.125 nan
4.275 nan
};
\addplot [black, forget plot]
table {%
4.125 nan
4.275 nan
};
\path [draw=black, fill=blue]
(axis cs:5.05,0.815)
--(axis cs:5.35,0.815)
--(axis cs:5.35,1.18)
--(axis cs:5.05,1.18)
--(axis cs:5.05,0.815)
--cycle;
\addplot [black, forget plot]
table {%
5.2 0.815
5.2 0.46
};
\addplot [black, forget plot]
table {%
5.2 1.18
5.2 1.49
};
\addplot [black, forget plot]
table {%
5.125 0.46
5.275 0.46
};
\addplot [black, forget plot]
table {%
5.125 1.49
5.275 1.49
};
\path [draw=black, fill=darkgreen0640]
;
\addplot [black, forget plot]
table {%
5.5 nan
5.5 nan
};
\addplot [black, forget plot]
table {%
5.5 nan
5.5 nan
};
\addplot [black, forget plot]
table {%
5.425 nan
5.575 nan
};
\addplot [black, forget plot]
table {%
5.425 nan
5.575 nan
};
\addplot [black, forget plot]
table {%
-0.15 1.33
0.15 1.33
};
\addplot [black, forget plot]
table {%
0.15 1.41
0.45 1.41
};
\addplot [black, forget plot]
table {%
1.15 0.93
1.45 0.93
};
\addplot [black, forget plot]
table {%
1.45 1.83
1.75 1.83
};
\addplot [black, forget plot]
table {%
2.45 1.06
2.75 1.06
};
\addplot [black, forget plot]
table {%
2.75 nan
3.05 nan
};
\addplot [black, forget plot]
table {%
3.75 1.19
4.05 1.19
};
\addplot [black, forget plot]
table {%
4.05 nan
4.35 nan
};
\addplot [black, forget plot]
table {%
5.05 0.97
5.35 0.97
};
\addplot [black, forget plot]
table {%
5.35 nan
5.65 nan
};

% manual legend
\addlegendimage{area legend, draw=black, fill=blue}
\addlegendentry{P\&B}

\addlegendimage{area legend, draw=black, fill=darkgreen0640}
\addlegendentry{MIP}

\end{axis}

\end{tikzpicture}
    \caption{Integrality gap}
    \label{fig:pp-int-gap}
  \end{subfigure}
  \caption{Computational results for \gls{acr:pab} ($n=15$)}
  \label{fig:benchmark-boxplots}
\end{figure}

Here, we exclude instances that the commercial solver cannot solve from the reported result. Contrary, our \gls{acr:pab} approach solves all sets of instances to a median integrality gap of less than $1.33\%$ within the given time limit and the median solve time to find the first integer feasible solution was $1{,}800.69$ seconds even for the large instances. Figure~\ref{fig:pp-solve-time} and Figure~\ref{fig:pp-int-gap} show that both results are reasonably stable versus outliers. Specifically, the \gls{acr:pab} approach yields solutions below $2\%$ integrality gap for \revone{almost} all instances and it takes a maximum of less than $2{,}000$ seconds to find feasible solutions.

\vspace{2pt}
\begin{result}
Our algorithmic framework with a \gls{acr:pab} approach solves larger instances than our algorithmic framework with a commercial solver. The difference in solvable instance sizes reaches a factor of $6$, i.e., increases from $500$ to $3{,}000$ freight requests. 
\end{result} 
\vspace{2pt}

In Table~\ref{tab:bap_integrality_gap}, we compare our \gls{acr:pab} with the \gls{acr:bap} algorithm within our algorithmic framework. We report median, minimum and maximum integrality gaps remaining after $90$ minutes for both algorithms. 

\begin{table}[!b]
    \centering
    \caption{\gls{acr:bap} and P\&B integrality gap results ($n=15$)}
    \begin{tabular}{ccccccc}
        \toprule
        \multirow{2}{*}{\textbf{Instance size}} & \multicolumn{3}{c}{\textbf{\gls{acr:bap} integrality gap}} & \multicolumn{3}{c}{\textbf{P\&B integrality gap}} \\
        \cmidrule(lr){2-4} \cmidrule(lr){5-7}
        & \textbf{Median} & \textbf{Min.} & \textbf{Max.} & \textbf{Median} & \textbf{Min.} & \textbf{Max.} \\
        \midrule
        \midrule
        250 & 1.18\% & 0.84\% & 1.51\% & 1.56\% & 1.01\% & 1.95\% \\
        500 & 0.77\% & 0.54\% & 0.88\% & 0.93\% & 0.69\% & 1.36\% \\
        1,000 & 0.91\% & 0.59\% & 1.19\% & 1.06\% & 0.71\% & 2.17\% \\
        2,000 & 1.12\% & 0.57\% & 1.86\% & 1.19\% & 0.58\% & 2.32\% \\
        3,000 & 1.10\% & 0.37\% & 2.20\% & 0.97\% & 0.46\% & 1.49\% \\
        \bottomrule
    \end{tabular}
    \label{tab:bap_integrality_gap}
\end{table}

The \gls{acr:bap} approach decreases the median integrality gaps \revone{of instances with $250$ -- $2{,}000$ requests} to smaller values than the \gls{acr:pab} approach. \revone{However, the \gls{acr:pab} yields median integrality gaps that are at maximum $0.38$ percentage points higher than the ones from the \gls{acr:bap} algorithm} with much less \revone{pricing effort and improves on the median solution quality of the large instances with $3{,}000$ freight requests by $0.13$ percentage points. Figure~\ref{fig:comp-bab-and-pab} visualizes the differences between the two approaches. Figure~\ref{fig:abs-diff-pricing-effort} demonstrates that the \gls{acr:bap} approach relies on generating new columns. In our experiments with large instances that contain $3{,}000$ freight request, the \gls{acr:bap} approach generates an average of $\num{1.25e5}$ columns while the \gls{acr:pab} algorithm only generates $\num{0.6e5}$ columns. Furthermore, Figure~\ref{fig:rel-diff-bounds} shows the average difference in primal and dual bounds between the two algorithm. As can be seen, the \gls{acr:bap} algorithm improves on the dual bounds of the \gls{acr:pab} algorithm consistently. However, the improvement of less than $0.1$\% is marginal and the average differences in the primal bounds mostly determine the overall differences between the two algorithms with respect to the integrality gaps. In this context, the \gls{acr:bap} finds better primal bounds than the \gls{acr:pab} algorithm for most instances. However, in large instances with $3{,}000$ freight requests the effect reverses, and the \gls{acr:pab} yields better primal bounds than the \gls{acr:bap} algorithm. We observe this trend because the dual bounds become stronger when more requests are considered —-- if the aggregate demand remains constant  \citep[cf. p.50 in][]{ChekuriEtAl2009} --— and branching becomes more complex. Both effects disproportionately benefit the \gls{acr:pab} approach.}

\begin{figure}[!t]
  \begin{subfigure}[t]{0.48\textwidth}
    \centering
    % This file was created with tikzplotlib v0.10.1.
\begin{tikzpicture}

\definecolor{darkgray176}{RGB}{176,176,176}
\definecolor{lightgray204}{RGB}{204,204,204}
\definecolor{midnightblue848107}{RGB}{8,48,107}
\definecolor{skyblue147196222}{RGB}{147,196,222}

\begin{axis}[
legend cell align={left},
legend style={
  fill opacity=0.8,
  draw opacity=1,
  text opacity=1,
  at={(0.03,0.97)},
  anchor=north west,
  draw=lightgray204
},
tick align=outside,
tick pos=left,
x grid style={darkgray176},
xlabel={Instance size [\# Requests]},
xmin=112.5, xmax=3137.5,
xtick style={color=black},
xtick={250,1000,2000,3000},
xticklabel={\pgfmathprintnumber{\tick}}, %use comma,1000 sep={,}
scaled x ticks=false,
y grid style={darkgray176},
ylabel={\# Pricing problems},
ymajorgrids,
ymin=-1385.9639521619, ymax=159330.945839321,
ytick style={color=black},
width=0.7\linewidth,
height=0.62\linewidth,
scale only axis,
trim axis left, trim axis right,
]
\path [draw=skyblue147196222]
(axis cs:250,16409.3734074215)
--(axis cs:250,27655.6932592452);

\path [draw=skyblue147196222]
(axis cs:500,29945.6166855587)
--(axis cs:500,44025.1833144413);

\path [draw=skyblue147196222]
(axis cs:1000,51420.6185753395)
--(axis cs:1000,77971.2480913272);

\path [draw=skyblue147196222]
(axis cs:2000,85253.8070688222)
--(axis cs:2000,124505.259597845);

\path [draw=skyblue147196222]
(axis cs:3000,103208.234908777)
--(axis cs:3000,152025.63175789);

\addplot [semithick, skyblue147196222, mark=-, mark size=3, mark options={solid}, only marks, forget plot]
table {%
250 16409.3734074215
500 29945.6166855587
1000 51420.6185753395
2000 85253.8070688222
3000 103208.234908777
};
\addplot [semithick, skyblue147196222, mark=-, mark size=3, mark options={solid}, only marks, forget plot]
table {%
250 27655.6932592452
500 44025.1833144413
1000 77971.2480913272
2000 124505.259597845
3000 152025.63175789
};
\path [draw=midnightblue848107]
(axis cs:250,5919.35012926914)
--(axis cs:250,8676.91653739753);

\path [draw=midnightblue848107]
(axis cs:500,8815.81629779272)
--(axis cs:500,15819.6503688739);

\path [draw=midnightblue848107]
(axis cs:1000,14536.6870666561)
--(axis cs:1000,28256.7796000106);

\path [draw=midnightblue848107]
(axis cs:2000,30755.451035265)
--(axis cs:2000,47658.6822980683);

\path [draw=midnightblue848107]
(axis cs:3000,45906.8538437265)
--(axis cs:3000,79661.9461562735);

\addplot [semithick, midnightblue848107, mark=-, mark size=3, mark options={solid}, only marks, forget plot]
table {%
250 5919.35012926914
500 8815.81629779272
1000 14536.6870666561
2000 30755.451035265
3000 45906.8538437265
};
\addplot [semithick, midnightblue848107, mark=-, mark size=3, mark options={solid}, only marks, forget plot]
table {%
250 8676.91653739753
500 15819.6503688739
1000 28256.7796000106
2000 47658.6822980683
3000 79661.9461562735
};
\addplot [semithick, skyblue147196222, mark=*, mark size=3, mark options={solid}]
table {%
250 22032.5333333333
500 36985.4
1000 64695.9333333333
2000 104879.533333333
3000 127616.933333333
};
\addlegendentry{B\&P}
\addplot [semithick, midnightblue848107, mark=*, mark size=3, mark options={solid}]
table {%
250 7298.13333333333
500 12317.7333333333
1000 21396.7333333333
2000 39207.0666666667
3000 62784.4
};
\addlegendentry{P\&B}
\end{axis}

\end{tikzpicture}
    \centering \caption{Difference in pricing effort}
    \label{fig:abs-diff-pricing-effort}
  \end{subfigure}
  \begin{subfigure}[t]{0.48\textwidth}
    \centering
    % This file was created with tikzplotlib v0.10.1.
\begin{tikzpicture}

\definecolor{darkgray176}{RGB}{176,176,176}
\definecolor{darkorange25512714}{RGB}{255,127,14}
\definecolor{lightgray204}{RGB}{204,204,204}
\definecolor{steelblue31119180}{RGB}{31,119,180}

\begin{axis}[
legend cell align={left},
legend style={
  fill opacity=0.8,
  draw opacity=1,
  text opacity=1,
  at={(0.03,0.03)},
  anchor=south west,
  draw=lightgray204
},
tick align=outside,
tick pos=left,
x grid style={darkgray176},
xlabel={Instance size [\# Requests]},
xmin=112.5, xmax=3137.5,
xtick style={color=black},
xtick={250,1000,2000,3000},
xticklabel={\pgfmathprintnumber{\tick}}, %use comma,1000 sep={,}
scaled x ticks=false,
y grid style={darkgray176},
ylabel={$\Delta$Bound (P\&B - B\&P) [\%]},
ymajorgrids,
ymin=-0.00815940856040926, ymax=0.00638527555632335,
ytick style={color=black},
ytick={-0.01,-0.008,-0.006,-0.004,-0.002,0,0.002,0.004,0.006,0.008},
yticklabel={\pgfmathparse{\tick*100}\pgfmathprintnumber[fixed,precision=1]{\pgfmathresult}},
scaled y ticks=false,
width=0.7\linewidth,
height=0.62\linewidth,
scale only axis,
trim axis left, trim axis right,
]
\path [draw=steelblue31119180]
(axis cs:250,0.000488660651259744)
--(axis cs:250,0.00551548164393586);

\path [draw=steelblue31119180]
(axis cs:500,-0.000303446405209696)
--(axis cs:500,0.00453242196438688);

\path [draw=steelblue31119180]
(axis cs:1000,-0.00197365213201494)
--(axis cs:1000,0.00540811542578077);

\path [draw=steelblue31119180]
(axis cs:2000,-0.00413614380170352)
--(axis cs:2000,0.00572415355101732);

\path [draw=steelblue31119180]
(axis cs:3000,-0.00749828655510324)
--(axis cs:3000,0.00278200059786555);

\addplot [semithick, steelblue31119180, mark=-, mark size=3, mark options={solid}, only marks, forget plot]
table {%
250 0.000488660651259744
500 -0.000303446405209696
1000 -0.00197365213201494
2000 -0.00413614380170352
3000 -0.00749828655510324
};
% \addlegendentry{Dual bound}
\addplot [semithick, steelblue31119180, mark=-, mark size=3, mark options={solid}, only marks, forget plot]
table {%
250 0.00551548164393586
500 0.00453242196438688
1000 0.00540811542578077
2000 0.00572415355101732
3000 0.00278200059786555
};
% \addlegendentry{Dual bound}
\path [draw=darkorange25512714]
(axis cs:250,-0.000706145852993413)
--(axis cs:250,0.000163928369079767);

\path [draw=darkorange25512714]
(axis cs:500,-0.000431877615299877)
--(axis cs:500,0.000119745974252951);

\path [draw=darkorange25512714]
(axis cs:1000,-0.000748003105311682)
--(axis cs:1000,0.000104223553785038);

\path [draw=darkorange25512714]
(axis cs:2000,-0.000521800632923909)
--(axis cs:2000,0.000134230836311535);

\path [draw=darkorange25512714]
(axis cs:3000,-0.000505244843209761)
--(axis cs:3000,9.34908316857007e-05);

\addplot [semithick, darkorange25512714, mark=-, mark size=3, mark options={solid}, only marks, forget plot]
table {%
250 -0.000706145852993413
500 -0.000431877615299877
1000 -0.000748003105311682
2000 -0.000521800632923909
3000 -0.000505244843209761
};
% \addlegendentry{Dual bound}
\addplot [semithick, darkorange25512714, mark=-, mark size=3, mark options={solid}, only marks, forget plot]
table {%
250 0.000163928369079767
500 0.000119745974252951
1000 0.000104223553785038
2000 0.000134230836311535
3000 9.34908316857007e-05
};
% \addlegendentry{Dual bound}
\addplot [black, forget plot]
table {%
112.5 0
3137.5 0
};
\addplot [semithick, steelblue31119180, mark=*, mark size=3, mark options={solid}]
table {%
250 0.0030020711475978
500 0.00211448777958859
1000 0.00171723164688291
2000 0.0007940048746569
3000 -0.00235814297861884
};
\addlegendentry{Primal bound}
\addplot [semithick, darkorange25512714, mark=*, mark size=3, mark options={solid}]
table {%
250 -0.000271108741956823
500 -0.000156065820523463
1000 -0.000321889775763322
2000 -0.000193784898306187
3000 -0.00020587700576203
};
\addlegendentry{Dual bound}

\end{axis}

\end{tikzpicture}
    \vspace{-0.5cm}
    \centering \caption{Relative difference in bounds}
    \label{fig:rel-diff-bounds}
  \end{subfigure}
  \caption{ Difference between \gls{acr:pab} and \gls{acr:bap}. Means are denoted by circles, and error bars indicate standard deviations ($n=15$)}
  \label{fig:comp-bab-and-pab}
\end{figure}

The difference \revone{between the two algorithms} of up to $0.38$ percentage points in the median integrality gaps is marginal due to the low integrality gaps in general. Note that the integrality gaps reported for the \gls{acr:pab} in Table~\ref{tab:bap_integrality_gap} differ from the ones previously reported in Table~\ref{tab:performance}. This effect occurs because in Table~\ref{tab:performance}, we compute the gap based on the tightest lower bound found with either the commercial solver or the \gls{acr:pab} approach. Because of the good performance of \gls{acr:pab} we presented our algorithmic framework focusing on the \gls{acr:pab} approach. 

\vspace{2pt}
\begin{result}
    Replacing the \gls{acr:pab} algorithm with a full \gls{acr:bap} algorithm in our algorithmic framework \revone{can improve median integrality gaps by up to $0.38$} percentage points. 
\end{result}
\vspace{2pt}

\revone{
To evaluate our algorithms performance across different problem instances, we scale our instances from covering a $5$ hours time period to $6$, $7$, $8$, and $9$ hours respectively and show the results in Figure~\ref{fig:scaled-network}. The scaled networks comprise up to $91$ vehicles and substantially enlarge the time-expanded graph.} 

\begin{figure}[!t]
  \begin{subfigure}[t]{0.49\textwidth}
    \centering
    % This file was created with tikzplotlib v0.10.1.
\begin{tikzpicture}

\definecolor{darkgray176}{RGB}{176,176,176}
\definecolor{gray}{RGB}{128,128,128}
\definecolor{lightgray204}{RGB}{204,204,204}
\definecolor{midnightblue848107}{RGB}{8,48,107}
\definecolor{powderblue198219239}{RGB}{198,219,239}

\begin{axis}[
grid style={dotted,gray},
height=0.6\linewidth,
legend columns=2,
legend cell align={left},
legend style={
  fill opacity=0.8,
  draw opacity=1,
  text opacity=1,
  at={(0.03,0.99)},
  anchor=north west,
  draw=lightgray204
},
scale only axis,
tick align=outside,
tick pos=left,
trim axis left, trim axis right,
width=0.6\linewidth,
x grid style={darkgray176},
xlabel={Duration [h]},
xmin=-0.9, xmax=4.9,
xtick style={color=black},
xtick={0,1,2,3,4},
xticklabels={5,6,7,8,9},
y grid style={gray},
ylabel={\# Pricing Problems},
ymajorgrids,
ymajorgrids,
ymin=40000, ymax=280000,
ytick style={color=black},
ytick={{ 40000, 80000, 120000, 160000, 200000, 240000, 280000 }}
]
\path [draw=black, fill=powderblue198219239]
(axis cs:-0.45,114463)
--(axis cs:-0.09,114463)
--(axis cs:-0.09,142234.5)
--(axis cs:-0.45,142234.5)
--(axis cs:-0.45,114463)
--cycle;
\addplot [black, forget plot]
table {%
-0.27 114463
-0.27 91321
};
\addplot [black, forget plot]
table {%
-0.27 142234.5
-0.27 172152
};
\addplot [black, forget plot]
table {%
-0.36 91321
-0.18 91321
};
\addplot [black, forget plot]
table {%
-0.36 172152
-0.18 172152
};
\path [draw=black, fill=midnightblue848107]
(axis cs:0.09,50402.5)
--(axis cs:0.45,50402.5)
--(axis cs:0.45,73390.5)
--(axis cs:0.09,73390.5)
--(axis cs:0.09,50402.5)
--cycle;
\addplot [black, forget plot]
table {%
0.27 50402.5
0.27 44969
};
\addplot [black, forget plot]
table {%
0.27 73390.5
0.27 102456
};
\addplot [black, forget plot]
table {%
0.18 44969
0.36 44969
};
\addplot [black, forget plot]
table {%
0.18 102456
0.36 102456
};
\path [draw=black, fill=powderblue198219239]
(axis cs:0.55,118372)
--(axis cs:0.91,118372)
--(axis cs:0.91,127994.5)
--(axis cs:0.55,127994.5)
--(axis cs:0.55,118372)
--cycle;
\addplot [black, forget plot]
table {%
0.73 118372
0.73 113164
};
\addplot [black, forget plot]
table {%
0.73 127994.5
0.73 128545
};
\addplot [black, forget plot]
table {%
0.64 113164
0.82 113164
};
\addplot [black, forget plot]
table {%
0.64 128545
0.82 128545
};
\addplot [black, mark=o, mark size=3, mark options={solid,fill opacity=0}, only marks, forget plot]
table {%
0.73 101579
0.73 143165
0.73 150527
0.73 145167
};
\path [draw=black, fill=midnightblue848107]
(axis cs:1.09,51171.5)
--(axis cs:1.45,51171.5)
--(axis cs:1.45,67398)
--(axis cs:1.09,67398)
--(axis cs:1.09,51171.5)
--cycle;
\addplot [black, forget plot]
table {%
1.27 51171.5
1.27 45066
};
\addplot [black, forget plot]
table {%
1.27 67398
1.27 79456
};
\addplot [black, forget plot]
table {%
1.18 45066
1.36 45066
};
\addplot [black, forget plot]
table {%
1.18 79456
1.36 79456
};
\addplot [black, mark=o, mark size=3, mark options={solid,fill opacity=0}, only marks, forget plot]
table {%
1.27 97054
};
\path [draw=black, fill=powderblue198219239]
(axis cs:1.55,109733.5)
--(axis cs:1.91,109733.5)
--(axis cs:1.91,152568)
--(axis cs:1.55,152568)
--(axis cs:1.55,109733.5)
--cycle;
\addplot [black, forget plot]
table {%
1.73 109733.5
1.73 78625
};
\addplot [black, forget plot]
table {%
1.73 152568
1.73 198994
};
\addplot [black, forget plot]
table {%
1.64 78625
1.82 78625
};
\addplot [black, forget plot]
table {%
1.64 198994
1.82 198994
};
\path [draw=black, fill=midnightblue848107]
(axis cs:2.09,50096)
--(axis cs:2.45,50096)
--(axis cs:2.45,81781.5)
--(axis cs:2.09,81781.5)
--(axis cs:2.09,50096)
--cycle;
\addplot [black, forget plot]
table {%
2.27 50096
2.27 47300
};
\addplot [black, forget plot]
table {%
2.27 81781.5
2.27 129153
};
\addplot [black, forget plot]
table {%
2.18 47300
2.36 47300
};
\addplot [black, forget plot]
table {%
2.18 129153
2.36 129153
};
\path [draw=black, fill=powderblue198219239]
(axis cs:2.55,110922.25)
--(axis cs:2.91,110922.25)
--(axis cs:2.91,182060.5)
--(axis cs:2.55,182060.5)
--(axis cs:2.55,110922.25)
--cycle;
\addplot [black, forget plot]
table {%
2.73 110922.25
2.73 81719
};
\addplot [black, forget plot]
table {%
2.73 182060.5
2.73 196166
};
\addplot [black, forget plot]
table {%
2.64 81719
2.82 81719
};
\addplot [black, forget plot]
table {%
2.64 196166
2.82 196166
};
\path [draw=black, fill=midnightblue848107]
(axis cs:3.09,69839)
--(axis cs:3.45,69839)
--(axis cs:3.45,102720.5)
--(axis cs:3.09,102720.5)
--(axis cs:3.09,69839)
--cycle;
\addplot [black, forget plot]
table {%
3.27 69839
3.27 58393
};
\addplot [black, forget plot]
table {%
3.27 102720.5
3.27 140897
};
\addplot [black, forget plot]
table {%
3.18 58393
3.36 58393
};
\addplot [black, forget plot]
table {%
3.18 140897
3.36 140897
};
\path [draw=black, fill=powderblue198219239]
(axis cs:3.55,106631.25)
--(axis cs:3.91,106631.25)
--(axis cs:3.91,159919.5)
--(axis cs:3.55,159919.5)
--(axis cs:3.55,106631.25)
--cycle;
\addplot [black, forget plot]
table {%
3.73 106631.25
3.73 76913
};
\addplot [black, forget plot]
table {%
3.73 159919.5
3.73 187080
};
\addplot [black, forget plot]
table {%
3.64 76913
3.82 76913
};
\addplot [black, forget plot]
table {%
3.64 187080
3.82 187080
};
\addplot [black, mark=o, mark size=3, mark options={solid,fill opacity=0}, only marks, forget plot]
table {%
3.73 255244
};
\path [draw=black, fill=midnightblue848107]
(axis cs:4.09,66547)
--(axis cs:4.45,66547)
--(axis cs:4.45,94184.5)
--(axis cs:4.09,94184.5)
--(axis cs:4.09,66547)
--cycle;
\addplot [black, forget plot]
table {%
4.27 66547
4.27 63084
};
\addplot [black, forget plot]
table {%
4.27 94184.5
4.27 108486
};
\addplot [black, forget plot]
table {%
4.18 63084
4.36 63084
};
\addplot [black, forget plot]
table {%
4.18 108486
4.36 108486
};
\addplot [black, mark=o, mark size=3, mark options={solid,fill opacity=0}, only marks, forget plot]
table {%
4.27 136705
};
\addplot [black, forget plot]
table {%
-0.45 122080
-0.09 122080
};
\addplot [black, forget plot]
table {%
0.09 56493
0.45 56493
};
\addplot [black, forget plot]
table {%
0.55 121885
0.91 121885
};
\addplot [black, forget plot]
table {%
1.09 59891
1.45 59891
};
\addplot [black, forget plot]
table {%
1.55 144699
1.91 144699
};
\addplot [black, forget plot]
table {%
2.09 63778
2.45 63778
};
\addplot [black, forget plot]
table {%
2.55 123854.5
2.91 123854.5
};
\addplot [black, forget plot]
table {%
3.09 78096
3.45 78096
};
\addplot [black, forget plot]
table {%
3.55 137587
3.91 137587
};
\addplot [black, forget plot]
table {%
4.09 71698
4.45 71698
};
% --- Manual legend entries
\addlegendimage{area legend,draw=black,fill=powderblue198219239}
\addlegendentry{B\&P}

\addlegendimage{area legend,draw=black,fill=midnightblue848107}
\addlegendentry{P\&B}

% Marker for sample counts (hollow circle with black edge)
\addlegendimage{only marks, mark=diamond*, draw=black, fill=black}
\addlegendentry{Solved}
\end{axis}

\begin{axis}[
axis y line=right,
grid style={dotted,gray},
height=0.6\linewidth,
scale only axis,
scaled y ticks=false,
tick align=outside,
trim axis left, trim axis right,
width=0.6\linewidth,
x grid style={darkgray176},
xmin=-0.9, xmax=4.9,
xtick pos=left,
xtick style={color=black},
xtick={0,1,2,3,4},
xticklabels={5,6,7,8,9},
y grid style={darkgray176},
y tick scale label code/.code={},
ylabel={\# Instances solved},
ymajorgrids,
ymin=0, ymax=24,
ytick pos=right,
ytick style={color=black},
ytick={{ 0, 4, 8, 12, 16, 20, 24}},
yticklabel style={anchor=west}
]
\addplot [draw=black, fill=black, mark=diamond*, only marks]
table{%
x  y
-0.27 15
0.27 15
0.73 15
1.27 15
1.73 15
2.27 15
2.73 14
3.27 15
3.73 12
4.27 15
};
\end{axis}

\end{tikzpicture}
    \centering \caption{Pricing effort in scaled networks}
    \label{fig:scaled-network-pricing-effort}
  \end{subfigure}\hfill
  \begin{subfigure}[t]{0.49\textwidth}
    \centering
    % This file was created with tikzplotlib v0.10.1.
\begin{tikzpicture}

\definecolor{darkgray176}{RGB}{176,176,176}
\definecolor{gray}{RGB}{128,128,128}
\definecolor{lightgray204}{RGB}{204,204,204}
\definecolor{midnightblue848107}{RGB}{8,48,107}
\definecolor{powderblue198219239}{RGB}{198,219,239}

\begin{axis}[
grid style={dotted,gray},
height=0.6\linewidth,
legend columns=2,
legend cell align={left},
legend style={
  fill opacity=0.8,
  draw opacity=1,
  text opacity=1,
  at={(0.03,0.99)},
  anchor=north west,
  draw=lightgray204
},
scale only axis,
scaled y ticks=false,
tick align=outside,
tick pos=left,
trim axis left, trim axis right,
width=0.6\linewidth,
x grid style={darkgray176},
xlabel={Duration [h]},
xmin=-0.9, xmax=4.9,
xtick style={color=black},
xtick={0,1,2,3,4},
xticklabels={5,6,7,8,9},
y grid style={gray},
y tick scale label code/.code={},
ylabel={Integrality Gap [\%]},
ymajorgrids,
ymajorgrids,
ymin=0, ymax=0.12,
ytick style={color=black},
ytick={{0,0.02,0.04,0.06,0.08, 0.10, 0.12 }},
yticklabels={0.0,2.0,4.0,6.0,8.0,10.0,12.0}
]
\path [draw=black, fill=powderblue198219239]
(axis cs:-0.45,0.00722740761388559)
--(axis cs:-0.09,0.00722740761388559)
--(axis cs:-0.09,0.0136650332831919)
--(axis cs:-0.45,0.0136650332831919)
--(axis cs:-0.45,0.00722740761388559)
--cycle;
\addplot [black, forget plot]
table {%
-0.27 0.00722740761388559
-0.27 0.00512067901649259
};
\addplot [black, forget plot]
table {%
-0.27 0.0136650332831919
-0.27 0.0190106374400118
};
\addplot [black, forget plot]
table {%
-0.36 0.00512067901649259
-0.18 0.00512067901649259
};
\addplot [black, forget plot]
table {%
-0.36 0.0190106374400118
-0.18 0.0190106374400118
};
\path [draw=black, fill=midnightblue848107]
(axis cs:0.09,0.00814404206937826)
--(axis cs:0.45,0.00814404206937826)
--(axis cs:0.45,0.0117997845370943)
--(axis cs:0.09,0.0117997845370943)
--(axis cs:0.09,0.00814404206937826)
--cycle;
\addplot [black, forget plot]
table {%
0.27 0.00814404206937826
0.27 0.00456867057273019
};
\addplot [black, forget plot]
table {%
0.27 0.0117997845370943
0.27 0.0148670011670583
};
\addplot [black, forget plot]
table {%
0.18 0.00456867057273019
0.36 0.00456867057273019
};
\addplot [black, forget plot]
table {%
0.18 0.0148670011670583
0.36 0.0148670011670583
};
\path [draw=black, fill=powderblue198219239]
(axis cs:0.55,0.0123640193321297)
--(axis cs:0.91,0.0123640193321297)
--(axis cs:0.91,0.0184507054017804)
--(axis cs:0.55,0.0184507054017804)
--(axis cs:0.55,0.0123640193321297)
--cycle;
\addplot [black, forget plot]
table {%
0.73 0.0123640193321297
0.73 0.00894219887770715
};
\addplot [black, forget plot]
table {%
0.73 0.0184507054017804
0.73 0.0254437838905061
};
\addplot [black, forget plot]
table {%
0.64 0.00894219887770715
0.82 0.00894219887770715
};
\addplot [black, forget plot]
table {%
0.64 0.0254437838905061
0.82 0.0254437838905061
};
\path [draw=black, fill=midnightblue848107]
(axis cs:1.09,0.0101681934275335)
--(axis cs:1.45,0.0101681934275335)
--(axis cs:1.45,0.027466656550503)
--(axis cs:1.09,0.027466656550503)
--(axis cs:1.09,0.0101681934275335)
--cycle;
\addplot [black, forget plot]
table {%
1.27 0.0101681934275335
1.27 0.00638888996707701
};
\addplot [black, forget plot]
table {%
1.27 0.027466656550503
1.27 0.0442347683049945
};
\addplot [black, forget plot]
table {%
1.18 0.00638888996707701
1.36 0.00638888996707701
};
\addplot [black, forget plot]
table {%
1.18 0.0442347683049945
1.36 0.0442347683049945
};
\path [draw=black, fill=powderblue198219239]
(axis cs:1.55,0.0168969406475213)
--(axis cs:1.91,0.0168969406475213)
--(axis cs:1.91,0.0324455992307978)
--(axis cs:1.55,0.0324455992307978)
--(axis cs:1.55,0.0168969406475213)
--cycle;
\addplot [black, forget plot]
table {%
1.73 0.0168969406475213
1.73 0.00895415951510789
};
\addplot [black, forget plot]
table {%
1.73 0.0324455992307978
1.73 0.052651355380123
};
\addplot [black, forget plot]
table {%
1.64 0.00895415951510789
1.82 0.00895415951510789
};
\addplot [black, forget plot]
table {%
1.64 0.052651355380123
1.82 0.052651355380123
};
\path [draw=black, fill=midnightblue848107]
(axis cs:2.09,0.0143388141562843)
--(axis cs:2.45,0.0143388141562843)
--(axis cs:2.45,0.0381396657234027)
--(axis cs:2.09,0.0381396657234027)
--(axis cs:2.09,0.0143388141562843)
--cycle;
\addplot [black, forget plot]
table {%
2.27 0.0143388141562843
2.27 0.00910671328895312
};
\addplot [black, forget plot]
table {%
2.27 0.0381396657234027
2.27 0.0450051838831484
};
\addplot [black, forget plot]
table {%
2.18 0.00910671328895312
2.36 0.00910671328895312
};
\addplot [black, forget plot]
table {%
2.18 0.0450051838831484
2.36 0.0450051838831484
};
\path [draw=black, fill=powderblue198219239]
(axis cs:2.55,0.0265978959153693)
--(axis cs:2.91,0.0265978959153693)
--(axis cs:2.91,0.0537693986119598)
--(axis cs:2.55,0.0537693986119598)
--(axis cs:2.55,0.0265978959153693)
--cycle;
\addplot [black, forget plot]
table {%
2.73 0.0265978959153693
2.73 0.0147522976336601
};
\addplot [black, forget plot]
table {%
2.73 0.0537693986119598
2.73 0.0658187655798503
};
\addplot [black, forget plot]
table {%
2.64 0.0147522976336601
2.82 0.0147522976336601
};
\addplot [black, forget plot]
table {%
2.64 0.0658187655798503
2.82 0.0658187655798503
};
\path [draw=black, fill=midnightblue848107]
(axis cs:3.09,0.0262175684966805)
--(axis cs:3.45,0.0262175684966805)
--(axis cs:3.45,0.0430562325840288)
--(axis cs:3.09,0.0430562325840288)
--(axis cs:3.09,0.0262175684966805)
--cycle;
\addplot [black, forget plot]
table {%
3.27 0.0262175684966805
3.27 0.00915296481295121
};
\addplot [black, forget plot]
table {%
3.27 0.0430562325840288
3.27 0.0601303812117522
};
\addplot [black, forget plot]
table {%
3.18 0.00915296481295121
3.36 0.00915296481295121
};
\addplot [black, forget plot]
table {%
3.18 0.0601303812117522
3.36 0.0601303812117522
};
\path [draw=black, fill=powderblue198219239]
(axis cs:3.55,0.022384243688763)
--(axis cs:3.91,0.022384243688763)
--(axis cs:3.91,0.047700797974724)
--(axis cs:3.55,0.047700797974724)
--(axis cs:3.55,0.022384243688763)
--cycle;
\addplot [black, forget plot]
table {%
3.73 0.022384243688763
3.73 0.016509624599866
};
\addplot [black, forget plot]
table {%
3.73 0.047700797974724
3.73 0.0701903835232936
};
\addplot [black, forget plot]
table {%
3.64 0.016509624599866
3.82 0.016509624599866
};
\addplot [black, forget plot]
table {%
3.64 0.0701903835232936
3.82 0.0701903835232936
};
\path [draw=black, fill=midnightblue848107]
(axis cs:4.09,0.0307639469087408)
--(axis cs:4.45,0.0307639469087408)
--(axis cs:4.45,0.0517978977969554)
--(axis cs:4.09,0.0517978977969554)
--(axis cs:4.09,0.0307639469087408)
--cycle;
\addplot [black, forget plot]
table {%
4.27 0.0307639469087408
4.27 0.0203096917413074
};
\addplot [black, forget plot]
table {%
4.27 0.0517978977969554
4.27 0.0601489509506329
};
\addplot [black, forget plot]
table {%
4.18 0.0203096917413074
4.36 0.0203096917413074
};
\addplot [black, forget plot]
table {%
4.18 0.0601489509506329
4.36 0.0601489509506329
};
\addplot [black, forget plot]
table {%
-0.45 0.0118391411709101
-0.09 0.0118391411709101
};
\addplot [black, forget plot]
table {%
0.09 0.00965274690987526
0.45 0.00965274690987526
};
\addplot [black, forget plot]
table {%
0.55 0.0149647328867268
0.91 0.0149647328867268
};
\addplot [black, forget plot]
table {%
1.09 0.0166338412495027
1.45 0.0166338412495027
};
\addplot [black, forget plot]
table {%
1.55 0.0264892906966388
1.91 0.0264892906966388
};
\addplot [black, forget plot]
table {%
2.09 0.0273339184127913
2.45 0.0273339184127913
};
\addplot [black, forget plot]
table {%
2.55 0.0366084560745872
2.91 0.0366084560745872
};
\addplot [black, forget plot]
table {%
3.09 0.0379610682390909
3.45 0.0379610682390909
};
\addplot [black, forget plot]
table {%
3.55 0.0363109430237
3.91 0.0363109430237
};
\addplot [black, forget plot]
table {%
4.09 0.0425609119064963
4.45 0.0425609119064963
};

% --- Manual legend entries
\addlegendimage{area legend,draw=black,fill=powderblue198219239}
\addlegendentry{B\&P}

\addlegendimage{area legend,draw=black,fill=midnightblue848107}
\addlegendentry{P\&B}

% Marker for sample counts (hollow circle with black edge)
\addlegendimage{only marks, mark=diamond*, draw=black, fill=black}
\addlegendentry{Solved}
\end{axis}

\begin{axis}[
axis y line=right,
grid style={dotted,gray},
height=0.6\linewidth,
scale only axis,
scaled y ticks=false,
tick align=outside,
trim axis left, trim axis right,
width=0.6\linewidth,
x grid style={darkgray176},
xmin=-0.9, xmax=4.9,
xtick pos=left,
xtick style={color=black},
xtick={0,1,2,3,4},
xticklabels={5,6,7,8,9},
y grid style={darkgray176},
y tick scale label code/.code={},
ylabel={\# Instances solved},
ymajorgrids,
ymin=0, ymax=24,
ytick pos=right,
ytick style={color=black},
ytick={{ 0, 4, 8, 12, 16, 20, 24}},
yticklabel style={anchor=west}
]
\addplot [draw=black, fill=black, mark=diamond*, only marks]
table{%
x  y
-0.27 15
0.27 15
0.73 15
1.27 15
1.73 15
2.27 15
2.73 14
3.27 15
3.73 12
4.27 15
};
\end{axis}

\end{tikzpicture}
    \centering \caption{Integrality gaps in scaled networks}
    \label{fig:scaled-network-integrality-gap}
  \end{subfigure}
  \caption{Results in larger networks ($n=15$, $|\SetOfRequests|=3{,}000$)}
  \label{fig:scaled-network}
\end{figure}

\revone{On the demand side, we re-sample $10{,}000$ passenger requests and adjust the total passenger travel demand to reflect the longer horizons. For freight requests, we re-use the previously sampled instances with $3{,}000$ requests and only extend their latest time of service-completion according to the adjusted end of the evaluated time period. 
Figure~\ref{fig:scaled-network-pricing-effort} reports the number of generated columns in either the \gls{acr:pab} or the \gls{acr:bap} algorithm when being applied to the scaled networks. As can be seen, the median number of generated columns increases at increasing network sizes for the \gls{acr:pab} algorithm. However, this effect reverses when scaling the evaluated time period to $9$ hours, i.e., from $6$ a.m. to $3$ p.m. This inverse effect occurs due to the extended graph size that circumvents the convergence of the \gls{acr:cg} step to solve the linear relaxation of the master problem to an $\epsilon$-precise solution. Instead, the \gls{acr:cg} terminates prematurely and the algorithm proceeds with the branching. Thus, the total time devoted to the \gls{acr:cg} remains approximately constant for $8$ and $9$ hours of operations, yet the algorithm solves fewer pricing problems because each individual pricing subproblem is more complex in the larger network. The \gls{acr:bap} approach does not exhibit this trend because the algorithm generates new columns repeatedly by design. 
Figure~\ref{fig:scaled-network-integrality-gap} shows the associated integrality gaps. In this context, the \gls{acr:pab} algorithm cannot close the median gaps below $2\%$ for networks covering more than the $5$ hours that we consider by default in our analysis. Nevertheless, the solution quality degrades slowly and the gaps remain below $6$\%. The performance of the \gls{acr:bap} algorithm in terms of median integrality gaps is similar. However, the \gls{acr:bap} algorithm does not yield integer feasible solutions for $3$ out of $15$ instances in which $9$ hours of operations are reflected in the time-expanded graph.}

\vspace{2pt}
\revone{
\begin{result}
    The integrality gaps from the \gls{acr:pab} algorithm increase in larger networks. However, the increase is slow and in large networks representing $9$ hours of \gls{acr:pt} operations, the integrality gaps remain below $6$\%.
\end{result}
}

\subsection{Managerial insights} \label{subsection:managerial-insights}
\noindent The acceptance of a request by the municipality depends on the relation between transportation costs of the request in the \gls{acr:pts} and penalty costs for rejection. The transportation cost in the \gls{acr:pts}, on the one hand, depends on the externality cost associated with loading and unloading operations. On the other hand, the penalty cost of rejecting a request depends on the externality costs of truck based delivery. In our base case scenario, we set the arc cost on the transit arcs as $\Cost_{\UnbracedArc}=0.3, \medspace \Arc \in \SetOfArcs^{\Transit}$ and the externality cost for truck based delivery to be \euro{0.3} per vehicle and kilometer. 

In our base case with $3{,}000$ freight requests, we obtain a freight request rejection rate of $42.7\%$ (see Figure~\ref{fig:rejected-share-base}). In the following, we discuss the characteristics of this base case in comparison to two extreme scenarios: an optimistic scenario with truck delivery costs of \euro{0.4} and $\Cost_{\UnbracedArc}=0.3, \medspace \Arc \in \SetOfArcs^{\Transit}$ that favors the acceptance of freight deliveries (see Figure~\ref{fig:rejected-share-opt}), as well as a pessimistic scenario with truck delivery costs of \euro{0.3} and $\Cost_{\UnbracedArc}=0.4, \medspace \Arc \in \SetOfArcs^{\Transit}$ that reduces the acceptance of freight deliveries (see Figure~\ref{fig:rejected-share-pess}). 

\begin{figure}[!tb]
\centering
\scriptsize

\definecolor{darkblue}{rgb}{0.00, 0.44, 0.75}  
\definecolor{lightblue}{rgb}{0.82, 0.82, 0.82}

\centering
\begin{subfigure}[t]{0.34\textwidth} 
    \centering
    \begin{tikzpicture}[baseline]
        \pie[
            rotate=90,
            sum=auto,
            radius=1.5,
            color={darkblue, lightblue},
            text=pin,    
            hide number 
        ]{
            57.3/57.3\%,
            42.7/42.7\%
        }
    \end{tikzpicture}
    \vspace{0.75em}
    \caption{Base case}
    \label{fig:rejected-share-base}
\end{subfigure}
\begin{subfigure}[t]{0.30\textwidth}
    \centering
    \begin{tikzpicture}[baseline]
          \pie[
            rotate=45,
            sum=auto,
            radius=1.5,
            color={darkblue, lightblue},
            text=pin,     
            hide number   
          ]{
            99.6/99.6\%,  
            0.4/0.4\%
          }
    \end{tikzpicture}
    \caption{Optimistic scenario}
    \label{fig:rejected-share-opt}
\end{subfigure}
\begin{subfigure}[t]{0.30\textwidth}
    \centering
    \begin{tikzpicture}[baseline]
        \pie[
            rotate=90,
            sum=auto,
            radius=1.5,
            color={darkblue, lightblue},
            text=pin,    
            hide number  
        ]{
            32.6/32.6\%,
            67.4/67.4\%
        }
    \end{tikzpicture}
    \caption{Pessimistic scenario}
    \label{fig:rejected-share-pess}
\end{subfigure}
\vspace{0.1cm}
\begin{tikzpicture}
    \draw[darkblue, fill=darkblue] (0,0) rectangle (0.3,0.3);
    \node[right] at (0.5,0.15) {\footnotesize Accepted freight requests};

    \draw[lightblue, fill=lightblue] (6,0) rectangle (6.3,0.3); 
    \node[right] at (6.5,0.15) {\footnotesize Rejected freight requests};
\end{tikzpicture}

\caption{Cargo-hitching penetration across different scenarios ($n=15$)}
\label{fig:rejected-share}
\end{figure}

\begin{figure}[!t]
    \begin{subfigure}{0.99\textwidth}
        \centering
        \includegraphics[width=0.3\linewidth]{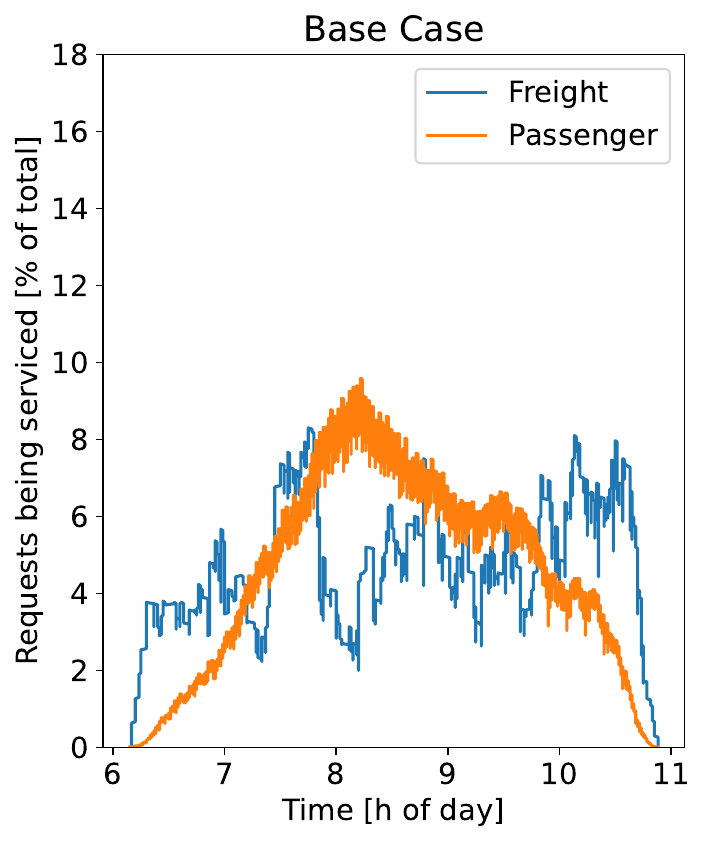}
        \includegraphics[width=0.3\linewidth]{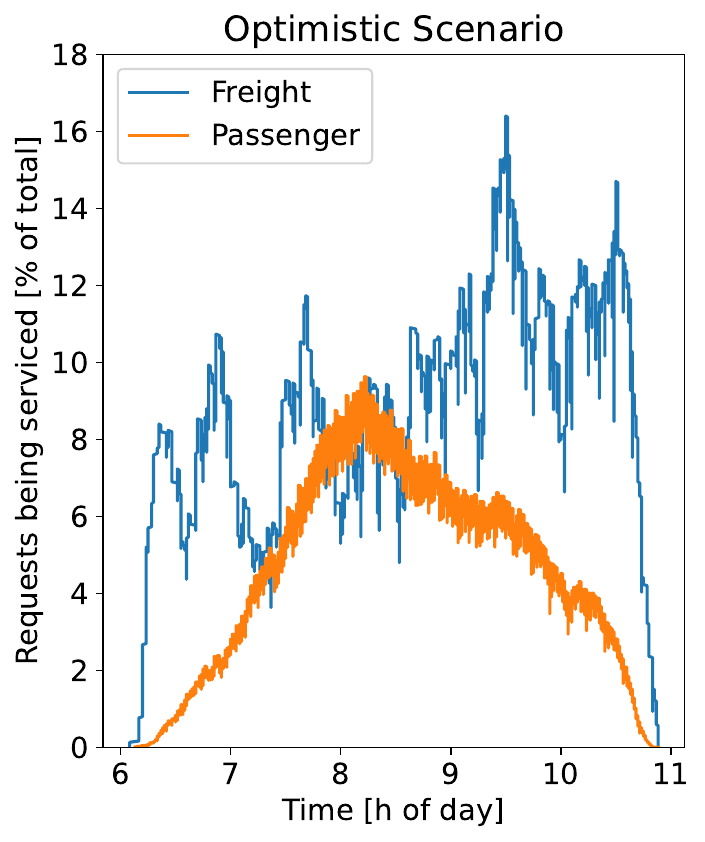}
        \includegraphics[width=0.3\linewidth]{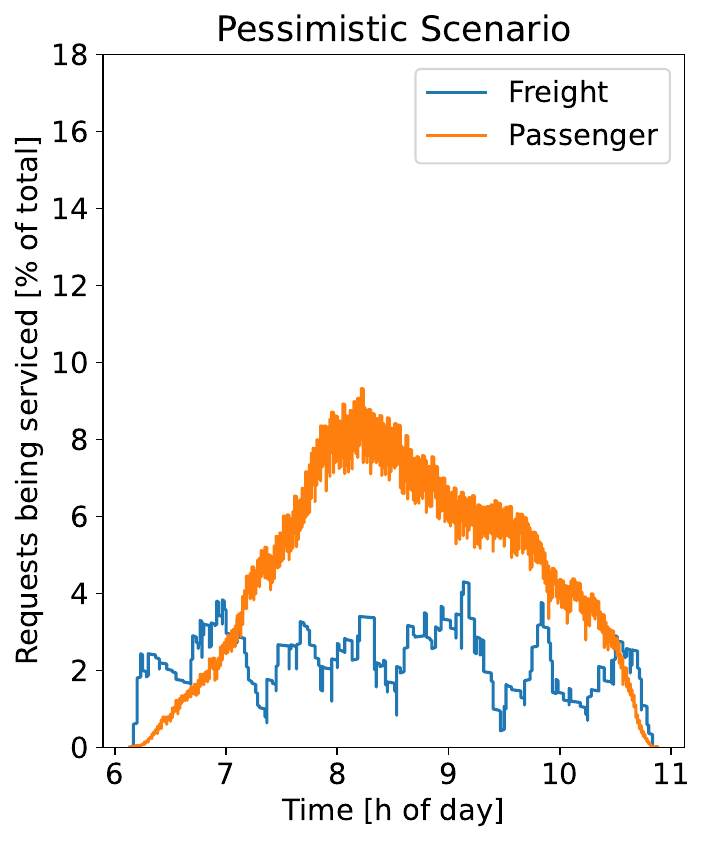}
        \captionsetup{width=0.9\columnwidth}
        \caption{\centering Temporal utilization}
        \label{fig:temp-utilization}
    \end{subfigure}
    \begin{subfigure}{0.99\textwidth}
        \centering
        \includegraphics[width=0.3\linewidth]{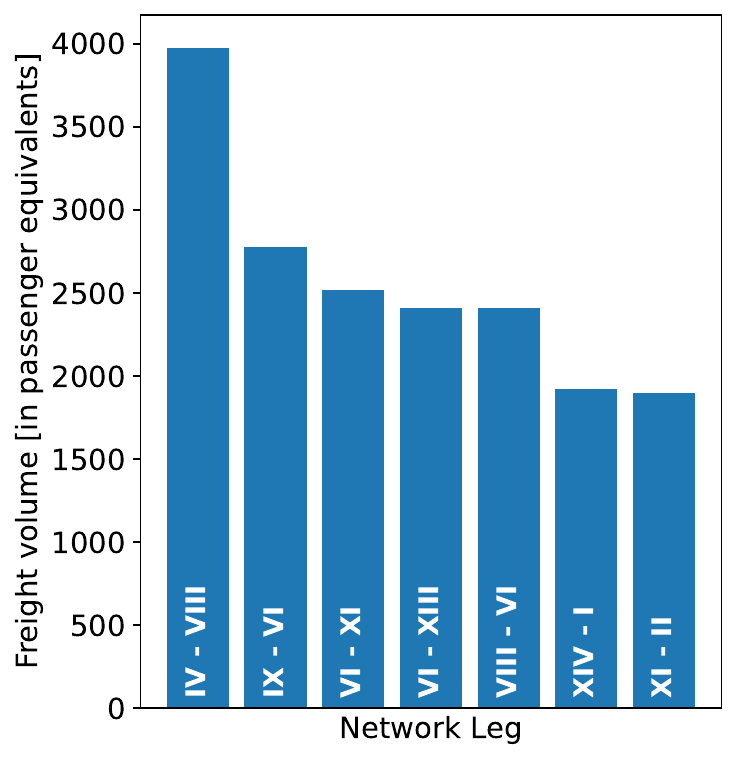}
        \includegraphics[width=0.3\linewidth]{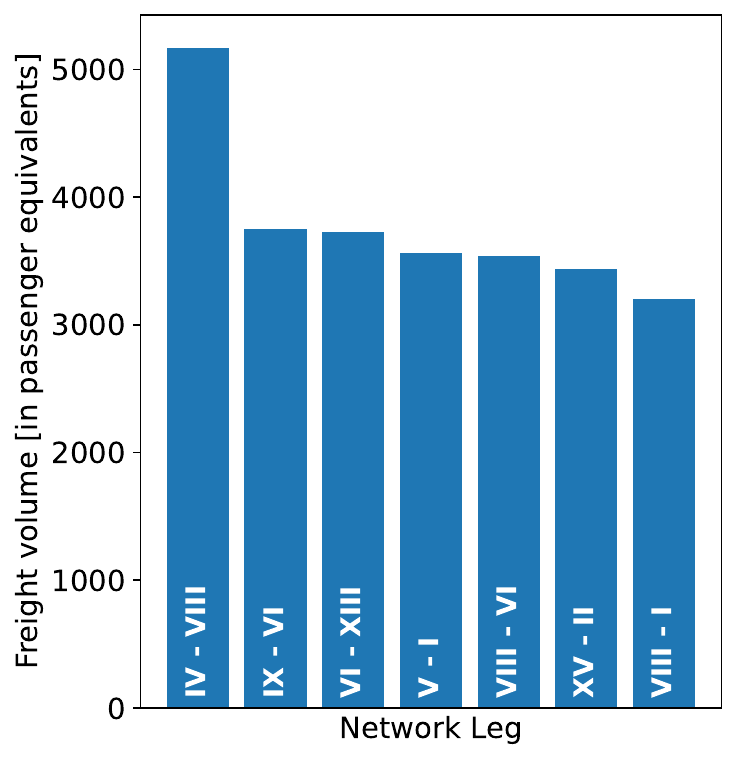}
        \includegraphics[width=0.3\linewidth]{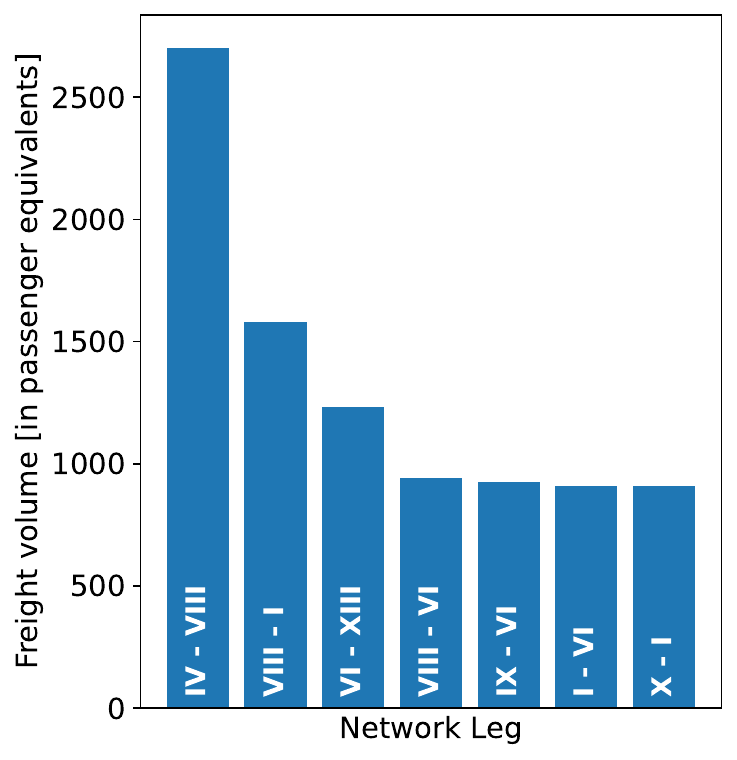}
        \captionsetup{width=0.9\columnwidth}
        \caption{\centering Spatial utilization}
        \label{fig:spatial-utilization}
    \end{subfigure}
    \begin{subfigure}{0.99\textwidth}
        \centering
        \includegraphics[width=0.3\linewidth]{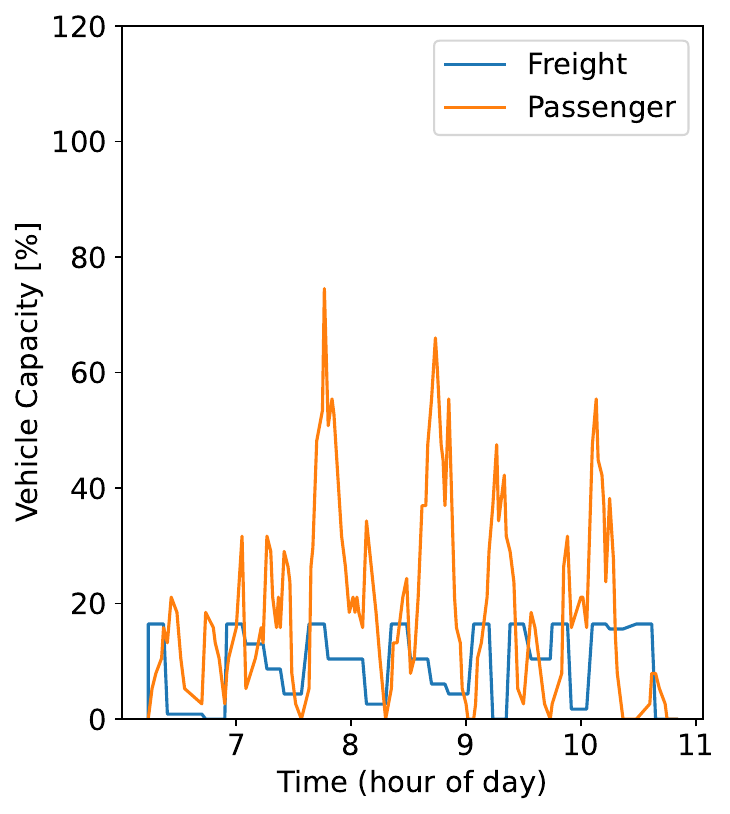}
        \includegraphics[width=0.3\linewidth]{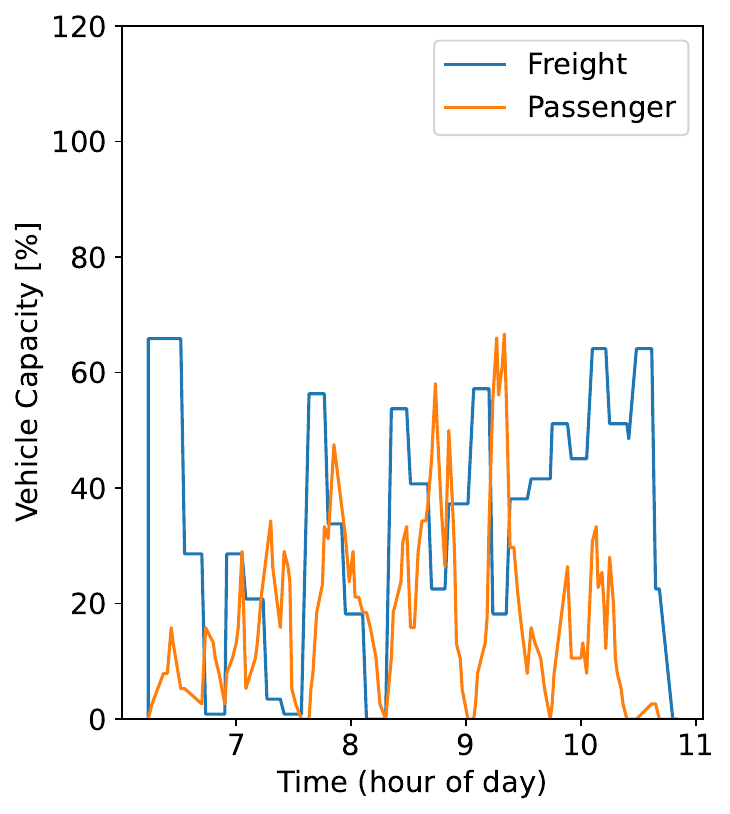}
        \includegraphics[width=0.3\linewidth]{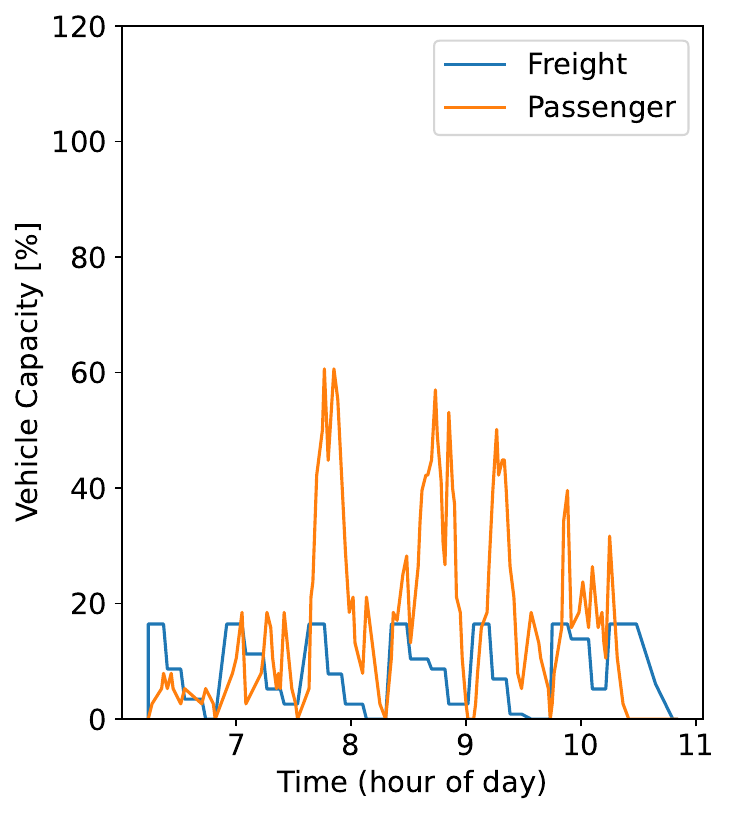}
        \captionsetup{width=0.9\columnwidth}
        \caption{\centering Vehicle utilization}
        \label{fig:vehicle-utilization}
    \end{subfigure}
    \caption{System utilization ($n=1$)}
    \label{fig:utilization}
    \vspace{0.3em}
    \begin{minipage}{0.95\textwidth}
        \centering \justify \scriptsize
        I:~Sendlinger Tor, II:~Scheidplatz, III:~Garching, Forschungszentrum, IV:~Neuperlach Zentrum, V:~Klinikum Großhadern, VI:~Odeonsplatz, VII:~Mangfallplatz, VIII:~Innsbrucker Ring, IX:~Arabellapark, X:~Olympia-Einkaufszentrum, XI:~Münchner Freiheit, XII:~Messestadt Ost, XIII:~Westendstraße, XIV:~Fürstenried West, XV:~Feldmoching
    \end{minipage}
    \vspace{-0.2cm}
\end{figure}

In the following, we analyze system characteristics for these three cases in Figure~\ref{fig:utilization} before providing a more granular analysis on the respective cost trade-off in Table~\ref{tab:sensitivity}. Figure~\ref{fig:utilization} details the system utilization for the three scenarios mentioned above, focusing on the system's utilization over time (Figure~\ref{fig:temp-utilization}), the system's spatial utilization (Figure~\ref{fig:spatial-utilization}), as well as the utilization per vehicle over time for one representative vehicle (Figure~\ref{fig:vehicle-utilization}). 

Focusing on the system's utilization over time, Figure~\ref{fig:temp-utilization} shows the share of freight requests and passengers in the system with respect to the overall amount of requests and passengers in the respective scenario. To visualize the systems dynamics accurately, we exclude idle requests from this visualization. Since passenger flows are prioritized over freight requests and are constrained in time and alternative paths, the passenger utilization shows a similar pattern in all three scenarios, exhibiting a typical commuter peak around 8 a.m. For the system's freight utilization, we observe different dynamics across all scenarios: in the base case, the freight utilization ranges constantly between 4\% and 8\% during most of the time horizon, showing slightly elevated utilizations before and after the passenger utilization peak. In the optimistic scenario, we try to push as much freight as possible through the system but are limited by its capacity constraints, i.e., the additional passenger flow that occupies transport volume and is prioritized over the freight requests. Accordingly, we observe a drop in freight utilization at the beginning of the passenger peak, and an additional peak in freight utilizaiton once the passenger peak declines. In the pessimistic scenario, we observe a significantly reduced freight utilization that stems from the shifted cost ratio; although additional transport capacity is available in the system it is cost-optimal to leave the freight requests to conventional truck-based delivery.

Figures~\ref{fig:spatial-utilization}~and~\ref{fig:vehicle-utilization} detail the impact of these flow volumes for the most utilized legs in the \gls{acr:pts} and a representative vehicle. As can be seen, the utilization related to freight requests in certain network legs scales with the overall utilization in the system (see Figure~\ref{fig:spatial-utilization}). Still, we observe bottleneck effects as the utilization on Leg~IV-VIII is significantly higher than all other utilizations across all scenarios.

\revone{To further assess the effects on Leg~IV--VIII, Figure~\ref{fig:spatial-utilization-passenger} displays the $20$ most utilized \gls{acr:pt} legs together with the corresponding freight volumes on the respective legs. The figure also reports the number of services on each leg, where a service occurs every time a vehicle traverses the leg. The spatial distribution of passenger flow on the network is stable across the different scenarios and is linked to the available number of services. 

\begin{figure}[!b]
    \centering
    \includegraphics[height=0.37\linewidth]{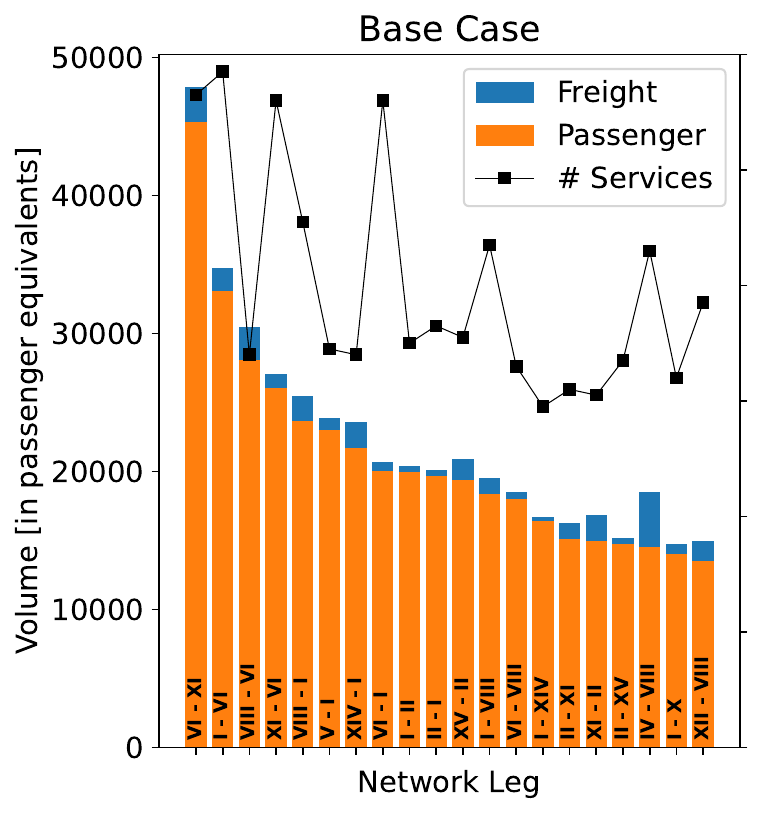}
    \includegraphics[height=0.37\linewidth]{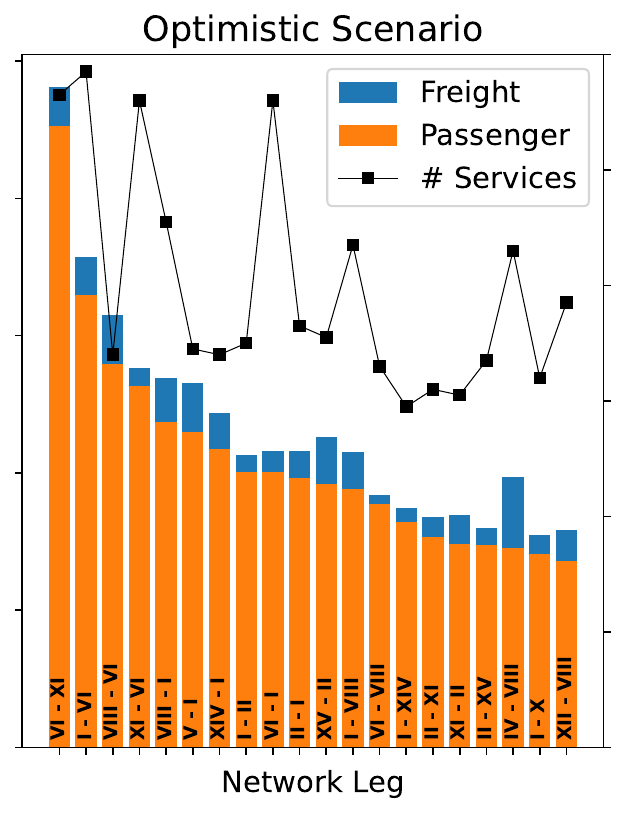}
    \includegraphics[height=0.37\linewidth]{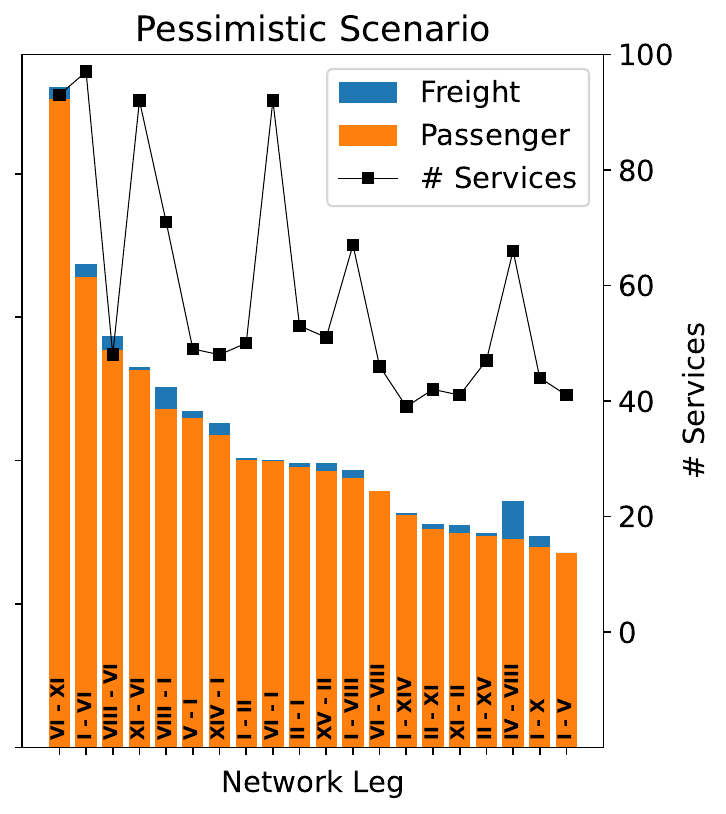}
    \captionsetup{width=0.9\columnwidth}
    \caption{\centering Spatial utilization by passenger flows ($n=1)$}
    \label{fig:spatial-utilization-passenger}
\end{figure}

This observation underlines that passenger demand largely takes precedence over freight requests and that the number of freight shipments exerts only a very limited influence on passenger service quality in a spatial dimension. Note that we also restrict the impact on passenger service quality in a temporal dimension by selecting the potential paths accordingly.
Furthermore, Figure~\ref{fig:spatial-utilization-passenger} shows that freight volume remains small relative to passenger volume in the system, even in the optimistic scenario. Moreover, Leg~IV--VIII carries fewer than $20{,}000$ passengers despite $66$ services on this leg. This negative correlation between passenger and freight volumes persists across the top~$4$ legs for freight transportation in the optimistic scenario. Leg~VIII--VI carries the highest freight volume while also ranking among the top~$20$ legs by passenger volume. Despite this demand, only $48$ services operate on this leg. Hence, Leg~VIII--VI plays a critical role in the network. Figure~\ref{fig:munich-network} supports this observation by showing that this leg is serviced only by a single line (i.e., U5) while being a central component and connecting districts east of the river with the city center.}

\begin{figure}[!b]
    \centering
    \includegraphics[width=0.3\linewidth]{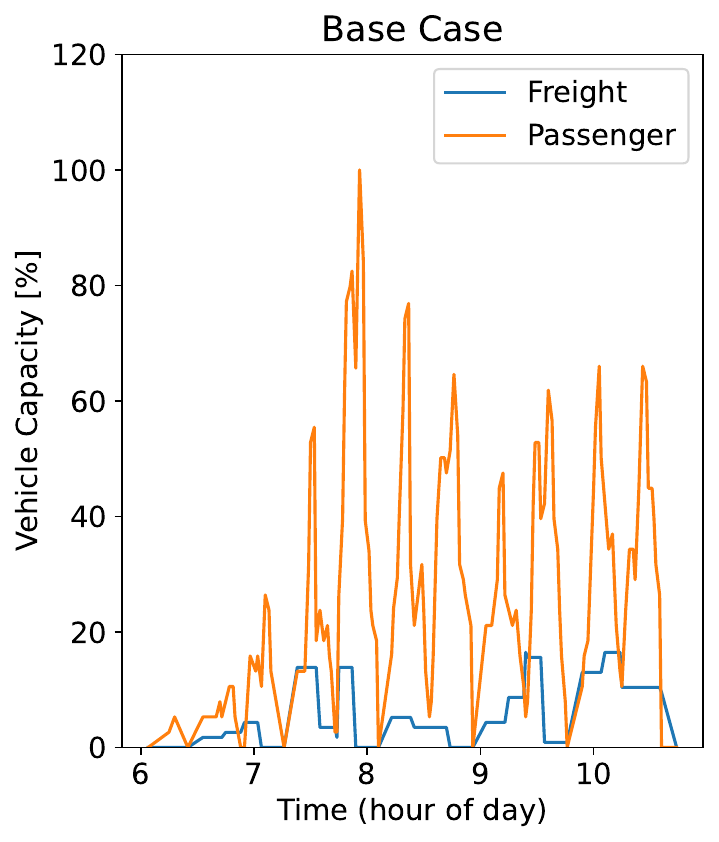}
    \includegraphics[width=0.3\linewidth]{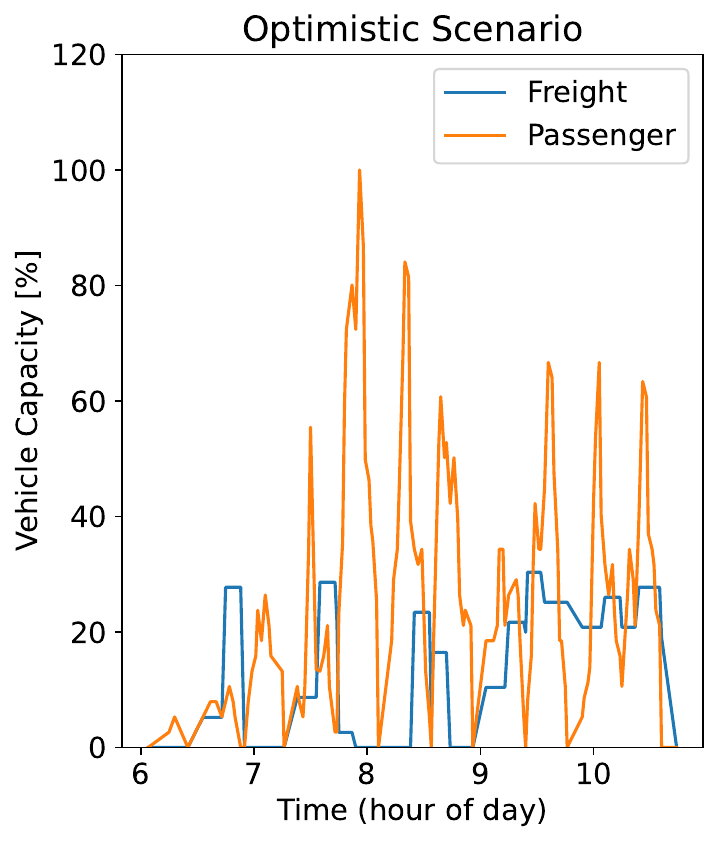}
    \includegraphics[width=0.3\linewidth]{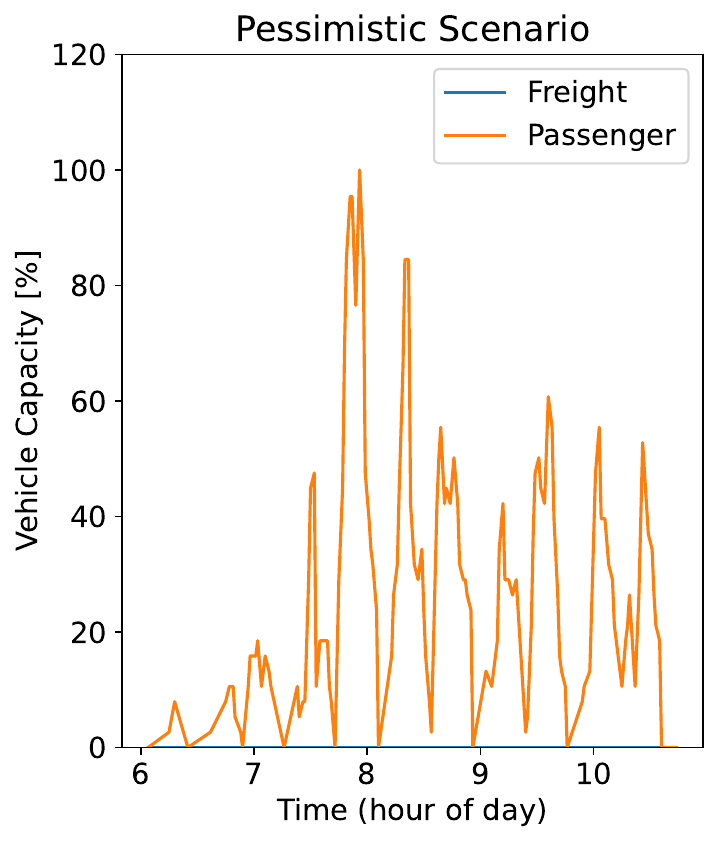}
    \captionsetup{width=0.9\columnwidth}
    \caption{\centering Temporal utilization of additional vehicle ($n=1$)}
    \label{fig:add_vehicle}
\end{figure}

From a microscopic perspective at vehicle level, we observe that the passenger utilization curve shows multiple peaks in a frequency of around $0.5$ hours (see Figure~\ref{fig:vehicle-utilization}). These peaks stem from the vehicles trajectory, going back and forth on its line: the vehicle frequently crosses the vaster city center area but is less utilized around the turnaround points. As can be seen, our algorithm determines a solution that routes freight through the network utilizing the available excess capacity of the vehicle. \revone{Figure~\ref{fig:vehicle-utilization} shows the utilization of one of the $85$ vehicles that operate on the \gls{acr:pt} network. This particular vehicle is heavily utilized for freight --- especially in the optimistic scenario. In this scenario, the passenger utilization remains moderate, so that even during peak hours only $60$\% of the vehicle capacity is occupied, and the algorithm leverages the excess capacity by allocating a substantial number of freight shipments to this vehicle. In particular during the first hour of operations and from $9{:}30$ a.m. to $11$ a.m. the vehicle predominantly transports freight. Comparing the optimistic scenario in Figure~\ref{fig:vehicle-utilization} to both the base case and the pessimistic scenario shows that the algorithm shifted passenger flow from the displayed vehicle into other vehicles of the fleet to extend the available capacity for freight. Moreover, the transported freight volumes in the base case and the pessimistic scenario hardly differ. Figure~\ref{fig:add_vehicle} shows a complementary vehicle of the \gls{acr:pt} fleet during the same operations in which the difference in total transported freight volume between the scenarios is also reflected on the vehicle level. Even more drastically, in the pessimistic scenario, the vehicle is not part of the cargo-hitching system and does not carry freight at all. However, the vehicle displays low transported freight volumes in the base case, and the transport volume increases even further in the optimistic scenario.} The beneficial effect of dynamically allocating capacity is particularly evident in \revone{Figure~\ref{fig:add_vehicle} at around $8$ a.m. of the optimistic scenario} when passengers almost fully utilize the vehicle, and no additional freight is transported. Nevertheless, when its capacity is not fully required for passengers, it repeatedly transports freight.

\vspace{2pt}
\begin{result}
Cargo-hitching offers a utilization increase at zero additional installed capacity, and our algorithmic framework provides solutions that predominantly utilize the \gls{acr:pts}'s off-peak hours to transport freight requests, leading to higher overall utilizations. 
\end{result}
\vspace{2pt}

Abstracting from the scenarios analyzed in Figure~\ref{fig:utilization}, it becomes evident that the potential savings that can be realized by cargo-hitching depend on the spare capacity left within the \gls{acr:pts} as well as the amount of freight requests that can be shipped through it. In this context, relying on \glspl{acr:htu} to realize cargo-hitching is particularly beneficial if the amount of freight requests is high but the spare capacity left in the \gls{acr:pts} fluctuates over the day due to passenger flow peaks. 

Beyond these general findings, the amount of freight requests transported through the \gls{acr:pts} is sensitive to the loading and unloading cost as well as the externality cost of truck-based deliveries. Table~\ref{tab:sensitivity} shows this trade-off and its impact on the share of rejected freight requests.

At an externality cost of truck-based delivery per vehicle and kilometer of \revone{\euro{0.1} the acceptance of cargo-hitching vanishes and all requests are rejected by the municipality. The lower bound reported in the literature for this parameter is \euro{0.05} \citep[cf.][]{DeLanghe2017}}. On the other hand, rejection rates diminish at higher penalty costs driven by higher externality costs from truck-based deliveries. As can be seen, \revone{at penalty costs resulting from an externality cost of truck-based transport of \euro{1.5} per vehicle and kilometer, rejection rates fall below a small threshold of $2$\% and the municipality accepts almost all requests}. 

\vspace{2pt}
\begin{result}
    At an externality cost of \euro{1.5} per vehicle and kilometer, cargo-hitching reaches full penetration if the cost for loading and unloading is less than \euro{2.0} per passenger equivalent unit.
\end{result}
\vspace{2pt}

Focusing on a broader interpretation of Table~\ref{tab:sensitivity}, we observe that rejection rates above the diagonal are low. This general observation indicates that the penetration of the concept is highly correlated with the direct comparison between the externality cost for truck-based delivery and the cost for loading and unloading represented by the respective arc cost factor. 

\revone{
To quantify the value of allowing a dynamic allocation of capacity, we focus on the large instances with $3{,}000$ requests and determine solutions with a static capacity allocation by modifying Constraints~\eqref{cons:cap-assignment} in the \gls{acr:rmp} of the \gls{acr:pab} algorithm such that they enforce equality between their left-hand sides and their right-hand sides. After deriving the integer feasible solutions with a static capacity allocation, we relax Constraints~\eqref{cons:cap-assignment} to their original version, set the solution with the static allocation as a primal bound, and branch again to derive integer feasible solutions with a dynamic capacity allocation (cf. Line~\ref{l:pricing-over-pab} of Algorithm~\ref{algo:price-and-branch}). 

\begin{table}[!b]
\centering
\footnotesize
\caption{Average share of rejected requests depending on externality cost ($n=15$, $|\SetOfRequests|=3{,}000$)}
\label{tab:my-table}
\begin{tabular}{l*{9}{S[table-format=1.3]}}
\toprule
\multirow{2}{*}{\begin{tabular}[c]{@{}c@{}}$\Cost_{\UnbracedArc}$, \\ $\Arc \in \SetOfArcs^{\Transit}$\end{tabular}} &
\multicolumn{9}{c}{Externality cost (Truck) [EUR per vehicle and kilometer]} \\
\cmidrule(lr){2-10}
& {0.10} & {0.20} & {0.30} & {0.40} & {0.50} & {0.75} & {1.00} & {1.25} & {1.50} \\
\midrule
\midrule
0.1 & 1.000 & 0.623 & 0.009 & 0.001 & 0.000 & 0.001 & 0.000 & 0.000 & 0.014 \\
0.2 & 1.000 & 0.841 & 0.057 & 0.002 & 0.000 & 0.000 & 0.000 & 0.000 & 0.000 \\
0.3 & 1.000 & 0.950 & \textbf{0.427} & 0.004 & 0.001 & 0.000 & 0.000 & 0.000 & 0.000 \\
0.4 & 1.000 & 1.000 & 0.674 & 0.073 & 0.002 & 0.000 & 0.001 & 0.000 & 0.000 \\
0.5 & 1.000 & 1.000 & 0.717 & 0.210 & 0.013 & 0.000 & 0.001 & 0.000 & 0.000 \\
0.6 & 1.000 & 1.000 & 0.767 & 0.607 & 0.068 & 0.000 & 0.000 & 0.000 & 0.000 \\
0.8 & 1.000 & 1.000 & 0.970 & 0.695 & 0.480 & 0.003 & 0.000 & 0.000 & 0.000 \\
1.0 & 1.000 & 1.000 & 1.000 & 0.729 & 0.683 & 0.060 & 0.003 & 0.001 & 0.000 \\
1.2 & 1.000 & 1.000 & 1.000 & 0.877 & 0.688 & 0.065 & 0.001 & 0.000 & 0.000 \\
1.4 & 1.000 & 1.000 & 1.000 & 1.000 & 0.703 & 0.467 & 0.058 & 0.001 & 0.000 \\
1.6 & 1.000 & 1.000 & 1.000 & 1.000 & 0.790 & 0.679 & 0.059 & 0.001 & 0.001 \\
1.8 & 1.000 & 1.000 & 1.000 & 1.000 & 1.000 & 0.679 & 0.067 & 0.058 & 0.000 \\
2.0 & 1.000 & 1.000 & 1.000 & 1.000 & 1.000 & 0.681 & 0.475 & 0.058 & 0.001 \\
\bottomrule
\end{tabular}
\label{tab:sensitivity}
\end{table}

Moreover, we expect the value of dynamic capacity allocation to increase if the available excess capacity in the \gls{acr:pts} shrinks. Therefore, we study different settings in which we consider elevated passenger demand by scaling $\Demand_{\Request} \leftarrow (1 + \PassDemandScale) \Demand_{\Request}, \medspace \Request \in \SetOfRequests^{\Passenger}$ with a scaling factor $\PassDemandScale \in \{0.0, 0.05, 0.1, 0.15\}$ where $\PassDemandScale = 0.0$ is the default parameterization. Note that for $\PassDemandScale \geq 0.2$, the total system capacity is insufficient and the problem is rendered infeasible by Constraint~\eqref{cons:service-level}, which yields a natural bound for the value range of $\PassDemandScale$.

Figure~\ref{fig:static-dynamic-costs} shows the value of allowing the dynamic allocation of capacity in terms of costs for varying $\PassDemandScale$ in the base case and in the optimistic scenario respectively. In the base case (cf. Figure~\ref{fig:static-benchmark-costs-base}), the impact on the total cost is marginal across all scales of passenger demand and only occasionally the algorithm finds marginally improved solutions with different rejection rates which leads to different cost component distributions. However, Figure~\ref{fig:static-benchmark-costs-opt} shows that in the optimistic scenario with elevated freight volume in the system, the picture changes. Allowing dynamic capacity allocation leads to increased design costs, decreased penalty costs, and accordingly increased routing costs. While the cost effects almost balance out for $\PassDemandScale=0.0$, the median total cost saving in the optimistic scenario with $\PassDemandScale=0.15$ is $1.3$\% and reaches up to $3.2$\%. The algorithm accomplishes this by accepting up to $100$\% of the requests that had to be rejected in the solution with static capacity allocation due to capacity bottlenecks.      
}

\begin{figure}[!t]
  \begin{subfigure}[t]{0.49\textwidth}
    \centering
    \input{Figure_14a}
    \vspace{-0.3cm}
    \centering \caption{Base Case}
    \label{fig:static-benchmark-costs-base}
  \end{subfigure}\hfill
  \begin{subfigure}[t]{0.49\textwidth}
    \input{Figure_14b}
    \vspace{-0.3cm}
    \centering \caption{Optimistic Scenario}
    \label{fig:static-benchmark-costs-opt}
  \end{subfigure}
  \caption{Relative cost difference between a static and a dynamic capacity allocation for varying $\PassDemandScale$ and different cost components. Negative values refer to a cost decrease in the dynamic allocation case compared to the static allocation case ($n=15$, $|\SetOfRequests|=3{,}000$)}
  \label{fig:static-dynamic-costs}
\end{figure}

\vspace{2pt}
\revone{
\begin{result}
    Allowing a dynamic allocation of capacity increases the acceptance rates of cargo-hitching and decreases total costs. The increase of acceptance rates and the decrease of total costs are negatively correlated to the available excess capacity. In the optimistic scenario, dynamic capacity allocation yields up to $100$\% increased acceptance rates and cost savings of up to $3.2$\%.  
\end{result}
}

\revtwo{As shown in Figure~\ref{fig:static-dynamic-costs}, total cost savings from dynamic capacity allocation are modest in the base settings, not exceeding $3.2$\%. This occurs because the economic gain from higher acceptance rates (utilizing off-peak capacity) is nearly offset by the increased operational costs required to route these additional shipments.

To identify conditions where dynamic allocation is more advantageous, we explore scenarios where successful service fulfillment is incentivised. We introduce a scaling factor $\PenaltyCostScale$ to increase the penalty costs for rejected requests: $\Cost^{\Request}_{\Penalty} \leftarrow (1 + \PenaltyCostScale) \Cost^{\Request}_{\Penalty}, \medspace \Request \in \SetOfRequests^{\Freight}$. This represents scenarios with stricter service mandates or regulatory penalties for conventional transport.

Figure~\ref{fig:sensitivity-penalty} illustrates the cost difference between static and dynamic allocation as $\PenaltyCostScale$ increases. In the base case, the difference is negligible. However, as the cost of rejection rises, the flexibility of dynamic allocation becomes critical. At $\PenaltyCostScale = 1.0$ (a $100$\% increase in penalty), median savings rise to $6$\% in the base scenario. In the optimistic scenario, this effect is even stronger, reaching median savings of $12$\%. This confirms that dynamic allocation provides significant value when service reliability and rejection avoidance are prioritized.
}

\vspace{2pt}
\revtwo{
    \begin{result}
         The comparative value of dynamic capacity allocation scales with the cost of service failure. While base savings are modest, dynamic allocation significantly outperforms static assignment when the penalty for rejecting freight is high. Specifically, doubling the rejection penalty increases cost savings to $12$\% in optimistic scenarios.
    \end{result}
}

\begin{figure}[!t]
    \centering
    \input{Figure_15}
    \caption{Relative cost difference between a static and a dynamic capacity allocation for varying $\PenaltyCostScale$ across different scenarios. Negative values refer to a cost decrease in the dynamic allocation case compared to the static allocation case ($n=15$, $|\SetOfRequests|=3{,}000$, $\PassDemandScale = 0.15$)}
    \label{fig:sensitivity-penalty}
\end{figure}

\section{Conclusion} \label{sec:conclusion}
We introduced the urban cargo-hitching problem with dynamic allocation of capacity, heterogeneous \gls{acr:pt} vehicles, and freight transshipments on a state-of-the-art partially time-expanded and spatially-expanded graph. Based on the expanded graph, we provided an algorithmic framework that relies on multiple preprocessing techniques and a \gls{acr:pab} algorithm to solve instances with up to $3{,}000$ freight requests, obtaining a median integrality gap of less than $1.56$\% within computational time of $90$ minutes. We further present a \gls{acr:bap} algorithm that allows to obtain even smaller optimality gaps of up to $0.38$ percentage points at the price of increased \revone{pricing effort}.

Our results for the subway network of Munich, Germany, indicate that the ratio of externality costs is the determining factor for high penetration rates of the cargo-hitching concept. We conducted a sensitivity analysis on this ratio and found that cargo-hitching is worthwhile if truck-based transport occurs at an externality cost of more than \euro{1.5} per vehicle and kilometer and loading and unloading costs of less than \euro{2.0} per passenger equivalent. Additionally, we showed that our framework increases the \gls{acr:pts} utilization with a focus on off-peak hours and enables decision-makers to evaluate the importance of single parts of the evaluated \gls{acr:pts}. We show that relying on \glspl{acr:htu} to realize cargo-hitching is particularly beneficial if the amount of freight requests is high but the spare capacity left in the \gls{acr:pts} \revone{is low and} fluctuates over the day due to passenger flow peaks. Moreover, we quantify the value of relying on \glspl{acr:htu} in terms of total cost to reach up to $3.2$\%, and even at lower total cost savings, \glspl{acr:htu} support the acceptance of cargo-hitching. \revtwo{Assuming stricter service mandates or regulatory penalties for conventional transport, the median savings due to \glspl{acr:htu} increase to $12$\% in optimistic scenarios.}

Our work provides a scalable algorithmic framework that lays the foundation for future work, e.g., by extending it to determine (capacitated) \gls{acr:ft} locations. In this context, future work may also focus on incorporating stochastic demand patterns to take informed strategic decisions on the respective network design. 

\section*{Acknowledgments}
This research did not receive any specific grant from funding agencies in the public, commercial, or not-for-profit sectors.

\singlespacing{
    \setstretch{0.92}
    
} 

% start appendix
\newpage
\onehalfspacing
\appendix
\section{List of abbreviations}\label{app:abbrevations}

% remove page counter in hard coded glossary
\makeatletter
\renewcommand*\glossaryentrynumbers[1]{}%
\makeatother
  
\begin{theglossary}\glossaryheader
\glsgroupheading{A}\relax \glsresetentrylist %
\glossentry{acr:agv}{}%
\glossentry{acr:alns}{}\glsgroupskip
\glsgroupheading{B}\relax \glsresetentrylist %
\glossentry{acr:bab}{}%
\glossentry{acr:bap}{}\glsgroupskip
\glsgroupheading{C}\relax \glsresetentrylist %
\glossentry{acr:cg}{}\glsgroupskip
\glsgroupheading{F}\relax \glsresetentrylist %
\glossentry{acr:ft}{}\glsgroupskip
\glsgroupheading{H}\relax \glsresetentrylist %
\glossentry{acr:hlp}{}%
\glossentry{acr:htu}{}\glsgroupskip
\glsgroupheading{L}\relax \glsresetentrylist %
\glossentry{acr:lsp}{}\glsgroupskip
\glsgroupheading{M}\relax \glsresetentrylist %
\glossentry{acr:mip}{}\glsgroupskip
\glsgroupheading{P}\relax \glsresetentrylist %
\glossentry{acr:pab}{}%
\glossentry{acr:pt}{}%
\glossentry{acr:pts}{}\glsgroupskip
\glsgroupheading{R}\relax \glsresetentrylist %
\glossentry{acr:rmp}{}\glsgroupskip
\glsgroupheading{S}\relax \glsresetentrylist %
\glossentry{acr:spp}{}\glsgroupskip
\glsgroupheading{V}\relax \glsresetentrylist %
\glossentry{acr:vrp}{}%
\end{theglossary}\glossarypostamble
\pagebreak

\section{Notation}\label{app:notation}
Table~\ref{tbl:problem-notation} provides a summary of notation which is introduced in the main body of the paper.
\begin{table}[!h]
	\begin{threeparttable}
		\caption{Notation}
		\label{tbl:problem-notation}
		\centering
        \footnotesize 
		\begin{tabular}{lr}
			\toprule
			Symbol & Meaning \\
			\midrule
			\multicolumn{2}{c}{Basic sets} \\
            \midrule
            $\Graph = (\SetOfVertices, \SetOfArcs)$ & Expanded, multi-layered, preprocessed and directed graph \\
            $\SetOfPaths$ & Set of valid paths in $\Graph$ \\
			$\SetOfVehicles$ & Set of vehicles \\
            $\SetOfRequests$ & Set of requests \\
            \midrule
			\multicolumn{2}{c}{Subsets} \\
            \midrule
            $\SetOfRequests^{\Passenger} \subseteq \SetOfRequests$ & Set of passenger requests \\
            $\SetOfRequests^{\Freight} \subseteq \SetOfRequests$ & Set of freight requests \\
			$\SetOfStops \subset \SetOfVertices^{\Exp}$ & Temporal vertices representing \gls{acr:pt} stops \\
			$\SetOfTerminals \subseteq \SetOfStops^{\Exp}$ & Temporal vertices representing \glspl{acr:ft} \\
			$\SetOfOrigins^{\Exp} \subset \SetOfVertices^{\Exp}$ & Temporal vertices representing request origins\\
			$\SetOfDestinations^{\Exp} \subset \SetOfVertices^{\Exp}$ & Temporal vertices representing request destinations\\
            $\SetOfArcs^{\VehicleLayer} \subset \SetOfArcs$ & Vehicle arcs ($\Cost_{\UnbracedArc} > 0, \quad \forall \Arc \in \SetOfArcs^{\VehicleLayer}$) \\
            $\SetOfArcs^{\PathSegment} \subset \SetOfArcs$ & Freight segments arcs ($\Cost_{\UnbracedArc} > 0, \quad \forall \Arc \in \SetOfArcs^{\PathSegment}$) \\
            $\SetOfArcs^{\Hold} \subset \SetOfArcs$ & Holding arcs ($\Cost_{\UnbracedArc} = 0, \quad \forall \Arc \in \SetOfArcs^{\Hold}$) \\
            $\SetOfArcs^{\Transit} \subset \SetOfArcs$ & Transit arcs ($\Cost_{\UnbracedArc} > 0, \quad \forall \Arc \in \SetOfArcs^{\Transit}$) \\
            $\SetOfArcs^{\Dummy} \subset \SetOfArcs$ & Dummy arcs ($\Cost_{\UnbracedArc} > 0, \quad \forall \Arc \in \SetOfArcs^{\Dummy}$) \\
            $\SetOfArcs^{\Access} \subset \SetOfArcs$ & Access arcs ($\Cost_{\UnbracedArc} \ge 0, \quad \forall \Arc \in \SetOfArcs^{\Access}$) \\
            $\SetOfArcs^{\Egress} \subset \SetOfArcs$ & Egress arcs ($\Cost_{\UnbracedArc} \ge 0, \quad \forall \Arc \in \SetOfArcs^{\Egress}$) \\
            \midrule
			\multicolumn{2}{c}{Indices} \\
            \midrule
			$\Vehicle$ & Vehicle $\Vehicle \in \SetOfVehicles$ \\
			$\Arc$ & Arc in $\Graph=(\SetOfVertices, \SetOfArcs)$ \\
			$\Request$ & Request $\Request \in \SetOfRequests$ \\
			$\Path$ & Path $\Path \in \SetOfPaths$ \\
			\midrule
			\multicolumn{2}{c}{Variables} \\
            \midrule
			$\FreightFlow^{\Request}_{\UnbracedArc}$ & Freight flow of $\Request$ traversing $\Arc$ \\
			$\PassengerFlow^{\Request}_{\Path}$ & Passenger flow of $\Request$ traversing $\Path$ \\
			$\HTUVehicle_{\Vehicle}$ & Number of \glspl{acr:htu} assigned to $\Vehicle$ \\
			$\HTUOps_{\UnbracedArc}$ & Number of \glspl{acr:htu} transporting freight on $\Arc$ \\ 
			\midrule
			\multicolumn{2}{c}{Parameters} \\
            \midrule
			$\ServiceLevel$ & Demand-weighted passenger service level \\
			$\Cost_{\Vehicle}$ & Scaled investment cost per \gls{acr:htu} with vehicle $\Vehicle$ \\
			$\Cost_{\UnbracedArc}$ & Cost per unit of flow on arc $\Arc$ \\	
			$\UnitCapacity_{\Vehicle}$ & Transportation unit capacity of vehicle $\Vehicle$ \\
			$\NumUnits_{\Vehicle}$ & Number of transportation units with vehicle $\Vehicle$ \\
			\midrule
			\multicolumn{2}{c}{Other} \\
            \midrule
            $\SetOfPaths(\Request) \subseteq \SetOfPaths$ & Paths of $\Request$ \\	
            $\Neighbors^{+}(i)$ / $\Neighbors^{-}(i) \subset \SetOfVertices^{\Exp}$ & Set of neighboring vertices of $i$ (w.r.t. freight) \\
            $\SetOfArcs^{\PotentialFreightArcs} \subseteq \SetOfArcs$ & Temporal arcs allowing freight flow \\
            $\PathArcRelation^{\Path}_{\UnbracedArc}$ & $1$ if path $\Path$ contains $\Arc$, $0$ otherwise \\
            $\VehicleArcRelation^{\Vehicle}_{\UnbracedArc}$ & $1$ if vehicle $\Vehicle$ operates arc $\Arc$, $0$ otherwise \\
            $\VertexDemand^{\Request}_{i}$ & Temporal vertex demand at $i$ for request $\Request$ \\
			\bottomrule
		\end{tabular}
	\end{threeparttable}
\end{table}
\pagebreak

\section{Partial pricing results} \label{app:partial-pricing}
Figure \ref{fig:partial-pricing} shows the integrality gaps after $90$ minutes of our \gls{acr:pab} approach for varying instance sizes and pricing strengths.
\begin{figure}[htbp]
    \centering
    \input{Figure_16.tex}
    \caption{Integrality gaps for varying pricing strengths ($n=15$)}
    \label{fig:partial-pricing}
\end{figure}
\FloatBarrier
\newpage

\section{Exemplary passenger paths}\label{app:passenger-paths}
\begin{figure}[hbt]
  \centering
  \setlength{\abovecaptionskip}{0pt} 
  \setlength{\belowcaptionskip}{0pt}  
  \setlength{\textfloatsep}{0pt}    
  \begin{subfigure}{0.5\textwidth}
    \centering
    \includegraphics[width=\linewidth]{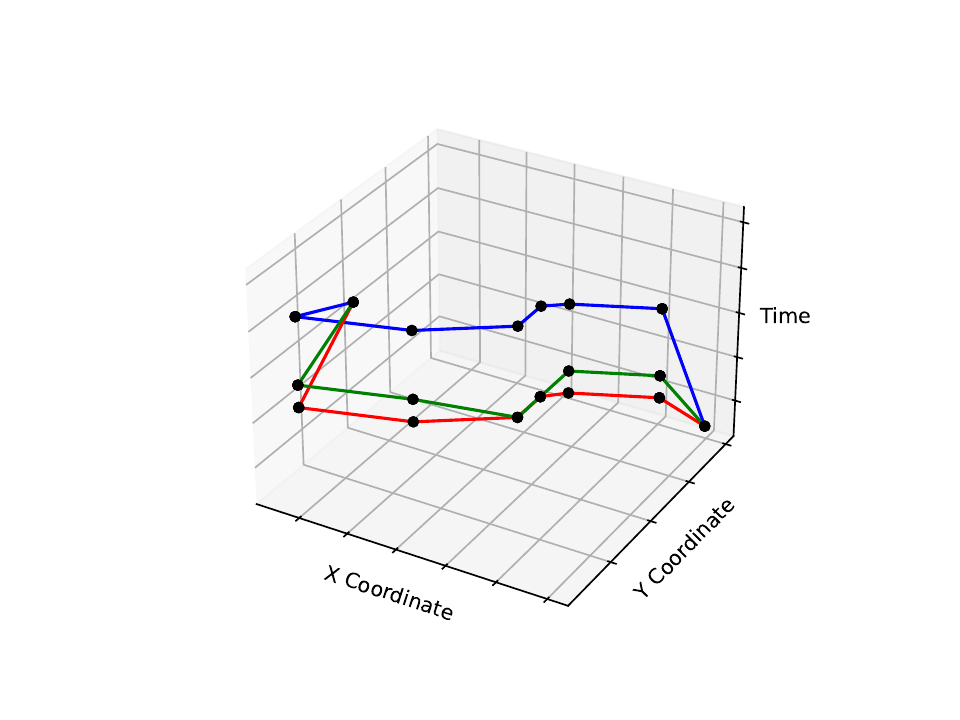}
    \caption{Non-disjoint paths}
    \label{fig:non-disjunct-paths}
  \end{subfigure}%
  \begin{subfigure}{0.5\textwidth}
    \centering
    \includegraphics[width=\linewidth]{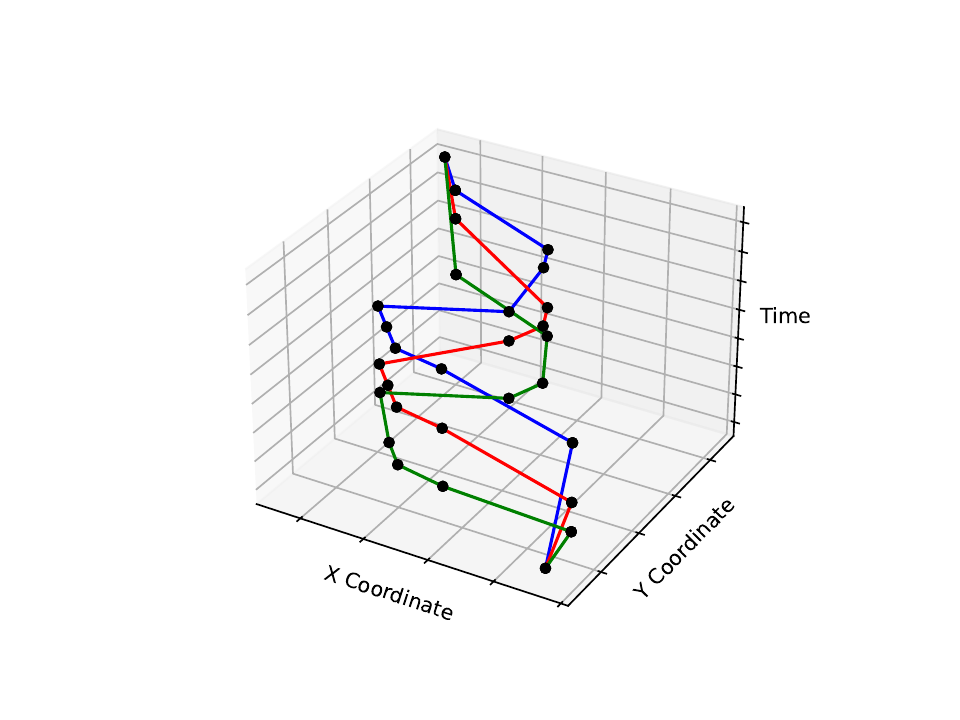}
    \caption{Disjoint paths}
    \label{fig:disjunct-paths}
  \end{subfigure}
  \caption{Exemplary passenger paths}
  \label{fig:app-pass-paths}
\end{figure}

We assume tight service time intervals and argue that passengers prefer convenience, i.e., short travel times. Thus, we compute the shortest paths by travel time for every request. In this context, we purposefully do not compute edge-disjoint paths as this appears to be a hard limitation, especially for requests that start or end at remote locations only connected to the \gls{acr:pts} by one single line. In these cases, the computation of disjoint paths would lead to a set of paths that mainly differ in the temporal dimension, and the resulting paths would not offer a choice between equally valued paths from a passenger's perspective. Furthermore, the network topology of a subway network leads to paths with limited overlap and, in many instances, even to disjoint paths without explicitly enforcing it. Figure~\ref{fig:non-disjunct-paths} demonstrates why we waive the computation of disjoint paths. In this example, two paths are almost edge-disjoint but share a single edge. Figure~\ref{fig:disjunct-paths} shows an example where our approach yields edge-disjoint paths.  

\newpage

\section{Implementation}\label{app:implementation}
Two comments on the graph expansion and algorithmic framework implementation are in order. 

First, we pre-compute feasible and realistic paths for passenger requests $\Request \in \SetOfRequests^{\Passenger}$ by computing the three shortest paths with respect to travel time such that the itinerary starts and ends in the given service time interval $[\RequestStart^{\Request}, \RequestEnd^{\Request}]$. Here, we assume a reasonable walking speed of $1$ m/s on access and egress arcs $\Arc \in \SetOfArcs^{\Access} \cup \SetOfArcs^{\Egress}$. 

Second, we prune the graph $\Graph$ when creating the sets $\SetOfArcs^{\Access}$ and $\SetOfArcs^{\Egress}$ as follows. We remove $(\Origin^{\Request}, \RequestStart^{\Request}) \text{ and } (\Destination^{\Request}, \RequestEnd^{\Request}), \Request \in \SetOfRequests^{\Passenger}$ including their incoming and outgoing arcs from graph $\Graph$ because our algorithmic framework does not requires them anymore after pre-computing the paths for passenger requests. Instead, we remove the respective components from both the graph $\Graph$ and the pre-computed paths. 
Furthermore, we only connect $i = (\Origin^{\Request}, \RequestStart^{\Request}) \in \SetOfOrigins, \Request \in \SetOfRequests^{\Freight}$ to a vertex $j = (\Stop, \TimeStep, \Hold) \in \SetOfStops^{\Hold}$ if 
\begin{itemize}
    \item[a)] the time difference remains above a certain threshold $\zeta(\Request), \Request \in \SetOfRequests^{\Freight}$ indicating that time suffices to relocate from $\Origin^{\Request}$ to $\Stop$, thus $\TimeStep - \RequestStart^{\Request} \ge \zeta(\Request)$
    \item[b)] if there is no earlier representation of the same stop $\Stop$ to which $i$ can be connected while respecting condition a), thus $\nexists \medspace (\Stop, \TimeStep^{\prime}, \Hold) \in \SetOfStops^{\Hold}: \TimeStep^{\prime} < \TimeStep \land \TimeStep^{\prime} - \RequestStart \ge \zeta(\Request)$
    \item[c)] the vertex $j$ represents one of the $\NumConnections$ closest \glspl{acr:ft} in the \gls{acr:pts} according to the distance between the represented stop $\Stop$ and the requests origin $\Origin^{\Request}$.
\end{itemize}
By applying the same reasoning, we prune the temporal arc set $\SetOfArcs^{\Egress}$. Hence, we condition the existence of an arc on 
\begin{itemize}
    \item[a)] $\RequestEnd^{\Request} - \TimeStep \ge \zeta(\Request)$
    \item[b)] the absence of a later representation of the same stop in the holding layer
    \item[c)] $\Stop$ being one of the $\NumConnections$ closest \glspl{acr:ft}. 
\end{itemize} 

We allow to temporarily store freight at \glspl{acr:ft} by assigning a zero cost to the arcs in $\SetOfArcs^{\Hold}$. In this context, connecting origins and destinations to the earliest and latest representations of stops reduces the cardinality of the arc set without sacrificing solution quality. We can post-process the solutions to avoid unnecessary long service times due to holding a request at the first or last stop of its path. 

\end{document}